\documentclass[11pt]{amsart}
\usepackage{amsmath}
\usepackage{amssymb}
\usepackage{mathtools}
\usepackage{bm}
\usepackage{esint}
\usepackage{graphicx}
\usepackage{verbatim}
\usepackage{enumitem}
\usepackage[normalem]{ulem}
\usepackage{stackengine}
\usepackage{algorithm}
\usepackage{algpseudocode}
\stackMath
\setlist[enumerate,1]{leftmargin=*}
\setlist[itemize,1]{leftmargin=*}
\usepackage{mathrsfs}
\usepackage[dvipsnames]{color}
\usepackage{todonotes}
\usepackage{float}
\usepackage[pdftex,pdfstartview=FitH,pdfborderstyle={/S/B/W 1},%
colorlinks=true, linkcolor=blue, urlcolor=blue, citecolor=blue,%
pagebackref=true]{hyperref}
\usepackage[nobysame,alphabetic,initials,msc-links]{amsrefs}
\usepackage{imakeidx}
\usepackage{tikz}
\usetikzlibrary{decorations.markings}
\usepackage{comment}
\usepackage{mathtools}

\calclayout

\definecolor{orange}{rgb}{1,0.5,0}

\usepackage{mathrsfs} %Script dla zapisu \mathscr{}

\DeclareMathAlphabet{\mathpzc}{OT1}{pzc}{L}{it} %stylizowany
\theoremstyle{definition}
\newtheorem{definition}{Definition}[section]
\newtheorem{theorem}[definition]{Theorem}

\newtheorem{proposition}[definition]{Proposition}
\newtheorem{corollary}[definition]{Corollary}
\newtheorem{lemma}[definition]{Lemma}

\newtheorem{remark}[definition]{Remark}

\newtheorem*{thm:vanishingCriterion}{Theorem~\ref{thm:vanishingCriterion}}
\newtheorem*{thm:AbramovMainNadler}{Theorem~\ref{thm:AbramovMainNadler}}
\newtheorem*{introthm:horocyclic-Lk}{Theorem~\ref{introthm:horocyclic-Lk}}

\numberwithin{equation}{section}

\def\geq{\geqslant}
\def\leq{\leqslant}
\def\R{\mathbb{R}}

\def\Z{\mathbb{Z}}
\def\N{\mathbb{N}}

\def\R{\mathbb{R}}
\def\id{\mathrm{id}}

\def\Q{\mathbb{Q}}

\def\s{\mathbf{s}}
\def\t{\mathbf{t}}
\def\e{\mathbf{e}}
\def\v{\mathbf{v}}
\def\w{\mathbf{w}}

\def\0{\mathbf{0}}

\def\epsilon{\varepsilon}

\def\mf{\mathfrak}
\def\Lie{\operatorname{Lie}}
\def\ad{\operatorname{ad}}

\def\diag{\operatorname{diag}}

\def\Kak{\operatorname{Kak}}
\def\HKak{\operatorname{HKak}}
\def\TKak{\operatorname{TKak}}

\def\SL{\operatorname{SL}}
\def\vol{\operatorname{vol}}
\def\ess{\operatorname{ess}}
\def\GL{\operatorname{GL}}
\def\Conj{\operatorname{Conj}}
\def\Conf{\operatorname{Config}}

\newcommand{\bea}{\begin{eqnarray}}
  \newcommand{\eea}{\end{eqnarray}}
  \newcommand{\beab}{\begin{eqnarray*}}
  \newcommand{\eeab}{\end{eqnarray*}}

  \newcommand{\be}{\begin{equation}}
  \newcommand{\ee}{\end{equation}}

\newcommand{\set}[1]{\left\lbrace #1 \right\rbrace}
\newcommand{\abs}[1]{\left| #1 \right|}
\newcommand{\mc}{\mathcal}
\newcommand{\norm}[1]{\abs{\abs{#1}}}

\newcommand{\of}{\circ}

\newcommand{\mbf}{\mathbf}
\newcommand{\ve}{\epsilon}
\newcommand{\card}[1]{\operatorname{Card}\left( #1 \right)}

\begin{document}

\title[Kakutani equivalence of abelian unipotent actions]{Kakutani equivalence of higher rank abelian unipotent actions I}

\author[A. Kanigowski]{Adam Kanigowski}
\thanks{A.K. acknowledges the support of the National Science Foundation under grant DMS-2554280}
\address{A.K.: School of Mathematics, Georgia Institute of Technology, Atlanta, GA 30332-0160, USA}
\email{akanigowski3@gatech.edu}

\author[P. Kunde]{Philipp Kunde}
%\thanks{P.K. acknowledges %financial support from a DFG Forschungsstipendium under Grant No. 405305501. 
%}
\address{P.K.: Department of Mathematics, Oregon State University, Corvallis, OR 97331-4605, USA}
\email{kundep@oregonstate.edu}

\author[K. Vinhage]{Kurt Vinhage}
\thanks{K.V. acknowledges support from the NSF (DMS-2405440) and  Alfred P. Sloan Foundation (grant number FG-2025-24955).}
\address{K.V.: Department of Mathematics, The University of Utah, Salt Lake City, UT 84112-0140, USA}
\email{vinhage@math.utah.edu}

\author[D. Wei]{Daren Wei}
\thanks{D.W. acknowledges support by NUS startup grants A-0009806-00-00 and A-0009806-01-00}
\address{D.W.: Department of Mathematics, National University of Singapore, 119076, Singapore}
\email{darenwei@nus.edu.sg}
\date{\today}

\setlength{\marginparwidth }{2 cm}

\begin{abstract}
We study Kakutani equivalence of higher-rank abelian unipotent actions. We construct a Kakutani invariant for general $\R^n$-actions and establish a quantitative vanishing criterion. By this criterion, we prove that ergodic abelian horocyclic actions on finite-volume quotients of connected semisimple Lie groups are loosely Kronecker. We also include a proof of the Abramov entropy formula for time changes of $\R^n$-actions, following the strategy developed by Nadler \cite{NadlerAbramovNflow78}.
\end{abstract}

\maketitle

\setcounter{tocdepth}{1}
\tableofcontents

\section{Introduction}

The classification of measure-preserving transformations and flows on probability spaces up to some definition of equivalence is one of the oldest problems in ergodic theory. The corresponding classification for group actions has received attention as well, but remains much less understood. We aim to study probability preserving abelian group actions up to {\bf Kakutani equivalence}, focusing on cases where classification questions are likely to have positive answers.

Kakutani equivalence is a subtle weakening of measurable isomorphism. For flows and transformations, two of the most natural notions of equivalence are measurable isomorphism and orbit equivalence. For orbit equivalence, the problem was solved by Dye, who showed that any two ergodic aperiodic measure preserving transformations of standard non atomic probability spaces are orbit equivalent. The same result holds for abelian (amenable) group actions: all free ergodic probability preserving actions of countably infinite amenable groups are orbit equivalent.

On the other hand, the story for measurable isomorphism is much more involved. Classification is possible in certain subclasses of systems. Famously, Ornstein showed that two Bernoulli systems are measurably isomorphic if and only if they have the same metric entropy. This also holds for Bernoulli shifts over countably infinite amenable groups \cite{OrnsteinWeissEntropyIsomorphism}. At the other end of the complexity spectrum, Halmos and von Neumann \cite{HalmosNeumann42} showed that ergodic systems with pure point spectrum are classified by their eigenvalues and are all isomorphic to translations on compact abelian groups (such systems are called {\bf Kronecker}). This also generalizes to abelian group actions. In her seminal work, Ratner proved strong rigidity of measurable isomorphisms for homogeneous unipotent flows on finite volume quotients of connected semisimple Lie groups as a consequence of her joining theorem \cite{ratner1990measure}*{Corollary~6}. 

These classifications are exceptional and could only be carried out for specific classes of systems. In fact, recent works have shown that such classifications are impossible in general for flows and transformations \cite{FRW11}.

The case of higher rank ($\R^k$ and $\Z^k$ actions, $k \ge 2$) is very different to rank one situation. Indeed, in the smooth category, the Katok-Spatzier conjecture and related directions point to hyperbolic actions being extensions of products of actions by homogeneous, classifiable systems. When asking about classification as measure-preserving actions, some rigidity phenomena are also expected, and have been shown in some cases \cite{KatokKatokSchmidt02}*{Theorem 5.2, 5.6}.

It is therefore natural to look for intermediate notions of equivalence which imply orbit equivalence and are implied by measurable isomorphism. One such notion is {\bf Kakutani equivalence}, introduced for flows and transformations by Kakutani in \cite{kakutani1943131} and developed by Katok \cite{katok1977monotone} and Ornstein, Rudolph and Weiss \cite{ornstein1982equivalence}. There are many equivalent definitions of Kakutani equivalence. For flows, the definition captures the following notion: one allows for measurable isomorphism {\bf and} a restricted class of time changes which have controlled derivatives along almost every orbit. 

\begin{definition}
    Let $\varphi_t$ be a flow preserving $\mu$, and $\psi_t(x) \coloneq\varphi_{\tau(x,t)}$ a time change which is $C^\infty$ along almost every orbit, with $\tau(x,\cdot)$ an increasing $C^{\infty}$ diffeomorphism and $\rho(x,\cdot)\coloneq\tau(x,\cdot)^{-1}$. We call the time change {\bf integrable} if $\int_X \partial_t\rho(x,0) \, d\mu(x) < \infty$. We say that $\varphi_t$ and $\psi_t$ are {\bf Kakutani equivalent} if there exists an integrable time change of $\varphi^t$ which is measurably conjugated to $\psi^t$.%
\end{definition} 

The corresponding invariant probability measure for $\psi$ has density with Radon-Nikodym derivative $\partial_t\rho(x,0)$. One may relax the $C^\infty$ condition along orbits to the existence of an integrable Radon-Nikodym derivative. This notion of equivalence does not have an obvious extension to $\R^k$ or $\Z^k$-actions. However, several developments contributed to finding a good generalization, we refer the reader to \cites{FeldmanNadler79,DelJuncoRudolph84, KammeyerRudolphRestrictedOrbit,Hasfura-BuenagaIdKakutani,JohnsonSahinRankOneLB98}. The detailed  definition of Kakutani equivalence for group actions will be given in Section \ref{sec:prelim}.

\subsection{Kakutani equivalence of $\R^k$-actions}
A definition of Kakutani equivalence for abelian group actions was introduced in a series of papers in the 1970s and '80s by del Junco, Feldman, Nadler, and Rudolph. We recall their definitions and results using a common notation.

Let $(X,\mathcal{B},\mu)$ be a standard Borel probability space, and  $\alpha:\R^n\curvearrowright(X,\mathcal{B},\mu)$ a measurable, measure-preserving, ergodic, and essentially free action. If $\mbf v \in \R^n$, we let $\alpha_\v: X \to X$ denote the ``time $\v$'' transformation of $X$. A \textbf{reparametrizing cocycle} over $\alpha$ is a measurable map $\rho:X\times \R^n\to\R^n$ such that for $\mu$-a.e. $x\in X$, the map $\rho(x,\cdot)$ is a homeomorphism of $\R^n$, and for every $\v,\w\in\R^n$ and $\mu$-a.e. $x\in X$
\[
\rho(x,\v+\w)=\rho(x,\v)+\rho(\alpha_\v(x),\w).
\]
In other words, a reparametrizing cocycle is an $\R^n$-valued cocycle for which $\rho(x,\cdot)$ is invertible for almost every $x \in X$. If $\rho$ is a reparametrizing cocycle, the \textbf{inverted cocycle} associated with $\rho$ is a measurable map $\tau:X\times \R^n \to \R^n$ such that for every $\v\in \R^n$ and $\mu$-a.e. $x\in X$
\[
\tau(x,\rho(x,\v))=\v.
\]

\medskip

That is, $\tau(x,\cdot)$ is the inverse of $\rho(x,\cdot)$ (which exists, since by assumption, $\rho(x,\cdot)$ is a homeomorphism of $\R^n$). Given a reparametrizing cocycle $\rho$ of $\alpha$ and its inverse $\tau$, the \textbf{time change of $\alpha$ with respect to $\tau$} is the action $\alpha^{\tau}$ defined by
    \[
    \alpha^{\tau}_{\v}(x) := \alpha_{\tau(x,\v)}(x)\qquad \text{ for all $\v\in\R^n$ and $\mu$-a.e. $x\in X$.}
    \]

\begin{remark}
    One of the reasons computations around time changes are subtle is that one does {\bf not} produce a new action by reparameterizing orbits through a cocycle. One must pass to the inverse $\tau$, which is a cocycle over $\alpha^\tau$
\end{remark}

\begin{remark}
    For $\R^n$-actions with $n \ge 2$, Rudolph \cite{RudolphSmoothOrbitEquivalence79} showed a strengthening of Dye's theorem: all free ergodic probability preserving $\R^n$-actions are measurably isomorphic up to a smooth time change. This flexibility implies that additional constraints are required to define a nontrivial equivalence relation.
\end{remark}

Feldman and Nadler \cite{FeldmanNadler79}, and Rudolph \cite{RudolphInner82} extended the definition of Kakutani equivalence to $\R^n$-actions by introducing \textbf{tempered time changes}. If $\rho$ is a reparameterizing cocycle over $\alpha$ and $\rho(x,\cdot)$ and its inverse $\tau(x,\cdot)$ are $C^1$ on a common invariant conull set, we write  
\[
d_\mbf{v}\rho(x)=D_{\v}\rho(x,0) \in \GL_n(\R),\qquad d_\mbf{v}\tau(x)=D_{\v}\tau(x,0).
\]
The chain rule gives $d_{\v}\tau(x)=\big(d_{\v}\rho(x)\big)^{-1}$.

\begin{definition}
    A reparameterizing cocycle with the preceding $C^1$ regularity is {\bf tempered} if 
    \[
    \norm{d_\mbf{v}\rho},\norm{d_{\mbf{v}}\tau} \in L^\infty(X,\mu).
    \]
    Two actions $\alpha_1$ and $\alpha_2$ are called {\bf Kakutani equivalent} if $\alpha_1$ is measurably isomorphic to a tempered time change of $\alpha_2$.
\end{definition}

This condition is among the most intuitive definitions and has the strongest assumptions. Rudolph \cite{RudolphInner82} showed that the tempered condition follows from much weaker ones and can be strengthened by requiring $d_\mbf{v}\rho$ and $d_{\mbf{v}}\tau$ to be almost constant. 

We discuss this history and equivalent definitions in Section \ref{sec:equivalentDef} (see Theorem \ref{thm:kakutani-equivalences}). One could alternatively use a notion of integrable time changes (cf.~Definition~\ref{def:integrableTimeChange}) for flow equivalence, which asks for an $L^1$-property of the derivatives; tempered time changes are a strengthening of this condition by requiring the time change to be differentiable along orbits with $L^\infty$ derivatives. %The resulting equivalence relation (cf.~Definition~\ref{def:kakutani}) is known as \textbf{Kakutani equivalence for $\R^n$-actions}: two actions are equivalent if they are measurably isomorphic up to a tempered time change. 
In Section \ref{sec:flow-compatibility}, we address the compatibility of this framework with the classical definition of Kakutani equivalence for flows.

Although the present work focuses on $\R^n$-actions, we note that the theory of Kakutani equivalence for $\Z^n$-actions is already well-established; we refer the reader to \cites{DelJuncoRudolph84,Hasfura-BuenagaIdKakutani,JohnsonSahinRankOneLB98,JohnsonSahinIsoEx00} for details. Most notably, in a seminal contribution, Kammeyer and Rudolph \cite{KammeyerRudolphRestrictedOrbit} extended this framework to the setting of discrete amenable group actions.

\begin{remark}
Due to the work of Gerber and the second author \cite{gerber2021anti}, from the perspective of classification, Kakutani equivalence of flows and transformations is ``unclassifiable'' in the same way that measurable isomorphism is.
\end{remark}

\subsection{Natural families of higher-rank abelian group actions}
\label{sec:examples}
In this paper, we will focus on higher rank abelian unipotent actions. Given an arbitrary diagonalizable element of a semisimple group, its horocyclic subgroup is always nilpotent, sometimes abelian, but very rarely 1-dimensional. One example of a horocyclic abelian subgroup is

\begin{equation} 
\label{eq:horo-ex}
U = \begin{pmatrix}
    1 & t_1 & t_2 & \dots & t_{n-1} \\
    0 & 1 & 0 & & 0 \\
     & & \ddots  & & \vdots \\
    0 & & \dots & & 1
\end{pmatrix} \subset \SL_n(\R).
\end{equation}

If $\Gamma \subset \SL_n(\R)$ is a lattice, $U$ acts on $\SL_n(\R) / \Gamma$ and is the horocyclic subgroup for the partially hyperbolic element 
\[
a = \diag(2^{n-1},1/2,1/2,\dots,1/2).
\]
When $n = 2$, this is the case of the horocycle flow on a hyperbolic surface, which is loosely Kronecker \cite{ratner1978horocycle}. Furthermore, in the ergodic finite volume semisimple setting, any nontrivial unipotent flow which is loosely Kronecker must come from the horocycle flow on some $\operatorname{SL}_2(\R)$-factor, and act trivially on the other factors \cite{kanigowski2021kakutani}. As a consequence, unipotent flows are loosely Kronecker if and only if they are horocyclic. In fact, a much more striking rigidity phenomenon occurs: if two unipotent flows are not loosely Kronecker, any Kakutani equivalence between them coincides with an affine conjugacy up to time change \cite{lindenstrauss2025time}.

We are interested in the studying the phenomena observed for $\R$ action in the higher rank setting \eqref{eq:horo-ex}. In the current work we focus on identifying higher rank loosely Kronecker actions. We show the first half of this phenomenon in Theorem \ref{introthm:horocyclic-Lk}.

\begin{remark}
Another natural class in which to study Kakutani equivalence is that of higher-rank hyperbolic group actions. One may ask whether the isomorphism rigidity of Katok-Katok-Schmidt \cite{KatokKatokSchmidt02} persists in Kakutani equivalence.
\end{remark}

\subsection{Statement of Results}

We prove three results in this paper. The first result is due to Nadler unpublished manuscript; we include a proof based on his unpublished manuscript in the appendix. We are grateful to Benjamin Weiss for sharing Nadler's original exposition with us. We also generalize the notion of loosely Kronecker and the criteria for this property, and prove that horocyclic actions are loosely Kronecker.

\subsubsection{Nadler's Abramov formula}
A central feature of classical Kakutani equivalence is the preservation of entropy classes; specifically, the classes of zero, finite positive, and infinite entropy are invariant under equivalence. This invariance stems from the Abramov entropy formula for suspension flows. Accordingly, extending this phenomenon to $\R^n$-actions requires a generalization of the Abramov entropy formula to time changes of $\R^n$-actions. 

This generalization is due to Nadler \cite{NadlerAbramovNflow78}. Since the manuscript remains unpublished, we include a comprehensive proof in Appendix~\ref{sec:AppendixAbramovFormula}. We follow Nadler's strategy and provide necessary technical details.
\begin{theorem}[Nadler's Abramov formula \cite{NadlerAbramovNflow78}]\label{thm:AbramovMainNadler}
    Let $\alpha$ be a measurable, measure-preserving, ergodic, and free $\R^n$-action on a standard non-atomic probability space $(X,\mathcal{B},\mu)$, and $\alpha^{\tau}$ a tempered time change of $\alpha$. Then 
    \[
     h_{\mu}(\alpha)= h_{\mu^{\tau}}(\alpha^{\tau})\cdot\int_X\lambda_\tau(x)\,d\mu.
    \]
    Here, $\mu^{\tau}$ is the invariant probability measure of $\alpha^{\tau}$ induced by the time change, with %the function $\Delta_\rho(x) = \Delta_\rho(x,\mbf 0) = \det (d_{\mbf v}\rho)$ yields 
    \[
    \frac{d\mu^{\tau}}{d\mu}(x)=\frac{\lambda_\tau(x)}{\int_X\lambda_\tau(y)\,d\mu(y)}
    \]
    as in \ref{item:a4} of Section \ref{sec:prelim}.%the Radon-Nikodym derivative of $\mu^{\tau}$ with respect to $\mu$.
\end{theorem}

\subsubsection{Minimal Kakutani entropy} Among all Kakutani equivalence classes for flows, one class is distinguished by its structural simplicity: \textbf{loosely Kronecker flows} (also known as zero-entropy loosely Bernoulli systems or standard flows).  This class is characterized by a fundamental minimality property\footnote{Cf.~\cite{ornstein1982equivalence}*{Theorem~4.1} and \cite{katok1977monotone}*{Theorem~1} for more details.}: a flow is loosely Kronecker if and only if every ergodic flow is Kakutani equivalent to an extension of it. Typical examples of this equivalence class include irrational flows on the two-dimensional torus and ergodic translation flows on compact translation surfaces, cf.~\cites{ornstein1982equivalence,katok1977monotone}. This concept was generalized to $\R^n$-actions by Feldman and Nadler \cite{FeldmanNadler79}, whose framework includes ergodic Kronecker $\R^n$-actions among the loosely Kronecker systems.

Following Ratner's strategy for proving the loosely Kronecker property \cite{ratner1981some}, in \S\ref{sec:KakutaniInvariant} we construct the {\bf Kakutani invariant} or {\bf Kakutani slow entropy} $h_\mu^{\operatorname{Kak}}(\alpha,s,\mathcal{L})$ for $\R^n$-actions and establish its basic properties. This invariant depends on a linear transformation $\mathcal{L}\in\GL_n(\R)$ and a  \textbf{scaling function} $s$ belonging to the class 
\[
\mathcal{F}=\{s:\R^+\to\R^+:\text{nondecreasing, unbounded functions}\}.
\]
The transformation $\mathcal{L}$ governs the geometry of the averaging sets via the F\o lner sequence $\{\mathcal{L}(C_R)\}_{R>0}$, where $C_R=[-R,R]^n$ denotes the standard hypercube. We will in fact show that if the scale $s$ is sufficiently small, the Kakutani slow entropy is independent of $\mc L$ (see Proposition \ref{prop:LinearInvariance}). Let $\mc S$ denote the set of slowly growing scaling functions, so 
\begin{equation}\label{eq:SFamily}
\mc S = \set { s  \in \mathcal{F}\cap C^\infty(\R_+,\R_+): \lim_{t \to \infty} s(t) = \infty, \lim_{t\to \infty} \frac{s(at)}{s(t)} = 1 \mbox{ for all }a >0}.
\end{equation}

Note that $t\mapsto\log(1+t)$ belongs to $\mc S$, whereas no power function $t \mapsto t^\theta$, $\theta > 0$, belongs to $\mathcal{S}$. For $s\in\mathcal{S}$, we denote the invariant by $h_{\mu,s}^{\operatorname{Kak}}(\alpha).$ Our first main result establishes a quantitative characterization of loosely Kronecker $\R^n$-actions, which generalizes \cite{ratner1981some}*{Theorem~1}:
\begin{theorem}\label{thm:vanishingCriterion}
   Let $\alpha$ be a zero-entropy ergodic $\R^n$-action. Then $\alpha$ is loosely Kronecker if and only if  
\[
h_{\mu,s}^{\operatorname{Kak}}(\alpha)=0\qquad  \text{for all $s\in \mathcal{S}$.}
\]
\end{theorem}

 \begin{remark}
    Theorem \ref{thm:vanishingCriterion} can also be considered a generalization of Ferenczi's trivial growth theorem: $h_{\mu,s}(f) = 0$ for all growth functions $s$ if and only if $f$ is Kronecker, where $h_{\mu,s}(f)$ is the slow entropy of $f$.
 \end{remark}

\subsubsection{Horocyclic actions}

Let $G$ be a semisimple Lie group and $\Gamma<G$ a lattice. Ratner's seminal papers \cites{ratner1978horocycle,ratner1979cartesian,ratner1981some} on the Kakutani equivalence of horocycle flows---the simplest unipotent flows on $G/\Gamma$---laid the foundation for the broader study of unipotent dynamics \cites{kanigowski2021kakutani,lindenstrauss2025time}. In this paper, we investigate a ``higher-rank'' analogue. while more detailed definition is provided in Definition~\ref{def:horocyclicAction}: 
\begin{definition}
    The \textbf{horocyclic action} is generated by the full unstable Lie subalgebra of a real semisimple element in $\Lie(G)$. 
\end{definition}
Our next theorem generalizes \cite{ratner1978horocycle}, and we view this result as a first step toward establishing the general theory of Kakutani equivalence for higher-rank abelian unipotent actions. 
\begin{theorem}
\label{introthm:horocyclic-Lk}
    If $G$ is a connected semisimple Lie group and $U \subset G$ is an abelian horocyclic subgroup, and $\Gamma < G$ is a  lattice, then the homogeneous action $\alpha : U \curvearrowright G / \Gamma $ defined by $\alpha(u)\cdot g\Gamma = (ug)\Gamma$ is loosely Kronecker whenever it is ergodic for Haar probability measure.
\end{theorem}

\subsection*{Structure of the paper}
\S\ref{sec:prelim} lays the groundwork for studying time changes of $\R^n$-actions. In \S\ref{sec:equivalentDef}, we formalize Kakutani equivalence for $\R^n$-actions, establishing the equivalence of several definitions and verifying consistency with the standard theory for $\R$-actions. \S\ref{sec:KakutaniInvariant} introduces a generalized Kakutani invariant for $\R^n$-actions; we establish its fundamental properties---including invariance and the generating property---and derive a vanishing criterion corresponding to loosely Kronecker systems. We also discuss the relationship between this generalized invariant and Ratner's classical invariant for flows. In \S\ref{sec:horocyclicAction}, we demonstrate that horocyclic actions are loosely Kronecker. Finally,  \S\ref{sec:AppendixAbramovFormula} presents a proof of the Abramov entropy formula for time changes of $\R^n$-actions, largely following the strategy developed by Nadler \cite{NadlerAbramovNflow78}.

\subsection*{Acknowledgement}
We are indebted to Benjamin Weiss for sharing with us the unpublished manuscript by D.~Nadler, which served as the primary guide for the proof presented in \S\ref{sec:AppendixAbramovFormula}. Although we were unable to contact Nadler directly, we fully acknowledge that the strategy for the Abramov entropy formula proof is due to him. We are grateful to Aaron Brown, Elon Lindenstrauss and Federico Rodriguez Hertz for enlightening discussions regarding this work.

\section{Preliminaries}\label{sec:prelim}
In this section, we recall some properties of time changes of measure preserving $\R^n$-actions, closely following \cite{FeldmanNadler79}. We use jointly measurable versions of the action and cocycles on a common invariant conull set. Lebesgue measure on $\R^n$ is denoted by $m$.

\subsection{Integrable time changes} Throughout this paper, except in \S\ref{sec:horocyclicAction}), we work under the following minimal assumptions.

\begin{definition}\label{def:integrableTimeChange}
A time change $\alpha^{\tau}$ is called \textbf{integrable} if the underlying system and cocycle satisfy the following conditions:
\begin{enumerate}[label=(A\arabic*)]
    \item\label{item:a1} $(X,\mathcal{B},\mu)$ is a standard Borel probability space,
    \item\label{item:a2} $\alpha:\R^n\curvearrowright(X,\mathcal{B},\mu)$ is a measurable, measure-preserving, ergodic, and free action,
    \item\label{item:a3} $\rho$ is a reparametrizing cocycle of $\alpha$ with inverse cocycle $\tau$ such that the maps $\tau_x(\mathbf{v}) := \tau(x,\mathbf{v})$ take the Lebesgue measure class to the Lebesgue measure class for almost every $x\in X$, and
        \item\label{item:a4} $0<\int_X\lambda_\tau \,d\mu<\infty$, where $\lambda_\tau$ is defined below.
\end{enumerate}
\end{definition}

For a Borel set $E\subset\R^n$, let 
\[
q_x(E)=m(\rho_x(E))=(\tau_x)_*m(E),\qquad \hat{\lambda}_{\tau}(x,\v)=\frac{dq_x}{dm}(\v).
\]

The measures $q_x$ are locally finite and equivalent to $m$. A jointly measurable version of their densities exists, but evaluation at $\v=0$ needs justification. The cocycle identity gives
\[
q_x(E+\w)=q_{\alpha_{\w}.x}(E).
\]
After choosing the version below, the recentering formula is 
\begin{equation}\label{eq:recenter}
    \hat{\lambda}_{\tau}(x,\w)=\lambda_{\tau}(\alpha_{\w}.x)\qquad \text{for $\mu\times m$-a.e. $(x,\w)$.}
\end{equation}

\begin{lemma}
\label{lem:lambda}
    There is a measurable function $\lambda_{\tau}:X\to(0,\infty)$, defined almost everywhere, such that for almost every $x$ and every Borel set $E\subset\R^n$,
    \begin{equation}\label{eq:rhoxeIntegralFormula}
       m(\rho_x(E))=\int_E\lambda_{\tau}(\alpha_{\v}.x)dm(\v). 
    \end{equation}
    The measure $\nu(B)\coloneq\int_B\lambda_{\tau}d\mu$ is invariant under $\alpha^{\tau}$ and equivalent to $\mu$. In addition, if $\rho_x$ is a $C^1$ diffeomorphism, then we have 
    \[
    \lambda_{\tau}(x)=\abs{\det\partial_{\v}\rho(x,\0)},\qquad \text{$\mu$-a.e. $x$.}
    \]
\end{lemma} 

\begin{proof}
    Let $C_r = [-r,r]^n$ and define
    \[
    \lambda_{\tau}(x) = \lim_{r\to0}\frac{q_x(C_r)}{m(C_r)}\in(0,\infty).
    \] 
    It follows from the Lebesgue differentiation and covariance of $q_x$, the above limit exists at $\alpha_{\v}.x$ for $m$-a.e. $\v$, for almost every $x$. Combining the invariance of $\mu$ and Fubini's theorem, we obtain that it exists for $\mu$-almost every $x$ and $m$-a.e. $\v$,
    \[
    \frac{dq_x}{dm}(\v)=\lim_{r\to0}\frac{q_x(\v+C_r)}{m(C_r)}=\lim_{r\to0}\frac{q_{\alpha_{\v}.x}(C_r)}{m(C_r)}=\lambda_{\tau}(\alpha_{\v}.x).
    \]
    This proves \eqref{eq:rhoxeIntegralFormula} and \eqref{eq:recenter}. It also follows from the covariance that equation~\eqref{eq:rhoxeIntegralFormula} holds on a common invariant conull set.

    Consider a measurable flow box $x=\alpha_{\v}.z$, we obtain the local product measure $d\eta(z)dm(\v)$ from the invariance of $\mu$, where $\eta$ is the transverse measure. It follows from the change of orbit coordinate $\mathbf{u}=\rho(z,\v)$, equation \eqref{eq:rhoxeIntegralFormula} shows that
    \[
    d\nu=d\eta(z)dm(\mathbf{u}).
    \]
    The action $\alpha^{\tau}$ translates $\mathbf{u}$. On overlapping boxes, the cocycle identity makes the change of the $\mathbf{u}$ coordinate a translation, and the transverse holonomy preserves $\eta$, as it does for $\mu$. Consequently the local product measures are invariant under $\alpha^{\tau}$. A countable family of measurable flow boxes covers a conull set for a free measurable $\R^n$-action, so this proves invariance globally. Positivity of $\lambda_{\tau}$ gives equivalence to $\mu$, and the final assertion is the change of variables formula.
\end{proof}

    For an integrable time change, define $c_{\tau}\coloneq\int_X\lambda_{\tau}d\mu$. Its invariant probability measure is
    \begin{equation}\label{eq:NewMeasureDensityEquation}
        \frac{d\mu^{\tau}}{d\mu}(x)=\frac{\lambda_{\tau}(x)}{c_{\tau}}.
    \end{equation}
    It is ergodic as $\alpha$ and $\alpha^{\tau}$ have the same orbits and equivalent measures.

\subsection{The mean vector map}

Fix an $\R^n$-valued cocycle $\rho$ over an ergodic $\R^n$-action $\alpha : \R^n \curvearrowright X$ which determines an integrable time change of $\alpha$. Whenever $\rho(\cdot,\v)$ is integrable for every $\v$, its mean vector map $J(\rho):\R^n\to\R^n$ is defined by
\[
J(\rho)(\v)=\int_X\rho(x,\v) \,d\mu.
\]
We first recall a classical result on ergodic directions. A vector $\v\in\mathbb{R}^n$ is said to be ergodic if the transformation $\alpha_{\v}$ is ergodic.  
\begin{theorem}[Theorem $1$ in \cite{PughShub71}]\label{thm:PughShub}
If $\mathbb{R}^n$ acts ergodically on $(X,\mu)$, $\mu(X)<\infty$ and $L^2(X,\mu)$ is separable, then all the elements of $\mathbb{R}^n$, off a countable family of affine hyperplanes, are ergodic.
\end{theorem}

The following proposition summarizes the properties of $J(\rho)$ under a mild growth assumption.
\begin{proposition}\label{prop:jRhoProp}
    Suppose that there exists a constant $C>0$ such that $\norm{\rho(x,\v)}\leq C\norm{\v}$ for $\mu$-a.e. $x\in X$ and all $\v\in\R^n$. Then
    \begin{enumerate}[label=(\roman*)]
        \item\label{item:JrhoLinear} $J(\rho)$ is an $\R$-linear map;
        \item\label{item:JrhoErgodicSum} for every $\epsilon>0$ and every $\delta\in(0,1)$, there exists $N_{\epsilon,\delta}\geq 0$ such that for all $N>N_{\epsilon,\delta}$,
\[
\mu(\{x\in X:\norm{\rho(x,\v)-J(\rho)(\v)}<\epsilon\norm{\v}\text{ for all $\norm{\v}\geq N$}\})>1-\delta;
\]
        \item\label{item:JrhoDet} $\int_X\lambda_\tau\,d\mu=\abs{\det J(\rho)}$.
    \end{enumerate}
\end{proposition}

\medskip

We first prove that $J(\rho)$ is linear.
\begin{proof}[Proof of Proposition~\ref{prop:jRhoProp}~\ref{item:JrhoLinear}]
For any $\v,\w\in\R^n$, the cocycle identity $\rho(x,\v+\w)=\rho(x,\v)+\rho(\alpha_{\v}.x,\w)$ and the $\alpha$-invariance of $\mu$ imply
\begin{equation*}
        \begin{aligned}
            J(\rho)(\v+\w)=&\int_X\rho(x,\v)d\mu+\int_X\rho(\alpha_{\v}.x,\w)d\mu=J(\rho)(\v)+J(\rho)(\w).
        \end{aligned}
    \end{equation*}
    This additivity implies that $J(\rho)(q\v)=qJ(\rho)(\v)$ for all $q\in\Q$ and $\v\in\R^n$.

    To extend this to $\R$, it remains to show that $J(\rho)$ is continuous. The hypothesis $\int_X\norm{\rho(x,\v)}d\mu\leq C\norm{\v}$ yields 
    \[
    \norm{J(\rho)(\v)}\leq \int_X\norm{\rho(x,\v)}d\mu\leq C\norm{\v}\quad \text{ for all $\v\in\R^n$.}
    \]
    Thus, $J(\rho)$ is Lipschitz continuous at $\0$, and by additivity, is globally Lipschitz. Since every continuous additive map is linear over $\R$, this proves  Proposition~\ref{prop:jRhoProp}~\ref{item:JrhoLinear}.
\end{proof}

\begin{proof}[Proof of Proposition~\ref{prop:jRhoProp}~\ref{item:JrhoErgodicSum}]
    By Theorem~\ref{thm:PughShub}, choose a basis $\v_1,\ldots,\v_n$ such that each transformation $\alpha_{\v_i}$ is ergodic. Fix $B\geq 1$ such that
    \[
    \sum_{i=1}^n\norm{\v_i}\leq B,\qquad \sum_{i=1}^n\abs{k_i}\norm{\v_i}\leq B\norm{\sum_{i=1}^nk_i\v_i},\qquad k_1,\ldots,k_n\in\R.
    \]
    For the given $\epsilon>0$ and $\delta\in(0,1)$, choose
    \[
    0<\eta<\min\left(10^{-1},\frac{\epsilon}{16B(C+1)}\right),\qquad 0<\kappa<\min\left(\delta,\frac{\eta}{16}\right).
    \]
    The pointwise ergodic theorem and the cocycle identity imply, for $\mu$-a.e. $x\in X$ and each $i$,
    \begin{equation*}
        \begin{aligned}
            \lim_{k\to\infty}\frac{\rho(x,k\v_i)}{k}=&\lim_{k\to\infty}\frac{1}{k}\sum_{j=0}^{k-1}\rho(\alpha_{j\v_i}.x,\v_i)
            =J(\rho)(\v_i).
        \end{aligned}
    \end{equation*}
    Let $X_0$ be the full measure set where the hypothesis $\norm{\rho(x,\v)}\leq C\norm{\v}$ holds. By Egorov's theorem and the linearity of $J(\rho)$ (established in Part \ref{item:JrhoLinear}), there is a set $W_{n+1}\subset X_0$ with $\mu(W_{n+1})>1-\frac{\kappa}{n+1}$ and an integer $N_0\geq1$ such that for all $x\in W_{n+1}$, $i=1,\ldots,n$, and $\abs{k}\geq N_0$,
    \begin{equation}\label{eq:ergodicAverage}
    \norm{\rho(x,k\v_i)-J(\rho)(k\v_i)}\leq \eta\abs{k}\norm{\v_i}.
    \end{equation}

\medskip

    We construct nested sets $W_1\subset\ldots\subset W_{n+1}$, proceeding from $i=n$ down to $i=1$. Apply the ergodic theorem for $\alpha_{\v_i}$ to the indicator of $X\setminus W_{i+1}$. Since $\mu(X\setminus W_{i+1})<\kappa<\eta/16$, we may discard a set of measure less than $\kappa/(n+1)$ from $W_{i+1}$ and choose an integer $N_i$ such that
    \begin{enumerate}[label=(\alph*)]
        \item\label{item:recurrenceA} $\mu(W_{i})\geq 1-\frac{(n+2-i)\kappa}{n+1}>1-\delta$;
        \item\label{item:recurrenceB} for every $x\in W_i$ and $M\geq N_{i}$, 
         \[
    \card{\{\ell\in\Z:\abs{\ell}\leq M,\alpha_{\ell\v_i}.x\notin W_{i+1} \}}\leq \frac{\eta}{8}(2M+1).
    \]
    \end{enumerate}

    Choose an integer $L$ so that 
    \[
    L\geq \max(2N_0,N_1,\ldots,N_n,2/\eta,100),
    \]
    and define 
    \[
    N_{\epsilon,\delta}=\frac{L}{\eta}.
    \]

    Fix $x\in W_1$ and $\v=\sum_{i=1}^nk_i\v_i$ with $\norm{\v}\geq N_{\epsilon,\delta}$. Let $I_{\v}$ be the set of indices with large coefficients:
    \[
    I_{\v}=\{i\in\{1,\ldots,n\}:\abs{k_i}\geq L\}.
    \]
    It follows from the recurrence property \ref{item:recurrenceB} that we can choose perturbed coefficients $k'_i=k'_i(x)\in\Z$ successively so that
   \begin{equation}\label{eq:kixExistence}
   \begin{gathered}
        \alpha_{\sum_{j=1}^{i}k'_j\v_j}.x\in W_{i+1},\qquad 1\leq i\leq n, \\
       k'_i=0 \text{ if } i\notin I_{\v},\qquad  \abs{k'_i-k_i}\leq \eta\abs{k_i}\text{ if }i\in I_{\v}.\\
        \end{gathered}
    \end{equation}
    The following is the detailed chosen procedure:
    \begin{itemize}
        \item Notice that if $\abs{k_i}<L$, the choice $k'_i=0$ is allowed by $W_i\subset W_{i+1}$.
        \item If we have $\abs{k_i}\geq L$, then apply \ref{item:recurrenceB} at the current base point with $M=\lceil(1+\eta)\abs{k_i}\rceil$. Notice that the interval $[k_i-\eta\abs{k_i},k_i+\eta\abs{k_i}]$ contains at least $\eta\abs{k_i}$ integers, whereas the number of good returns in $[-M,M]$ is larger than
    \[
    1-\frac{\eta}{8}\big(2(1+\eta)\abs{k_i}+3\big)\geq1-\frac{\eta}{3}\abs{k_i}.
    \]
    Thus an admissible $k'_i$ exists. Moreover, we always have $\abs{k'_i}\geq (1-\eta)L\geq N_0$ whenever $k'_i\neq0$.
    \end{itemize}
    
    Let $\v'=\sum_{i=1}^nk'_i\v_i$ and $\w_i=\sum_{j=1}^ik'_j\v_j$ with $\w_0=\0$. We estimate
    \begin{equation}\label{eq:decomposeTauJ}
        \begin{aligned}
            \norm{\rho(x,\v)-J(\rho)(\v)} 
            \leq&\norm{\rho\left(x,\v\right)-\rho\left(x,\v'\right)}+\norm{\rho\left(x,\v'\right)-J(\rho)\left(\v'\right)} \\
            & +\norm{J(\rho)\left(\v'\right)-J(\rho)\left(\v\right)}.
        \end{aligned}
    \end{equation}
    The coefficient bounds give
    \[
    \norm{\v-\v'}\leq LB+\eta B\norm{\v}\leq 2\eta B\norm{\v}.
    \]
    Since $\alpha_{\v'}.x\in W_{n+1}\subset X_0$, the cocycle identity and the hypothesis $\norm{\rho(x,\v)}\leq C\norm{\v}$ imply
    \begin{equation*}
        \max(\norm{\rho(x,\v)-\rho(x,\v')}, \norm{J(\rho)(\v-\v')})\leq 2C\eta B\norm{\v}.
    \end{equation*}

   For the middle term of \eqref{eq:decomposeTauJ}, every $\alpha_{\w_{i-1}}.x$ belongs to $W_i\subset W_{n+1}$, so \eqref{eq:ergodicAverage} gives
    \begin{equation*}
        \begin{aligned}
            \norm{\rho\left(x,\v'\right)-J(\rho)\left(\v'\right)}&
            \leq\sum_{i=1}^{n}\norm{\rho\left(\alpha_{\w_{i-1}}.x,k'_{i}\v_{i}\right)-J(\rho)\left(k'_{i}\v_{i}\right)}\\
            &\leq\eta\sum_{i=1}^n\abs{k'_i}\norm{\v_i}\leq2\eta B\norm{\v}.
        \end{aligned}
    \end{equation*}
    Combining these bounds gives an error at most $(4C+2)\eta B\norm{\v}<\epsilon\norm{\v}$. Since $\mu(W_1)>1-\delta$, this proves Proposition~\ref{prop:jRhoProp}~\ref{item:JrhoErgodicSum}.
\end{proof}

We finally prove  Proposition~\ref{prop:jRhoProp}~\ref{item:JrhoDet}:
\begin{proof}[Proof of Proposition~\ref{prop:jRhoProp}~\ref{item:JrhoDet}]
    By Lemma~\ref{lem:lambda},
    \[
    m(\rho(x,C_N))=\int_{C_N}\lambda_{\tau}(\alpha_{\v}.x)\,dm(\v).
    \]
    The pointwise ergodic theorem therefore gives
    \[
    \frac{m(\rho(x,C_N))}{m(C_N)}\longrightarrow \int_X\lambda_{\tau}d\mu>0,\qquad \text{$\mu$-a.e. $x$}.
    \]
    By \ref{item:JrhoErgodicSum}, first letting $\delta\downarrow0$ for each $\epsilon=1/j$ and then letting $j\to\infty$, the cocycle is uniformly sublinear relative to its mean for $\mu$-a.e. $x$. Fix such an $x$ for which the volume limit and the growth bound $\norm{\rho(x,\v)}\leq C\norm{\v}$ also hold. Uniform sublinearity and the growth bound near the origin imply
    \[
    \sup_{\v\in C_1}\norm{N^{-1}\rho(x,N\v)-J(\rho)\v}\to0.
    \]
    If $J(\rho)$ is singular, $N^{-1}\rho(x,C_N)$ would lie in a neighborhood of $J(C_1)$ whose width goes to zero. Its volume would tend to zero, contradicting to the positive volume limit above. Thus $J(\rho)$ is nonsingular.

    Define 
    \[
    F_N(\v)\coloneq N^{-1}J(\rho)^{-1}\rho(x,N\v),\qquad \epsilon_N\coloneq\sup_{\v\in C_1}\norm{F_N(\v)-\v}_{\infty}.
    \]
    Then $\epsilon_N\to0$ and $F_N(C_1)\subset (1+\epsilon_N)C_1$. For all sufficiently large $N$, we have $\epsilon_N<1$. If $\w\in(1-\epsilon_N)C_1$, the continuous map $\v\mapsto \w-(F_N(\v)-\v)$ maps $C_1$ to itself.  Brouwer's fixed point theorem gives $F_N(\v)=\w$ for some $\v\in C_1$ and hence $(1-\epsilon_N)C_1\subset F_N(C_1)$. Taking volumes gives
    \begin{equation*}
    (1-\epsilon_N)^n\abs{\det J(\rho)}\leq \frac{m(\rho(x,C_N))}{m(C_N)}\leq (1+\epsilon_N)^n\abs{\det J(\rho)}.
    \end{equation*}
    Letting $N\to\infty$ proves Proposition~\ref{prop:jRhoProp}~\ref{item:JrhoDet}.
\end{proof}

\section{Kakutani equivalence for \texorpdfstring{$\R^n$}{Rn}-actions}\label{sec:equivalentDef}

As discussed in the introduction, the higher-rank setting ($n\geq 2$) requires restricting the class of permissible time changes to avoid triviality of the equivalence relation. In this section, we explore the necessary additional restrictions needed to develop the theory of Kakutani equivalence for $\R^n$-actions. We also discuss the equivalence of various definitions and the compatibility of higher-rank definitions with $\R$-actions.

\subsection{Equivalence of various definitions}
This section presents various equivalent definitions of Kakutani equivalence for $\R^n$-actions, as established in \cites{FeldmanNadler79,RudolphInner82}. Although most of the results presented here have already appeared in \cites{FeldmanNadler79,RudolphInner82} in varying forms, we plan to provide a unified presentation for completeness and future development.

We first specify cohomology. Given two reparametrizing cocycles $\rho$ and $\tilde{\rho}$, we say they are cohomologous if there exists a measurable map $p:X\to\R^n$ such that
\[
\tilde{\rho}(x,\v)=\rho(x,\v)+p(x)-p(\alpha_{\v}.x).
\]
Here both cocycles are required to define integrable time changes. The cohomologous relation indeed induces an orbit preserving isomorphism
\[
H(x)\coloneq\alpha_{-p(x)}^{\tau}.x.
\]
from $(X,\mu^{\tilde{\tau}},\alpha^{\tilde{\tau}})$ to $(X,\mu^{\tau},\alpha^{\tau})$. It follows from $\alpha^{\tau}$-coordinate of $H(\alpha_{\v}.x)$ is $\tilde{\rho}(x,\v)-p(x)$ that $H$ is a measurable bijection on each orbit, and also intertwines the action follows. Apply the measurable flow box calculation in Lemma~\ref{lem:lambda} to this translation of orbit coordinates gives
\[
H_*(\lambda_{\tilde{\tau}}\mu)=\lambda_{\tau}\mu.
\]
This in particular gives that $c_{\tau}=\int_X\lambda_{\tau}d\mu=\int_X\lambda_{\tilde{\tau}}d\mu=c_{\tilde{\tau}}$ and hence $H$ transports the normalized probability measures as asserted.

When both cocycles have integrable values, cohomology also preserves their corresponding mean maps, even if $p$ is not integrable. For any fixed $\v$, the Birkhoff averages of $\tilde{\rho}(x,\v)-\rho(x,\v)$ equal 
\[
\frac{p(x)-p(\alpha_{k\v}.x)}{k},
\]
which goes to zero in probability by invariance of $\mu$ and finiteness of $p$. Its almost sure Birkhoff limit is therefore zero, so its integral is zero.

We give the definitions used in the equivalence theorem.

\begin{definition}\label{def:temperedTimeChange}
\begin{enumerate}[label=(\roman*)]
    \item \label{def:barelyLinearTimeChange} An integrable time change $\alpha^{\tau}$ is called \textbf{tempered} if there exists an $\alpha$-invariant set $A\subset X$ of full measure such that for all $x\in A$, the maps $\v\mapsto\tau(x,\v)$ and $\v\mapsto\rho(x,\v)$ are $C^1$, and the norms $\norm{\partial_{\v} \tau(x,\0)}$ and $\norm{\partial_{\v} \rho(x,\0)}$ are bounded functions of $x$ on $A$. We say that $\alpha^\tau$ is \textbf{cohomologically tempered} if it is cohomologous to a tempered time change.
    \item\label{item:Rtempered} An integrable time change $\alpha^{\tau}$ is said to be \textbf{$\ve$-tempered} for some $\ve>0$ if it is tempered, and for every $x\in A$,
    \[
    \norm{\partial_{\v} \tau(x,\0)-\int_X\tau(y,\cdot)d\mu^{\tau}(y)}<\ve, \quad \norm{\partial_{\v} \rho(x,\0)-\int_X\rho(y,\cdot)d\mu(y)}<\ve.
    \]
    We say that an integrable time change is {\bf $\ve$-cohomologically tempered} if it is cohomologous to an $\ve$-tempered time change.
    \item  \label{def:integrableLipTimeChange}   An integrable time change $\alpha^{\tau}$ is said to be \textbf{integrably Lipschitz} if there exists a constant $C>0$ such that for all $\v\in\R^n$
    \[
    \int_X\norm{\tau(x,\v)}d\mu^{\tau}\leq C\norm{\v},\quad \text{and}\quad \int_X\norm{\rho(x,\v)}d\mu\leq C\norm{\v}.
    \]
    \item \label{def:abvPinning} Let $\v\in \R^n$. A set $A\subset X$ of positive measure is called an \textbf{``$(a,b,\v')$''-pinning set} for $\v$ (where $a,b\in\R$ and $\v'\in\R^n$) if for any $k\in\N$ and any $x\in A$ such that $\alpha_{k\v}.x\in A$, we have
        \[
        \norm{\rho(x,k\v)-k\v'}<ak+b.
        \]
    \item An integrable time change $\alpha^{\tau}$ is called \textbf{linear pinning} if there exist constants $a,b\in\R$ and a basis $\{\v_1,\ldots,\v_n\}$ of $\R^n$ such that for each $i=1,\ldots,n$, the vector $\v_i$ admits an $(a,b,\v'_i)$-pinning set for some vector $\v'_i\in\R^n$.
\end{enumerate} 
\end{definition}

\begin{remark}
    In \cite{RudolphInner82}, cohomologically tempered is called barely linear.
\end{remark}

We aim to prove the following theorem, which establishes four equivalent ways to define Kakutani equivalence of $\R^n$-actions:

\begin{theorem}
\label{thm:kakutani-equivalences}
    Let $\alpha^{\tau}$ be an integrable time change of an ergodic action $\alpha$. The following are equivalent:
    \begin{enumerate}[label=(\alph*)]
        \item $\alpha^{\tau}$ is  cohomologically tempered;
        \item $\alpha^{\tau}$ is $\ve$-cohomologically tempered for every $\epsilon > 0$;
        \item $\alpha^{\tau}$ is cohomologous to an integrably Lipschitz time change;
        \item $\alpha^{\tau}$ is linear pinning.
    \end{enumerate}
\end{theorem}

\begin{definition}\label{def:kakutani}
    We say that two $\R^n$-actions $\varphi$ and $\alpha$ are \textbf{Kakutani equivalent} if $\varphi$ is measurably isomorphic to a time change $\alpha^{\tau}$ of $\alpha$ that satisfies the equivalent conditions of Theorem~\ref{thm:kakutani-equivalences}. Furthermore, they are \textbf{evenly Kakutani equivalent} if the time change can be chosen with a tempered representative satisfying $\int_X\rho(x,\mbf{v}) \,d\mu=\mbf{v}$ for all $\mbf{v} \in \R^n$, where $\rho$ is the inverse cocycle of $\tau$. The mean map is independent of the representative whenever both cocycles have integrable values.
\end{definition}    

    % We make some initial observations before proceeding with the proof of Theorem \ref{thm:kakutani-equivalences}. Note that tempered time changes are automatically integrable. The condition in \ref{item:Rtempered} admits a natural dynamical interpretation: it requires that the local geometry of the time change stays uniformly close to its global average behavior. The derivative $\partial_{\v}\tau(x,\0)$ captures the instantaneous rate of change of the reparametrization at the point $x$. On the other hand, by the pointwise ergodic theorem, the integral term $\int_X\tau(y,\cdot)\,d\mu^{\tau}(y)$ represents the asymptotic average growth of the cocycle along generic orbits. Thus, an $R$-tempered time change essentially demands that the instantaneous distortion at almost every point does not deviate from the long-term mean velocity by more than $R$.

%\medskip

The condition in \ref{item:Rtempered} requires the derivatives of the reparametrizing cocycles to stay uniformly close to their mean maps. We break the proof of Theorem \ref{thm:kakutani-equivalences} into a series of propositions.

\begin{proposition}\label{prop:temperedImpliesIntLip}
    If $\alpha^{\tau}$ is a tempered time change of $\alpha$, then it is an integrably Lipschitz time change of $\alpha$.
\end{proposition}
\begin{proof}
Let $A\subset X$ be the shift-invariant set of full measure from Definition~\ref{def:temperedTimeChange} where the derivatives are defined and bounded. Recall the cocycle identity $\rho(x,\v+\w)=\rho(x,\w)+\rho(\alpha_{\w}.x,\v)$. Differentiating this equation with respect to the second variable at $\v=\0$ yields for all $\w\in\R^n$ and $x\in A$,
\[
\partial_{\v}\rho(x,\w)=\partial_{\v}\rho(\alpha_{\w}.x,\0).
\]
Since $\alpha^{\tau}$ is tempered, the norm $\norm{\partial_{\v} \rho(x,\0)}$ is bounded by some constant $K$ for $x\in A$. The shift-invariance of $A$ implies that $\alpha_{\w}.x\in A$ for any $\w\in \R^n$ and $x\in A$. Consequently, for all $\w\in\R^n$ and $x\in A$,
\[
\norm{\partial_{\v}\rho(x,\w)}\leq K.
\]

Since $\rho(x,\0)=\0$, the mean value theorem implies for all $\v\in \R^n$ and $x\in A$,
\[
\norm{\rho(x,\v)}\leq K\norm{\v}.
\]
Integrating this inequality over $X$ gives the result for $\rho$. The same argument for the cocycle $\tau$ over $\alpha^{\tau}$, integrated against $\mu^{\tau}$, gives the other inequality.
\end{proof}

\begin{proposition}[\cite{RudolphInner82}*{Corollary~3.18}]\label{prop:intLipImpliesLinPin}
    If $\alpha^{\tau}$ is an integrably Lipschitz time change of $\alpha$, then it is a linear pinning time change of $\alpha$.
\end{proposition}

\medskip

The following theorem, which is a reformulation of a result by Rudolph in our current notation, establishes the fundamental connection between linear pinning and cohomologically tempered time changes:
\begin{theorem}[\cite{RudolphInner82}*{Theorem~3.4}]\label{thm:linearPinImpliesBarLin}
    If $\alpha^{\tau}$ is a linear pinning time change, then $\alpha^{\tau}$ is a cohomologically tempered time change.
\end{theorem}

\medskip

By definition, a cohomologically tempered time change $\alpha^{\tau}$ is cohomologous to a tempered time change of $\alpha$. A natural question arises:  can we impose stricter control on the properties of this target tempered time change? \cite{RudolphInner82}*{Corollary~3.17} answers this question affirmatively. We reformulate this corollary as follows:
\begin{proposition}[\cite{RudolphInner82}*{Corollary~3.17}]\label{prop:barelyLinTempControl}
    Suppose that $\alpha^{\tau}$ is a cohomologically tempered time change and $\epsilon>0$. Then $\alpha^{\tau}$ is cohomologous to a $C^{\infty}$-$\epsilon$-tempered time change $\alpha^{\bar{\tau}}$.
\end{proposition}

\medskip

Since the proof of Proposition~\ref{prop:barelyLinTempControl} is technical, we provide a sketch below for the reader's convenience. We begin with the following lemma, which reformulates results from \cite{RudolphInner82} in our current notation:
\begin{lemma}[\cite{RudolphInner82}*{Lemma~3.5}]\label{lem:rhoBoundedFinite}
    Let $\alpha^{\tau}$ be an integrable time change of $\alpha$. For any $T>0$, $\epsilon>0$, there exist a constant $B>0$ and a set $A\subset X$ with $\mu(A)>1-\epsilon$, such that for any $\v\in \R^n$ with $\norm{\v}<T$, if both $x$ and $\alpha_{\v}.x$ are in $A$, then
    \[
    \norm{\rho(x,\v)}<B.
    \]
\end{lemma}

\medskip

\begin{proof}[Sketch of the proof of Proposition~\ref{prop:barelyLinTempControl}]
Without loss of generality, we may assume that $\alpha^{\tau}$ is a tempered time change. Our goal is to prove that $\alpha^{\tau}$ is cohomologous to a time change $\alpha^{\bar{\tau}}$ satisfying the stricter control conditions. 

By Proposition~\ref{prop:jRhoProp}~\ref{item:JrhoErgodicSum}, for any $\epsilon'>0$, there exists $N_{\epsilon'}\in\N$ and a set $A_{\epsilon'}$ with $\mu(A_{\epsilon'})>1-\epsilon'$ such that for every $x\in A_{\epsilon'}$ and every $\norm{\v}\geq N_{\epsilon'}$
\[
\norm{\rho(x,\v)-J(\rho)(\v)}<\epsilon'\norm{\v}.
\]
Combining this with Lemma~\ref{lem:rhoBoundedFinite}, we deduce the existence of a positive constant $C'$ and a set $A'_{\epsilon'}$ with $\mu(A'_{\epsilon'})>1-2\epsilon'$ such that for every pair $x,\alpha_{\v}.x\in A'_{\epsilon'}$, we have
\begin{equation}\label{eq:barelyLinearJRho}
\norm{\rho(x,\v)-J(\rho)(\v)}<\epsilon'\norm{\v}+C'.    
\end{equation}

%We assume $\epsilon'$ is small enough such that $\sqrt{\epsilon'}<\norm{J(\rho)}$.

By Proposition~\ref{prop:jRhoProp}~\ref{item:JrhoDet}, $J(\rho)$ is nonsingular. We choose $\epsilon'<1/10$ sufficiently small below and write $\epsilon''=4\sqrt{n}\epsilon'.$

The remaining argument proceeds in two steps.
\begin{enumerate}[label=(\roman*)]
    \item \textbf{Construction of a measurable transfer function.} We first seek a measurable map $p:X\to\R^n$ such that $\alpha^{\tau}$ is cohomologous to a Lipschitz time change $\alpha^{\tilde{\tau}}$. That is, we require that for $\mu$-a.e. $x\in X$ and all $\v\in\R^n$
\[
\norm{\tilde{\rho}(x,\v)-J(\tilde{\rho})(\v)}\leq\epsilon''\norm{\v},
\]
where the new cocycle is given by
\[
\tilde{\rho}(x,\v)=\rho(x,\v)+p(x)-p(\alpha_{\v}.x).
\]
The parameter $\epsilon''$ is chosen so that $\epsilon''\norm{J(\rho)^{-1}}<1/2$. We summarize the construction of $p$ below and refer the reader to \cite{RudolphInner82}*{pp.~365--367} for details. 

The construction begins by eliminating the constant $C'$ in \eqref{eq:barelyLinearJRho} (cf.~\cite{RudolphInner82}*{p.~365}). The main idea is to use the strong Rokhlin lemma to find a transverse set with large self-return time, ensuring that $C'$ is negligible compared to the minimal self-return time under the action of $\alpha$. As in \cite{RudolphInner82}*{p.~365}, there is a Borel tower base $F$ such that
\[
V=\bigcup_{\v\in C_{C'/\epsilon'}}\alpha_{\v}(F),
\]
satisfies
\begin{enumerate}
    \item $\mu(V)>1-\epsilon'$;
    \item every $y\in V$ has a unique representation $y=\alpha_{\v}.x$ with $x\in F$ and $\v\in C_{C'/\epsilon'}$.
\end{enumerate}
For each Borel set $E\subset F$, define
\[
\eta_F(E)\coloneq\frac{1}{\mu(V)}\mu\left(\bigcup_{\v\in C_{C'/\epsilon'}}\alpha_{\v}(E)\right)
\]
The uniqueness of the tower coordinates implies that $\eta_F$ is a probability measure on $F$. The bass $F$ can be chosen so that 
\[
\eta_F(F\cap A'_{\epsilon'})>1-3\epsilon'.
\]
In these coordinates, $\mu|V$ is the product measure
\[
\frac{\mu(V)}{m(C_{C'/\epsilon'})}d\eta_F(x)dm(\v).
\]

If $x,\alpha_{\v}.x\in F\cap A'_{\epsilon'}$ and $\v\neq\0$, uniqueness implies $\norm{\v}\geq C'/\epsilon'$, and hence
\[
\norm{\rho(x,\v)-J(\rho)(\v)}<2\epsilon'\norm{\v}.
\]
The saturation of $F\cap A'_{\epsilon'}$ has full measure by ergodicity, since its tower has positive measure. We work on this invariant conull set.

The next step is an ``averaging'' procedure for the error term $\rho(x,\v)-J(\rho)(\v)$. Fix $x$ in this conull set and consider the set of times 
\[
G(x)=\{\v\in\R^n:\alpha_{\v}.x\in F\cap A'_{\epsilon'}\}.
\]
For each $\v\in G(x)$, let $C(x,\v)$ be the $n$-dimensional cube centered at $J(\rho)(\v)-\rho(x,\v)$ with side lengths $4\epsilon'\norm{\v}$ and edges parallel to the coordinate axes. \cite{RudolphInner82}*{Lemma~3.10} guarantees that the intersection $\bigcap_{\v\in G(x)}C(x,\v)$ is a nonempty closed rectangle, possibly degenerate, which we denote by $P(x)$.

Let $p(x)$ be the midpoint of $P(x)$. The map $x\mapsto p(x)$ is measurable. Furthermore, \cite{RudolphInner82}*{Lemma~3.11 and Corollary~3.12} give, on the invariant conull set, $x\in X$ and $\v\in\R^n$,
\begin{equation}\label{eq:lipschitzBound}
    \norm{p(x)-p(\alpha_{\v}.x)+\rho(x,\v)-J(\rho)(\v)}<\epsilon''\norm{\v},\qquad \v\neq\0.
\end{equation}
The cocycle identity gives the same Lipschitz bound for the difference $\tilde{\rho}_x-J(\rho)$ at any two points. Since $\epsilon''\norm{J(\rho)^{-1}}<1/2$, the contraction mapping theorem shows that $\tilde{\rho}_x$ is a bi-Lipschitz bijection of $\R^n$. This gives an integrable Lipschitz time change, cohomologous to $\rho$ and therefore $J(\tilde{\rho})=J(\rho)$.

\item \textbf{Smoothing via kernel function.} The second step is a regularization procedure, corresponding to \cite{RudolphInner82}*{Corollary~3.13--Lemma~3.16}. The idea is to use a kernel function to smooth the transfer function. Let $\kappa:\R^n\to\R^+$ be a $C^{\infty}$ integrating kernel of norm $1$ with bounded support. Define
\begin{equation*}
    \bar{p}(x)=\int_{\R^n}\kappa(\v)\cdot(p(\alpha_{\v}.x)-\rho(x,\v)+J(\rho)(\v))d\v.
\end{equation*}
It follows from \eqref{eq:lipschitzBound} that $\bar{p}(x)$ is well defined. We define the final cocycle $\bar{\rho}$ by
\[
\bar{\rho}(x,\v)=\rho(x,\v)+\bar{p}(x)-\bar{p}(\alpha_{\v}.x).
\]
It follows from \cite{RudolphInner82}*{Corollary~3.13} that on the invariant conull set and for every $\v\in\R^n$,
\begin{equation*}
    \norm{\bar{\rho}(x,\v)-J(\rho)(\v)}\leq\epsilon''\norm{\v}.
\end{equation*} 
Moreover, \cite{RudolphInner82}*{Lemma~3.14 and Lemma~3.15} imply that the maps $\v\mapsto\bar{\rho}(x,\v)$ and $\v\mapsto\bar{\tau}(x,\v)$ are $C^{\infty}$. The first derivative satisfies
\[
\norm{\partial_{\v}\bar{\rho}(x,\0)-J(\rho)}\leq\epsilon''.
\]
The same smallness condition gives global invertibility and 
\[
\norm{\partial_{\v}\bar{\tau}(x,\0)-J(\rho)^{-1}}\leq2\norm{J(\rho)^{-1}}^2\epsilon''.
\]
These uniform bounds make the time change integrable. By choosing $\epsilon''$ sufficiently small and notice that $J(\bar{\rho})=J(\rho)$ and $J(\bar{\tau})=J(\rho)^{-1}$, and   noting that $\alpha^{\bar{\tau}}$ is cohomologous to $\alpha^{\tau}$ (cf.~\cite{RudolphInner82}*{Lemma~3.16}), we finish the proof of Proposition~\ref{prop:barelyLinTempControl}.
\end{enumerate}
\end{proof}

Notice that linear pinning is preserved under cohomology: for sufficiently large $M$, restricting the pinning sets to where the transfer function has norm at most $M$ preserves positive measure and replaces the bound $ak+b$ by $ak+b+2M$. This together with Proposition~\ref{prop:temperedImpliesIntLip}, Proposition~\ref{prop:intLipImpliesLinPin}, Theorem~\ref{thm:linearPinImpliesBarLin}, and Proposition~\ref{prop:barelyLinTempControl} prove Theorem \ref{thm:kakutani-equivalences}.

\medskip

With the formal definition in place, we now address the invariance of entropy. A fundamental feature of Kakutani equivalence is that it preserves entropy classes (zero, positive, or infinite). This fact is a direct consequence of Nadler's generalization of the Abramov entropy formula. 

As discussed in the introduction, we provide a detailed proof of this result in Appendix~\ref{sec:AppendixAbramovFormula}, but we restate the formula here for convenience.
\begin{thm:AbramovMainNadler}[Nadler's Abramov formula \cite{NadlerAbramovNflow78}]
    Let $\alpha$ be a measurable, measure-preserving, ergodic, and free $\R^n$-action on a standard non-atomic probability space $(X,\mathcal{B},\mu)$, and $\alpha^{\tau}$ a tempered time change of $\alpha$. Then 
    \[
     h_{\mu}(\alpha)= h_{\mu^{\tau}}(\alpha^{\tau})\cdot\int_X\lambda_\tau(x)\,d\mu.
    \]
\end{thm:AbramovMainNadler}
This formula confirms that Kakutani equivalence for $\R^n$-actions also preserves the entropy class, as the integral factor is finite and non-zero.

\subsection{Compatibility with Kakutani equivalence for flows}
\label{sec:flow-compatibility}
Recall that standard Kakutani equivalence for flows requires only that the time change be integrable (cf.~Definition~\ref{def:integrableTimeChange}). We now show that in the one-dimensional case, the inherent ordering of $\R$ guarantees that any such time change is cohomologous to an integrably Lipschitz one. This guarantees that our generalized definition for $\R^n$-actions is compatible with the classical theory.

\begin{proposition}\label{prop:flowCompatiRn}
    Let $\alpha_t$ be an ergodic flow acting on $(X,\mathcal{B},\mu)$. If $\alpha_t^{\tau}$ is an integrable time change of $\alpha_t$, then there exists an integrably Lipschitz time change $\alpha_t^{\bar{\tau}}$ which is cohomologous to $\alpha_t^{\tau}$.
\end{proposition}
\begin{proof}
    It follows from the ergodicity that the orientation of the $\R$-action is constant for $\mu$-a.e.. Reversing time order in the target space if necessary, we may assume that the time change is increasing. It follows from the well-known Ambrose-Kakutani suspension flow representation theorem that every ergodic flow is measurably isomorphic to a suspension flow defined over an ergodic automorphism. Without loss of generality, we assume that $\alpha_t$ is the suspension flow $T_t^f$, where $T$ is an ergodic automorphism of a probability space $(Y,\mathcal{C},\nu)$ and $f\in L_+^1(Y,\nu)$ is the given roof function. 
    
    It follows from \cite{katok1977monotone}*{Proposition~2.3} that any flow Kakutani equivalent to $T_t^f$ is measurably isomorphic to a suspension flow over the same transformation $T$. By using the same cross section and the given orbit correspondence, we obtain that there exists a roof function $g\in L_+^1(Y,\nu)$ such that the suspension flow $T^g_t$ is measurably isomorphic to the time changed flow $\alpha_t^{\tau}$.

    Let $\epsilon=10^{-2}\min(\int_Yfd\nu,\int_Ygd\nu)$. By applying \cite{katok1977monotone}*{Proposition~3.4}, we can find functions $\bar{f}$ and $\bar{g}$ that are cohomologous to $f$ and $g$, respectively, over $T$, satisfying the uniform bounds
    \[
    \ess\sup_Y\abs{\bar{f}-\int_Yfd\nu}<\epsilon,\quad \text{and}\quad 
    \ess\sup_Y\abs{\bar{g}-\int_Ygd\nu}<\epsilon.
    \]
    Since $\epsilon$ is small relative to the integrals, both $\bar{f}$ and $\bar{g}$ are bounded away from zero and infinity.

    Define the map $\psi: Y^{\bar{f}}\to Y^{\bar{g}}$ between the suspension spaces by
    \[
    \psi(y,s)=\left(y,s\bar{g}(y)/\bar{f}(y)\right).
    \]
    This map $\psi$ is a nonsingular bijection that maps the orbits of $T_t^{\bar{f}}$ to orbits of $T_t^{\bar{g}}$. Let $\alpha^{\bar{\tau}}_t$ be the time change of $\alpha_t$ corresponding to $\psi$ (viewing $\alpha_t$ as isomorphic to $T_t^{\bar{f}}$). By construction, $\alpha^{\bar{\tau}}_t$ is cohomologous to $\alpha^{\tau}_t$. Furthermore, since the scaling factor $\bar{g}(y)/\bar{f}(y)$ is essentially bounded and bounded away from zero, the time change $\alpha^{\bar{\tau}}_t$ is bi-Lipschitz, and hence integrably Lipschitz. This completes the proof of Proposition~\ref{prop:flowCompatiRn}.
\end{proof}

\section{The Kakutani invariant for \texorpdfstring{$\R^n$}{Rn}-actions}\label{sec:KakutaniInvariant}
This section is dedicated to constructing a Kakutani invariant for $\R^n$-actions. Our approach builds upon Ratner's theory for $\R$-actions \cite{ratner1981some} and the $\bar{f}$-metric generalization by Feldman and Nadler \cite{FeldmanNadler79}. We begin by verifying fundamental properties such as invariance and the generating property. Subsequently, we derive a vanishing criterion for this invariant by adapting the seminal ``nesting arguments'' of \cites{ornstein1982equivalence,KatokSataevStandardIET}. In this context, we also acknowledge the foundational work of Johnson and \c{S}ahin \cites{JohnsonSahinRankOneLB98,JohnsonSahinIsoEx00}, who developed analogous nesting machinery for Kakutani equivalence in the setting of $\Z^n$-actions. The section concludes with a discussion on the compatibility of this generalized invariant with the classical one for flows.

\subsection{Hausdorff uniform structure and loosely Kronecker actions}\label{sec:fMetricMatching}
To formulate the generalization of the $\bar{f}$-metric for $\R^n$-actions, we explicitly define the necessary notation:
\begin{itemize}
    \item $D$ denotes a closed cell in $\mathbb{R}^n$ with volume $m(D)$. We denote the orbit segment over $D$ by 
    \[
    \alpha_D.x=\{\alpha_{\v}.x:\v\in D\};
    \]
    \item $\mathcal{P}$ is a finite measurable partition of $X$. For every $x\in X$ and every $\v\in \R^n$, the $(\mathcal{P},\alpha)$-name is defined by
\[
\mathcal{P}(x)(\v)=\mathcal{P}(\alpha_{\v}(x))=P\quad \text{if $\alpha_{\v}(x)\in P$,}
\]
where $P$ is an atom of the partition $\mathcal{P}$ containing the point;
    \item $\delta$ denotes the discrete metric on the partition atoms (i.e., $\delta(P,Q)=1$ if $P\neq Q$ and $0$ otherwise);
    \item $\mathcal{D}_D$ is the set of $C^{\infty}$ self-diffeomorphisms of $D$ that coincide with the identity in a neighborhood of the boundary $\partial D$.
\end{itemize}

We are now ready to define the generalization of the $\bar{f}$-metric:

\begin{definition}\label{def:epsilonPmatchable}
Let $x,y\in X$, $\epsilon>0$, and $D$ a closed cell in $\mathbb{R}^n$. The orbit segments $\alpha_D.x$ and $\alpha_D.y$ are said to be \textbf{$(\epsilon,\mathcal{P})$-matchable} if there exists $h\in\mathcal{D}_D$ such that 
\[
d_D\big(\mathcal{P}(x)\circ h,\mathcal{P}(y)\big)+\norm{h'-\id}_{\infty}<\epsilon,
\]
where $d_D$ is defined by
\[
d_D\big(\mathcal{P}(x)\circ h,\mathcal{P}(y)\big)=\frac{1}{m(D)}\int_D\delta\big(\mathcal{P}(x)( h(\v)),\mathcal{P}(y)(\v)\big)d\v.
\]
Such a map $h$ is referred to as an \textbf{$(\epsilon,\mathcal{P})$-match} for $\alpha_D.x$ and $\alpha_D.y$.
\end{definition}

\begin{definition}\label{def:fbarmetric}
We define the quantity
\[
f^{\alpha}_D(x,y,\mathcal{P})=\inf\{\epsilon>0:\text{$\alpha_D.x$ and $\alpha_D.y$ are $(\epsilon,\mathcal{P})$-matchable}\}.
\]
\end{definition}

As noted in \cite{FeldmanNadler79}, $f_{D}^{\alpha}(x,y,\mathcal{P})$ is not a metric but endows the space of $(D,\mathcal{P})$-names with a Hausdorff uniform structure. The following approximate triangle inequality, valid when $f_{D}^{\alpha}(y,z,\mathcal{P})<1$, will be sufficient for our purposes:
\begin{equation}\label{eq:approxTriIneq}
f_{D}^{\alpha}(x,z,\mathcal{P})\leq f_{D}^{\alpha}(y,z,\mathcal{P})+(1-f_{D}^{\alpha}(y,z,\mathcal{P}))^{-n}f_{D}^{\alpha}(x,y,\mathcal{P}).
\end{equation}
By composition and inversion of matching maps, for sufficiently small $\epsilon>0$ depending only on $n$, two points in the same $\epsilon$-ball, or the centers of two intersecting $\epsilon$-balls, have $f_D^{\alpha}$-distance less than $3\epsilon$ in either order.

\medskip

Recalling that $C_N=[-N,N]^n$, we define the property of being \textbf{loosely Kronecker} for $\R^n$-actions as follows:
\begin{definition}[\cite{FeldmanNadler79}]\label{def:LK}
A process $(\alpha,\mathcal{P})$ is said to be \textbf{loosely Kronecker} if for every $\epsilon>0$, there exists an $M>0$ such that for all $N>M$, there exists a set $E_N\subset X$ with $\mu(E_N)>1-\epsilon$ such that 
\[
f_{C_N}^{\alpha}(x,y,\mathcal{P})<\epsilon\text{ for all }x,y\in E_N.
\]
The $\R^n$-action $\alpha$ is called \textbf{loosely Kronecker} if the process $(\alpha,\mathcal{P})$ is loosely Kronecker for all finite measurable partitions of $X$.
\end{definition}

\subsection{Construction of Kakutani invariant for \texorpdfstring{$\R^n$}{Rn}-actions}\label{sec:defKakInvBasicProp}

In this section, we utilize the generalized $\bar{f}$-``metric'' (cf.~$f$-metric in Definition~\ref{def:fbarmetric}) to construct a Kakutani invariant for $\R^n$-actions. The construction proceeds analogously to Ratner's method for $\R$-actions \cite{ratner1981some}:
\begin{enumerate}[label=(\Roman*)]
\item \textbf{Covering Numbers:} Let $\mathcal{P}$ be a finite measurable partition of $X$, $\epsilon>0$, $R>1$ and $\mathcal{L}\in\GL_n(\R)$.
\begin{itemize}
    \item We define the ball of radius $\epsilon$ at $x\in X$ over the shape~$\mathcal{L}(C_R)$~by
    \[
    B^{\alpha}_{R,\mathcal{L}}(x,\epsilon,\mathcal{P})=\{y\in X: f^{\alpha}_{\mathcal{L}(C_R)}(x,y,\mathcal{P})<\epsilon\}.
    \]
    \item An \textbf{$(\epsilon,R,\mathcal{L},\mathcal{P})$-cover} of $X$, denoted by $\mathscr{C}^{\alpha}_{R,\mathcal{L}}(\epsilon,\mathcal{P})$, is a family of such balls whose union has measure greater than $1-\epsilon$.
    \item We define the covering number as 
    \[
K^{\alpha}_{R,\mathcal{L}}(\epsilon,\mathcal{P})=\inf\operatorname{Card}(\mathscr{C}^{\alpha}_{R,\mathcal{L}}(\epsilon,\mathcal{P})),
\]
where the infimum is taken over all $(\epsilon,R,\mathcal{L},\mathcal{P})$-covers of $X$.
\end{itemize}

\item\label{item:Fconstruction} \textbf{Scaling Family:} Let $\mathcal{F}$ denote the family of all nondecreasing functions $s:\mathbb{R}^+\to\mathbb{R}^+$ satisfying $\lim_{t\to+\infty}s(t)=+\infty$.
\item \textbf{The Invariant:} For $s\in\mathcal{F}$, we define the following quantities:
    \begin{equation*}% \label{eq:invariant}
    \begin{aligned}
    \beta(\alpha,s,\epsilon,\mathcal{L},\mathcal{P})&=\liminf_{R\to\infty}\frac{\log K^{\alpha}_{R,\mathcal{L}}(\epsilon,\mathcal{P})}{s(\abs{\det\mathcal{L}}R)};\\
    h_{\mu,s}^{\Kak}(\alpha,\mathcal{L},\mathcal{P})&=\limsup_{\epsilon\to0}\beta(\alpha,s,\epsilon,\mathcal{L},\mathcal{P});\\
    h_{\mu,s}^{\Kak}(\alpha,\mathcal{L})&=\sup_{\mathcal{P}}h_{\mu,s}^{\Kak}(\alpha,\mathcal{L},\mathcal{P}).
    \end{aligned}
    \end{equation*}
    The quantity $h_{\mu,s}^{\Kak}(\alpha,\mathcal{L})$ is referred to as the Kakutani invariant of $\alpha$ with respect to $s$ and $\mathcal{L}$.
\end{enumerate}

\medskip

The following proposition establishes conditions under which the dependence on the linear transformation $\mathcal{L}$ may be omitted.
\begin{proposition}\label{prop:LinearInvariance}
Let $\alpha:\R^n\curvearrowright(X,\mathcal{B},\mu)$ be an ergodic action. Suppose that $s\in\mathcal{F}$ satisfies
\[
\lim_{R\to\infty}\frac{s(aR)}{s(R)}=1\quad \text{for all $a>0$.}
\]
Then for any $\mathcal{L}_1,\mathcal{L}_2\in\GL_n(\R)$,
\[
h_{\mu,s}^{\Kak}(\alpha,\mathcal{L}_1)=h_{\mu,s}^{\Kak}(\alpha,\mathcal{L}_2).
\]
Furthermore, for any $\mathcal{L}_1,\mathcal{L}_2\in\GL_n(\R)$, there exist constants $C,A>1$ such that for every finite measurable partition $\mathcal{P}$, every $R>0$ and all sufficiently small $\epsilon>0$, 
\[
K_{R/C,\mathcal{L}_1}^{\alpha}(A\sqrt{\epsilon},\mathcal{P}) \leq K_{R,\mathcal{L}_2}^{\alpha}(\epsilon,\mathcal{P}) \leq K_{CR,\mathcal{L}_1}^{\alpha}(\epsilon^2/A^2,\mathcal{P}).
\]
\end{proposition}
\begin{proof}
It suffices to show that $h_{\mu,s}^{\Kak}(\alpha,\mathcal{L})=h_{\mu,s}^{\Kak}(\alpha,\id)$ for any $\mathcal{L}\in\GL_n(\R)$. We focus on establishing this equality.

Since $\mathcal{L}\in\GL_n(\R)$, there exists a constant $r_{\mathcal{L}}>1$ such that for every $R>0$, 
\begin{equation*}
C_{R/r_{\mathcal{L}}}\subset \mathcal{L}(C_R)\subset C_{Rr_{\mathcal{L}}}.
\end{equation*}
Let $\mathcal{P}$ be a finite measurable partition of $X$ and $\epsilon>0$. To restrict a matching map to either smaller cell, multiply its displacement from the identity by a smooth cutoff that vanishes near the boundary and equals one outside a relative $\sqrt{\epsilon}$ collar. Since the original displacement is $O(\epsilon R)$, the derivative error and the additional unmatched proportion are $O(\sqrt{\epsilon})$. For small $\epsilon$ this gives a matching self diffeomorphism, and hence, for a constant $A>1$ depending only on $\mathcal{L}$,
\begin{equation*}
K_{R/r_{\mathcal{L}},\id}^{\alpha}(A\sqrt{\epsilon},\mathcal{P})\leq K_{R,\mathcal{L}}^{\alpha}(\epsilon,\mathcal{P})\leq K_{Rr_{\mathcal{L}},\id}^{\alpha}(\epsilon^2/A^2,\mathcal{P}).
\end{equation*}
Dividing by $s(R)$ and adjusting the terms, we obtain
\begin{equation*}
\begin{aligned}
\frac{\log K_{R/r_{\mathcal{L}},\id}^{\alpha}(A\sqrt{\epsilon},\mathcal{P})}{s(R/r_{\mathcal{L}})}\cdot\frac{s(R/r_{\mathcal{L}})}{s(R)}&\leq \frac{\log K_{R,\mathcal{L}}^{\alpha}(\epsilon,\mathcal{P})}{s(\abs{\det\mathcal{L}}R)}\cdot\frac{s(\abs{\det\mathcal{L}}R)}{s(R)}\\
&\leq \frac{\log K_{Rr_{\mathcal{L}},\id}^{\alpha}(\epsilon^2/A^2,\mathcal{P})}{s(Rr_{\mathcal{L}})}\cdot\frac{s(Rr_{\mathcal{L}})}{s(R)}.
\end{aligned}
\end{equation*}
By the assumption that $\lim_{R\to\infty}\frac{s(aR)}{s(R)}=1$ for any $a>0$, the ratio $\frac{s(aR)}{s(R)}$ converges to $1$. Taking the limit inferior as $R\to\infty$, the inequality yields
\[\beta(\alpha,s,A\sqrt{\epsilon},\id,\mathcal{P})\leq\beta(\alpha,s,\epsilon,\mathcal{L},\mathcal{P})\leq \beta(\alpha,s,\epsilon^2/A^2,\id,\mathcal{P}).
\]
Letting $\epsilon\to 0$ establishes the equality. The same argument for $\mathcal{L}_1(C_R)$ and $\mathcal{L}_2(C_R)$ gives the last assertion.
\end{proof}
\begin{remark}
Under the assumption that $\lim_{R\to+\infty}\frac{s(aR)}{s(R)}=1$ for all $a>0$, Proposition~\ref{prop:LinearInvariance} allows us to remove the dependence on $\mathcal{L}$ and define the invariant simply by $h_{\mu,s}^{\Kak}(\alpha)$.
\end{remark}

\subsection{Fundamental properties of the Kakutani invariant}
In this subsection, we establish key properties of the invariant $h_{\mu,s}^{\Kak}(\alpha,\mathcal{L})$. More precisely, we will start with its invariance under higher rank Kakutani equivalence and then proceeding to the crucial generating properties.

\subsubsection{Invariance}
We now establish the invariance of $h_{\mu,s}^{\Kak}(\alpha,\mathcal{L})$ under higher rank Kakutani equivalence. This result generalizes seminal result \cite{ratner1981some}*{Theorem~3} to $\R^n$-actions.

\begin{proposition}\label{prop:KakutaniEquivalenceInvariant}
Let $\alpha$ and $\varphi$ be ergodic $\R^n$-actions on $(X,\mathcal{B},\mu)$ and $(Y,\mathcal{C},\nu)$, respectively. Suppose that $\alpha$ and $\varphi$ are Kakutani equivalent, then for any $\mathcal{L}_1,\mathcal{L}_2\in\GL_n(\R)$ and any $s\in \mathcal{F}$ satisfying 
\[
\lim_{R\to\infty}\frac{s(aR)}{s(R)}=1 \quad \text{for all $a>0$}
\]
we have
\[
h_{\mu,s}^{\Kak}(\alpha,\mathcal{L}_1)=h_{\nu,s}^{\Kak}(\varphi,\mathcal{L}_2).
\]
\end{proposition}

\medskip

The following lemma is a variation of \cite{FeldmanNadler79}*{Lemma~3.8} adapted to our setting.
\begin{lemma}\label{lem:matchingExtension}
Let $\mathcal{L}\in\GL_n(\R)$. For each sufficiently small $\epsilon>0$, there exists $\delta>0$ such that for any $C^{\infty}$ embedding
\[
h:(1-\delta)\mathcal{L}(C_N)\to \mathcal{L}(C_N)
\]
satisfying $\norm{h'-\id}_{\infty}<\delta$, there exists $h_1\in \mathcal{D}_{\mathcal{L}(C_N)}$ such that 
\begin{equation*}
\begin{gathered}
\norm{h'_1-\id}_{\infty}<\epsilon,\\
h_1(\v)=h(\v) \text{ if $\v\in(1-\epsilon)\mathcal{L}(C_N)$.}
\end{gathered}
\end{equation*}
\end{lemma}
\begin{proof}[Proof of Lemma~\ref{lem:matchingExtension}]
Fix $\epsilon\in(0,1)$ and choose $\delta>0$ sufficiently small. Let $h$ be a $C^{\infty}$ embedding from $(1-\delta)\mathcal{L}(C_N)$ to $\mathcal{L}(C_N)$ with $\norm{h'-\id}_{\infty}<\delta$. Consider the conjugated map 
\[
\tilde{h}=\mathcal{L}^{-1}\circ h\circ \mathcal{L}.
\]
This is a $C^{\infty}$ embedding from $(1-\delta)C_N$ to $C_N$ with 
\[
\norm{\tilde{h}'-\id}_{\infty}\leq \norm{\mathcal{L}^{-1}}\norm{h'-\id}_{\infty}\norm{\mathcal{L}}<\norm{\mathcal{L}}\norm{\mathcal{L}^{-1}}\delta.
\]
Define  $C'_{\mathcal{L}}=\norm{\mathcal{L}}\norm{\mathcal{L}^{-1}}$. Together with the derivative estimate, the inclusion $\tilde{h}((1-\delta)C_N)\subset C_N$ gives $\norm{\tilde{h}-\id}_{\infty}=O_{\mathcal{L}}(\delta N)$. Multiplying $\tilde{h}-\id$ by a smooth cutoff supported in $(1-\delta)C_N$ and equal to one on $(1-\epsilon)C_N$ gives, for sufficiently small $\delta$, a map $\tilde{h}_1\in\mathcal{D}_{C_N}$ with
\[
\norm{\tilde{h}'_1-\id}_{\infty}<\epsilon/C'_{\mathcal{L}},\qquad \tilde{h}_1(\v)=\tilde{h}(\v),\quad \v\in(1-\epsilon)C_N.
\]
Define 
\[
h_1=\mathcal{L}\circ\tilde{h}_1\circ\mathcal{L}^{-1}.
\]
Since $\mathcal{L}\in\GL_n(\R)$, it follows that  $h_1\in\mathcal{D}_{\mathcal{L}(C_N)}$, $\norm{h'_1-\id}_{\infty}<\epsilon$, and $h_1(\v)=h(\v)$ for all $\v\in(1-\epsilon)\mathcal{L}(C_N)$. This concludes the proof of Lemma~\ref{lem:matchingExtension}.
\end{proof}

\medskip

We now proceed to the proof of Proposition~\ref{prop:KakutaniEquivalenceInvariant}:
\begin{proof}[Proof of Proposition~\ref{prop:KakutaniEquivalenceInvariant}]
Let $\alpha$ and $\varphi$ be two ergodic $\R^n$-actions on $(X,\mathcal{B},\mu)$ and $(Y,\mathcal{C},\nu)$, respectively. Since $\alpha$ and $\varphi$ are Kakutani equivalent, there exists a tempered time change $\alpha^{\tilde{\tau}}$ of $\alpha$ that is measurably isomorphic to $\varphi$. 

Let $\mathcal{P}$ be a finite measurable partition of $Y$ and $\epsilon\in(0,1)$. Let $\mathcal{L}=J(\tilde{\rho})$ and choose $0<\eta<\epsilon/(10n)$. Apply Lemma~\ref{lem:matchingExtension} with collar $\eta$, and choose $\bar{\epsilon}>0$ sufficiently small in terms of $\epsilon$, $\mathcal{L}$ and $\eta$.

\medskip

By Proposition~\ref{prop:barelyLinTempControl}, there exists a $C^{\infty}$-$\bar{\epsilon}$-tempered time change $\alpha^{\tau}$ of $\alpha$ satisfying the following properties:
\begin{enumerate}[label=(\roman*)]
    \item $\alpha^{\tau}$ is cohomologous to $\alpha^{\tilde{\tau}}$;
    \item\label{item:psiConjugacy} $\alpha^{\tau}$ is measurably isomorphic to $\varphi$ via a measurable isomorphism $\psi$:
\[
\psi\circ\alpha_{\v}^{\tau}(x)=\varphi_{\v}\circ\psi(x),\qquad \text{for $\mu$-a.e. $x\in X$ and all $\v\in\R^n$.}
\]
\end{enumerate}

Let 
\[
\bar{\mathcal{P}}=\psi^{-1}(\mathcal{P}).
\]
Clearly, $\bar{\mathcal{P}}$ is a finite measurable partition of $X$.

\medskip

Let $\mathscr{C}$ be a $(\bar{\epsilon},R,\id,\bar{\mathcal{P}})$-cover of $X$ for the action $\alpha$. For any ball $B_{R,\id}^{\alpha}(x,\bar{\epsilon},\bar{\mathcal{P}})$ in $\mathscr{C}$ and any point $y$ in that ball, we have
\begin{equation*}
f^{\alpha}_{C_R}(x,y,\bar{\mathcal{P}})<\bar{\epsilon}.
\end{equation*}
Fix such points $x$ and $y$. Let $h\in\mathcal{D}_{C_R}$ be a matching function for $x,y$ such that
\begin{equation}\label{eq:hControl}
d_{C_R}(\bar{\mathcal{P}}(x)\circ h,\bar{\mathcal{P}}(y))+\norm{h'-\id}_{\infty}<2\bar{\epsilon}.
\end{equation}
We define the induced map $\tilde{h}:\rho(y,C_R)\to\rho(x,C_R)$ by
\[
\tilde{h}(\v)=\rho\Big(x,h\big(\tau(y,\v)\big)\Big),\qquad \v\in \rho(y,C_R).
\]
Since $\alpha^{\tau}$ is $\bar{\epsilon}$-tempered, combining this with \eqref{eq:hControl} yields, for a constant $C_{\mathcal{L}}$ depending only on $\mathcal{L}$,
\begin{equation}\label{eq:h1control}
\begin{aligned}
\tilde{h}((1-C_{\mathcal{L}}\bar{\epsilon})J(\rho)(C_R))\subset J(\rho)(C_R),\qquad 
\norm{\tilde{h}'-\id}_{\infty}<C_{\mathcal{L}}\bar{\epsilon}.
\end{aligned}
\end{equation}

Note that since $\alpha^{\tau}$ is cohomologous to $\alpha^{\tilde{\tau}}$, we have $J(\rho)=J(\tilde{\rho})$. Thus,  by Lemma~\ref{lem:matchingExtension},  equation~\eqref{eq:h1control}, and the choice of $\bar{\epsilon}$, there exists $\bar{h}\in\mathcal{D}_{J(\rho)(C_R)}$ satisfying
\begin{equation}\label{eq:h2Control}
\begin{gathered}
\norm{\bar{h}'-\id}_{\infty}<\epsilon/3,\\
\frac{1}{m(J(\rho)(C_R))}m\left(\{\v\in J(\rho)(C_R)\cap\rho(y,C_R):\tilde{h}(\v)=\bar{h}(\v)\}\right)>1-\epsilon/3.
\end{gathered}
\end{equation}

Let
\[
A_R\coloneq \{\w\in C_R:\bar{\mathcal{P}}(\alpha_{h(\w)}.x)=\bar{\mathcal{P}}(\alpha_{\w}.y)\},\qquad F_R=(1-\eta)J(\rho)(C_R)\cap\rho(y,A_R).
\]
By \eqref{eq:hControl}, the bound for $\rho(y,\cdot)$ and the choice of $\eta$. $\bar{\epsilon}$,
\[
m(F_R)>(1-\frac{\epsilon}{2})m(J(\rho)(C_R)).
\]
On $F_R$, the maps $\bar{h}$ and $\tilde{h}$ agree, and the changed names match. This implies, by the definition of the matching map $\bar{h}$, that
\begin{equation*}
f_{J(\rho)(C_R)}^{\alpha^{\tau}}(x,y,\bar{\mathcal{P}})<\epsilon.
\end{equation*}
By the property \ref{item:psiConjugacy} on p.~\pageref{item:psiConjugacy}, this transfers to the action $\varphi$:
\[
f_{J(\rho)(C_R)}^{\varphi}(\psi(x),\psi(y),\mathcal{P})<\epsilon.
\]

Since $d\mu^{\tau}/d\mu=\abs{\det\partial_{\v}\rho(x,\0)}/\abs{\det J(\rho)}<2$ for sufficiently small $\bar{\epsilon}$, we have $\mu^{\tau}(\cup\mathscr{C})>1-\epsilon$. Thus the images $\psi(x)$ and $\psi(y)$ belong to the same ball of an $(\epsilon,R,J(\rho),\mathcal{P})$-cover of $Y$ for the action $\varphi$. Denoting this cover by $\gamma$, we have
\[
\operatorname{Card}(\gamma)\leq\operatorname{Card}(\mathscr{C}).
\]
Consequently, taking the infimum over all covers $\mathscr{C}$ of $X$,
\begin{equation*}
K_{R,J(\rho)}^{\varphi}(\epsilon,\mathcal{P})\leq K_{R,\id}^{\alpha}(\bar{\epsilon},\bar{\mathcal{P}}).
\end{equation*}
Dividing by $s(R)$ and adjusting terms, this is equivalent to
\begin{equation*}
\frac{\log K_{R,J(\rho)}^{\varphi}(\epsilon,\mathcal{P})}{s(\abs{\det J(\rho)}R)}\cdot\frac{s(\abs{\det J(\rho)}R)}{s(R)}\leq \frac{\log K_{R,\id}^{\alpha}(\bar{\epsilon},\bar{\mathcal{P}})}{s(R)}.
\end{equation*}
Using the assumption $\lim_{R\to\infty}\frac{s(aR)}{s(R)}=1$ for all $a>0$, taking limits yields
\begin{equation*}
\beta(\varphi,s,\epsilon,J(\rho),\mathcal{P})\leq \beta(\alpha,s,\bar{\epsilon},\id,\bar{\mathcal{P}}).
\end{equation*}
Since $\mathcal{P}$ was arbitrary, we obtain
\begin{equation*}
h_{\nu,s}^{\Kak}(\varphi,J(\rho))\leq h_{\mu,s}^{\Kak}(\alpha,\id).
\end{equation*}

Finally, combining this inequality with Proposition~\ref{prop:LinearInvariance} ensures that for any $\mathcal{L}_1,\mathcal{L}_2\in\GL_n(\R)$
\begin{equation*}
h_{\nu,s}^{\Kak}(\varphi,\mathcal{L}_1)\leq h_{\mu,s}^{\Kak}(\alpha,\mathcal{L}_2).
\end{equation*}
The reverse inequality follows by symmetry, concluding the proof of Proposition~\ref{prop:KakutaniEquivalenceInvariant}.
\end{proof}

\subsubsection{Generating property}
We now turn to the generating property of the Kakutani invariant. The following result extends \cite{ratner1981some}*{Theorem~2} to the context of $\R^n$-actions.
\begin{proposition}\label{prop:generatingPartition}
Let $\alpha$ be an ergodic $\R^n$-action on $(X,\mathcal{B},\mu)$, and $\{\mathcal{P}_k\}_{k=1}^{\infty}$ an increasing sequence of finite measurable partitions of $X$ such that $\vee_{k=1}^{\infty}\mathcal{P}_k$ generates the $\sigma$-algebra $\mathcal{B}$. Then for any $s\in \mathcal{F}$ and any $\mathcal{L}\in\GL_n(\R)$,
\[
h_{\mu,s}^{\Kak}(\alpha,\mathcal{L})=\sup_kh_{\mu,s}^{\Kak}(\alpha,\mathcal{L},\mathcal{P}_k).
\]  
\end{proposition}
\begin{proof}[Proof of Proposition~\ref{prop:generatingPartition}]
Fix sufficiently small $\epsilon>0$ and a finite measurable partition  $\mathcal{P}=\{P_1,\ldots,P_r\}$.

Since $\vee_{k=1}^{\infty}\mathcal{P}_k$ generates the $\sigma$-algebra $\mathcal{B}$ and $\{\mathcal{P}_k\}_{k=1}^{\infty}$ is an increasing sequence, there exists $k\in\N$ and a finite measurable partition  $\mathcal{Q}=\{Q_1,\ldots,Q_r\}$ such that 
\begin{equation}\label{eq:partitionRefinement}
    \mathcal{Q}\leq\mathcal{P}_k \text{ and } \sum_{i=1}^r\mu(P_i\triangle Q_i)<\frac{\epsilon}{100\cdot 2^n}.
\end{equation}
Denote $B=\cup_{i=1}^r(P_i\triangle Q_i)$. Applying the Pointwise Ergodic Theorem (cf.~ \cite{LindenstraussPointwise01}*{Theorem~1.3}) to the $\R^n$-action $\alpha$,  the F\o lner set $\mathcal{L}(C_R)$ and the function $\chi_{B}$, we obtain that there exist $R_0>0$ and a set $Y\subset X$ with $\mu(Y)>1-\frac{\epsilon}{100\cdot 2^n}$ such that for all $R\geq R_0$ and $x\in Y$,
\begin{equation}\label{eq:ergodicPortionDiff}
   m(\{\v\in \mathcal{L}(C_R):\alpha_{\v}.x\in B\})<\frac{2\epsilon}{100\cdot 2^n}m(\mathcal{L}(C_R)).
\end{equation}

Fix an $R\geq R_0$ and let $\mathscr{C}$ be an $(\frac{\epsilon}{100\cdot 2^n},R,\mathcal{L},\mathcal{P}_k)$ cover of $X$. Define 
\[
Z=Y\cap\cup\mathscr{C}.
\]
It follows from the definitions of $Y$ and $\mathscr{C}$ that $\mu(Z)>1-\frac{2\epsilon}{100\cdot 2^n}$. Denote the restriction of the cover $\mathscr{C}$ to $Z$ by $\gamma$. It is clear that 
\[
\operatorname{Card}(\gamma)\leq\operatorname{Card}(\mathscr{C}).
\]

Suppose that $x$ and $y$ belong to the same element of $\gamma$. Then the estimates following \eqref{eq:approxTriIneq} gives
\[
f_{\mathcal{L}(C_R)}^{\alpha}(x,y,\mathcal{P}_k)\leq \frac{1}{100\cdot 2^n}\epsilon+2\cdot\frac{1}{100\cdot 2^n}\epsilon=\frac{3}{100\cdot 2^n}\epsilon.
\]
The above inequality implies that there exists a matching function $h$ between the orbit segments $\alpha_{\mathcal{L}(C_R)}.x$ and $\alpha_{\mathcal{L}(C_R)}.y$, and a set $A\subset  \mathcal{L}(C_R)$ with 
\[
m(A)>(1-\frac{3}{100\cdot 2^n}\epsilon)m(\mathcal{L}(C_R))
\]
so that 
\begin{equation*}
    \begin{gathered}
        \norm{h'-\id}_{\infty}<\frac{3}{100\cdot 2^n}\epsilon,\\
    \mathcal{P}_k(\alpha_{h(\v)}.x)=\mathcal{P}_k(\alpha_{\v}.y)\text{ for every $\v\in A$.}
    \end{gathered}
\end{equation*}
This, together with the definition of $B$, equations~\eqref{eq:partitionRefinement} and \eqref{eq:ergodicPortionDiff}, $x,y\in Z\subset Y$ and $\abs{\det{(h^{-1})'}}\leq 2^n$, gives a subset $A'\subset A$ with 
\[
m(A')>(1-\frac{(5+2^{n+1})\epsilon}{100\cdot 2^n})m(\mathcal{L}(C_R))
\]
such that 
\[
\mathcal{P}(\alpha_{h(\v)}.x)=\mathcal{P}(\alpha_{\v}.y) \text{ for any $\v\in A'$}.
\]
Since $\norm{h'-\id}_{\infty}<\frac{3}{100\cdot 2^n}\epsilon$, the above inequality yields
\[
f_{\mathcal{L}(C_R)}^{\alpha}(x,y,\mathcal{P})<\epsilon.
\]
Hence we only need $\operatorname{Card}(\gamma)$ many $B_{R,\mathcal{L}}^{\alpha}(x,\epsilon,\mathcal{P})$-balls to cover $Z$. Combining this with $\operatorname{Card}(\gamma)\leq\operatorname{Card}(\mathscr{C})$, we obtain
\[
K_{R,\mathcal{L}}^{\alpha}(\epsilon,\mathcal{P})\leq K_{R,\mathcal{L}}^{\alpha}\left(\frac{\epsilon}{100\cdot 2^n},\mathcal{P}_k\right).
\]
Dividing by $s(\abs{\det\mathcal{L}}R)$ and taking limits gives
\begin{equation*}
    \begin{aligned}
        \beta(\alpha,s,\epsilon,\mathcal{L},\mathcal{P})\leq\beta\left(\alpha,s,\frac{\epsilon}{100\cdot 2^n},\mathcal{L},\mathcal{P}_k\right)\leq \sup_kh_{\mu,s}^{\Kak}(\alpha,\mathcal{L},\mathcal{P}_k).
    \end{aligned}
\end{equation*}
Since $\epsilon$ and $\mathcal{P}$ were arbitrary, this implies
\[
h_{\mu,s}^{\Kak}(\alpha,\mathcal{L})\leq\sup_kh_{\mu,s}^{\Kak}(\alpha,\mathcal{L},\mathcal{P}_k).
\]
The reverse inequality follows immediately from the definition of $h_{\mu,s}^{\Kak}(\alpha,\mathcal{L})$, which concludes the proof of Proposition~\ref{prop:generatingPartition}.
\end{proof}

\subsection{Vanishing criterion for the Kakutani invariant}
The quantity $h_{\mu,s}^{\Kak}(\alpha,\mathcal{L})$ provides a precise characterization of loosely Kronecker systems. The following theorem establishes this connection, generalizing the result of Ratner \cite{ratner1981some}*{Theorem~1} to the setting of $\R^n$-actions. Recall that $\mathcal{S}$ is defined in \eqref{eq:SFamily}.
\begin{thm:vanishingCriterion}
Let $\alpha$ be a zero-entropy ergodic $\R^n$-action. Then $\alpha$ is loosely Kronecker if and only if there exists $\mathcal{L}\in\GL_n(\R)$ such that 
\[
h_{\mu,s}^{\Kak}(\alpha,\mathcal{L})=0\qquad  \text{for all $s\in \mathcal{S}$.}
\]
\end{thm:vanishingCriterion}

\subsubsection{Auxiliary lemmas}
We begin by establishing the following lemma, which ensures the existence of a basis with specific ergodicity properties. This result is necessary for the proof of  Theorem~\ref{thm:vanishingCriterion}.
\begin{lemma}\label{lemma:goodBasis}
Let $\alpha$ be an ergodic $\R^n$-action. Given a sequence of distinct positive real numbers $\{R_i\}_{i=1}^{\infty}$, there exists a basis $\mathcal{E}=\{\e_1,\ldots,\e_n\}$ of $\R^n$ such that for all $i\in\N$, $q\in\mathbb{Q}^+$, and $j=1,\ldots,n$, the transformation $\alpha_{qR_i\e_j}$ is ergodic. 
\end{lemma}
\begin{proof}[Proof of Lemma~\ref{lemma:goodBasis}]
By Theorem \ref{thm:PughShub}, the set of vectors in $\mathbb{R}^n$ such that the transformation $\alpha_{\v}$ is not ergodic is contained in a countable union of affine hyperplanes. In particular, this set, denoted by $\mathcal{R}$, has Lebesgue measure zero.  Define the set of ``bad'' vectors as
\[
\tilde{\mathcal{R}}=\bigcup_{i=1}^{\infty}\bigcup_{q\in\mathbb{Q}^+}\{\v\in\R^n:qR_i\v\in\mathcal{R}\}.
\]
Since each set in the union has measure zero, we obtain that $m(\tilde{\mathcal{R}})=0$. The complement $\R^n\setminus\tilde{\mathcal{R}}$ has full measure and therefore contains a basis for $\R^n$. This concludes the proof of Lemma~\ref{lemma:goodBasis}.
\end{proof}

\medskip

The following definition is required to adapt the ``nesting argument'' to the context of $\R^n$-actions. It relies on the notion of $(\epsilon,\mathcal{P})$-matchability introduced in  Definition~\ref{def:epsilonPmatchable}.
\begin{definition}\label{def:nestingDef}
Let $\mathcal{P}$ be a finite measurable partition of $X$. For any $x,y\in X$, we say that the orbit segments $\alpha_{C_R}.x$ and $\alpha_{C_R}.y$ are $(\theta,\epsilon,M,c)$-$\mathcal{P}$-matchable if there exist families of sub-cubes with disjoint interiors 
\begin{equation*}
C_i(x),C_i(y)\subset C_R,\text{ for }i=1,\ldots,k,
\end{equation*}
with paired cubes of equal side length, and an index set $J\subset\{1,\ldots,k\}$ such that
\begin{enumerate}[label=(\roman*)]
  \item \label{item:NestingCon1}$m\left(\bigcup_{i=1}^kC_i(x)\right),m\left(\bigcup_{i=1}^kC_i(y)\right)>(1-\theta)m(C_R)$;
  \item \label{item:NestingCon2}$m\left(\bigcup_{i\in J}C_i(x)\right),m\left(\bigcup_{i\in J}C_i(y)\right)\geq c\cdot m(C_R)$;
  \item \label{item:NestingCon3}for all $i\in J$, the orbit segments $\alpha_{C_i(x)}.x$ and $\alpha_{C_i(y)}.y$ are $(\epsilon,\mathcal{P})$-matchable;
  \item \label{item:NestingCon4}for all $i\notin J$, the volumes satisfy $m(C_i(x)),m(C_i(y))\geq2^nM^n$.
\end{enumerate}
The names on paired cubes are compared after translating their centers to the origin.
\end{definition}

We now state the ``nesting argument'' lemma.
\begin{lemma}\label{lem:nestingArgument}
Let $\mathcal{P}$ be a finite measurable partition of $X$, $K$ a positive integer and $\epsilon>0$ sufficiently small depending only on $n$. Suppose there exists an unbounded increasing sequence $\{R_i\}_{i\in\N}\subset\R^+$ such that:
\begin{enumerate}[label=(\alph*)]
    \item\label{item:NestingAssumptionA} $K_{R_i,\id}^{\alpha}(\epsilon,\mathcal{P})\leq K$ for all $i\in\N$;
    \item\label{item:NestingAssumptionB} $\alpha_{qR_i\e_j}$ is an ergodic transformation for every $q\in\Q\setminus\{0\}$, $i\in\N$ and $j=1,\ldots,n$. Here $\{\e_1,\ldots,\e_n\}$ is the standard basis of $\R^n$.
\end{enumerate} 
Then, for any $\theta>0$ and $M>0$, there exist $n_0=n_0(\epsilon,M,\theta)$ and a subset $Y\subset X$ with $\mu(Y)>1-\theta/2$ such that for all $R\geq n_0$ and all $x,y\in Y$, the orbit segments $\alpha_{C_R}.x$ and $\alpha_{C_R}.y$ are $\big(\theta,3\epsilon,M,1/(4K)\big)$-$\mathcal{P}$-matchable. Moreover, with density $1/(32K^2)$, the paired cubes can be chosen with a common radius $R_i>M$ and centers $\v+\w$ on the $x$-orbit and $\v$ on the $y$-orbit, where $\w$ belongs to a fixed finite set independent of $R,x,y$. The radius $R_i$ may be taken arbitrarily large.
\end{lemma}
\begin{proof}[Proof of Lemma~\ref{lem:nestingArgument}]
Let $\mathcal{P}$, $K$, $\epsilon$, and $\{R_i\}_{i\in\N}$ be fixed as in the assumption of the lemma.

Fix $\theta,M>0$. Choose an $R_i>M$. By assumption~\ref{item:NestingAssumptionA}, we have $K_{R_i,\id}^{\alpha}(\epsilon,\mathcal{P})\leq K$. Therefore, we can cover a portion of the space $X$ of measure at least $1-\epsilon$ using $K$ matching balls, which we denote by
\[
B_j^i=B_{R_i,\id}^{\alpha}(x_j,\epsilon,\mathcal{P}) \text{ for $j=1,\ldots,K$}.
\]
Since $\epsilon<1/2$, the union of these balls has measure greater than $1/2$. By the Pigeonhole Principle, there exists at least one matching ball, denoted by $B^i$, with
\[
\mu(B^{i})>\frac{1}{2K}.
\]
Choose $N_0>2R_i$ sufficiently large such that for all $R>N_0$, the boundary effect is small:
\[
\frac{(2R)^n-(2R-4R_i)^n}{(2R)^n}<10^{-2}\min\{\theta,10^{-2}\}.
\]

We apply the pointwise ergodic theorem to the action generated by $\{2R_i\e_1,\ldots,2R_i\e_n\}$ (which acts ergodically by assumption~\ref{item:NestingAssumptionB}) and the function $\chi_{B^i}$. We obtain a set $Y\subset X$ with $\mu(Y)>1-\theta/2$ and a constant $N_1>0$ such that for every $x\in Y$ and every $R>N_1$,
\begin{equation}\label{eq:latticeErgodicThm}
\card{\{\v\in C_{R-R_i}\cap(2R_i\Z)^n:\alpha_{\v}.x\in B^i\}} \geq\frac{3}{8K}\card{C_{R-R_i}\cap(2R_i\Z)^n},
\end{equation}
where $(2R_i\Z)^n$ denotes the lattice points $\{2R_i\v:\v\in\Z^n\}$.

Let $R\geq\max(N_0,N_1,M/\theta)$. For any $z\in Y$, define the set of ``good'' lattice points:
\begin{equation*}
    \begin{aligned}
     A_z=\{\v\in C_{R-R_i}\cap(2R_i\Z)^n:\alpha_{\v}.z\in B^i\}.
    \end{aligned}
\end{equation*}
Then for $x,y\in Y$, if $\v\in A_x$ and $\v'\in A_y$, then both orbit segments $\alpha_{C_{R_i}}.\alpha_{\v}.x$ and $\alpha_{C_{R_i}}.\alpha_{\v'}.y$ belong to the same matching ball $B^i$. By the approximate triangular inequality (cf.~\eqref{eq:approxTriIneq}), this implies:
\begin{equation}\label{eq:pairCon1}
\begin{aligned}
\alpha_{C_{R_i}}.\alpha_{\v}.x&\text{ and }\alpha_{C_{R_i}}.\alpha_{\v'}.y\text{ are }(3\epsilon,\mathcal{P})\text{-matchable}.
\end{aligned}
\end{equation}

We now construct the families of sub-cubes required by Definition~\ref{def:nestingDef}. Let the families $\{C_k(x)\}$ and $\{C_k(y)\}$ consist of the cubes $C_{R_i}(\v)$ for all $\v$ in the lattice $C_{R-R_i}\cap(2R_i\Z)^n$.

First, notice that these cubes fill $C_R$ except for a small neighborhood $\mathcal{N}$ of the boundary $\partial C_R$. By our choice of $R$,
\begin{equation}\label{eq:nestingThetaCheck}
\frac{m(\mathcal{N})}{m(C_R)}\leq\frac{(2R)^n-(2R-4R_i)^n}{(2R)^n}<\frac{\theta}{2}.
\end{equation}
This satisfies condition~\ref{item:NestingCon1} of Definition~\ref{def:nestingDef}.

Next, we identify the matching index set $J$. Let 
\[
J=\{1,\ldots,\min(\card{A_x},\card{A_y})\}.
\]
We can arbitrarily pair the cubes indexed by $A_x$ with the cubes indexed by $A_y$ (bijectively mapping the first $\card{J}$ elements). By \eqref{eq:pairCon1}, the paired orbit segments are $(3\epsilon,\mathcal{P})$-matchable, satisfying condition~\ref{item:NestingCon3} of Definition~\ref{def:nestingDef}.

From \eqref{eq:latticeErgodicThm}, we have the volume estimates for $z=x,y$:
\begin{equation*}
\begin{aligned}
\card{A_z}\cdot m(C_{R_i})\geq&\frac{3}{8K}\card{C_{R-R_i}\cap(2R_i\Z)^n}\cdot m(C_{R_i})\\
\geq&\frac{1}{4K}\cdot m(C_R).
\end{aligned}
\end{equation*}
Since this holds for both $x$ and $y$, the union of the paired cubes has measure at least $\frac{1}{4K}m(C_R)$. This verifies the condition~\ref{item:NestingCon2} of Definition~\ref{def:nestingDef}.

Finally, since $m(C_{R_i})=(2R_i)^n>(2M)^n$, condition \ref{item:NestingCon4} of Definition~\ref{def:nestingDef} is satisfied. Thus the orbit segments $\alpha_{C_R}.x$ and $\alpha_{C_R}.y$ are $(\theta,3\epsilon,M,\frac{1}{4K})$-$\mathcal{P}$-matchable. By picking $n_0=1+\max(N_0,N_1,M/\theta)$, we obtain the assertion.

\medskip

For the additional assertion, by the pointwise ergodic theorem, there are a finite set $W\subset(2R_i\Z)^n$ and a set $U$ with $\mu(U)>1-\mu(B^i)/8$ such that
\[
\frac{1}{\operatorname{Card}(W)}\sum_{\w\in W}\chi_{B^i}(\alpha_{\w}.z)>\mu(B^i)/2,\qquad z\in U.
\]
Let $\Lambda_R=(2R_i\Z)^n\cap C_{R-R_i-\sqrt{R}}$. By the pointwise ergodic theorem, there exist $R_0>R_i$ and a set $Y\subset X$ with 
$\mu(Y)>1-\theta/2$ such that, for every $R\geq R_0$ and $z\in Y$,
\[
\frac{\card{\{\v\in\Lambda_R:\alpha_{\v}.z\in B^i\}}}{\card{\Lambda_R}}>\frac{3\mu(B^i)}{4},\qquad \frac{\card{\{\v\in\Lambda_R:\alpha_{\v}.z\notin U\}}}{\card{\Lambda_R}}<\frac{\mu(B^i)}{4}.
\]
Thus, for any $x,y\in Y$, more than a proportion $\mu(B^i)/2$ of the centers $\v\in\Lambda_R$ satisfy both 
\[
\alpha_{\v}.y\in B^i,\qquad \alpha_{\v}.x\in U.
\]
Averaging over $W$, there exists $\w\in W$ such that more than a proportion $\mu(B^i)^2/4$ of the centers $\v\in\Lambda_R$ satisfy
\[
\alpha_{\v}.y\in B^i,\qquad \alpha_{\v+\w}.x\in B^i.
\]
For every $\v\in\Lambda_R$, pair the cube $C_{R_i}(\v+\w)$ on the $x$-orbit with $C_{R_i}(\v)$ on the $y$-orbit. Since $W$ is fixed and finite, all these cubes lie in $C_R$ for sufficiently large $R$. Each family has disjoint interiors, and 
\[
\frac{\card{\Lambda_R}m(C_{R_i})}{m(C_R)}\longrightarrow1.
\]
Thus each family covers more than $(1-\theta)m(C_R)$, and the union of the good paired cubes in each family has volume greater than
\[
\frac{\mu(B^i)^2}{8}m(C_R)\geq \frac{1}{32K^2}m(C_R).
\]
All cubes have radius $R_i>M$.
\end{proof}

The following generalizes \cite{KatokSataevStandardIET}*{Proposition~1} and \cite{ornstein1982equivalence}*{Proposition~7.4} to the setting of $\R^n$-actions. Analogous criteria for $\Z^n$-actions were established by Johnson and \c{S}ahin in \cite{JohnsonSahinRankOneLB98}*{Lemma~2.2.4} and \cite{JohnsonSahinIsoEx00}*{Theorem~3.4}.

\begin{proposition}
\label{prop:LK-Ferenczi}
    Let $\alpha$ be an ergodic $\R^n$-action. Suppose that for an increasing generating sequence of finite measurable partitions $\{\mathcal{P}_{\ell}\}_{\ell\in\N}$ and every $\epsilon> 0$, there exist a positive integer $K = K(\ell,\ve)$ and an unbounded increasing sequence $\{R_i\}_{i\in\N}\subset\R^+$ such that $K_{R_i,\id}^{\alpha}(\epsilon,\mathcal{P}_{\ell}) \leq K$ for all $i\in\N$. Then $\alpha$ is loosely Kronecker.
\end{proposition}

\begin{proof}[Proof of Proposition~\ref{prop:LK-Ferenczi}]
Fix $\ell\in\N$ and sufficiently small $\epsilon>0$. Let $\mathcal{P}_{\ell}$, $K$, and $\{R_i\}_{i\in\N}$ be defined as in the hypothesis. 

By Lemma~\ref{lemma:goodBasis}, there exists a good basis $\mathcal{E}=\{\e_1,\ldots,\e_n\}$ of $\R^n$ such that $\alpha_{qR_i\e_j}$ is an ergodic transformation for every $q\in\Q\setminus\{0\}$, $i\in\N$ and $j=1,\ldots,n$. The quantitative shape comparison in Proposition~\ref{prop:LinearInvariance} allows this change of coordinates, replacing $R_i$ by $R_i/C$ and $\epsilon$ by $A\sqrt{\epsilon}$, with a fixed integer $C>1$ and $A>0$ uniform for these bases. Relabeling the error, we may therefore assume that the good basis is the standard basis.

Define $c_0=0$ and 
\[
c_{k+1}=c_k+\frac{1}{32K^2}(1-c_k-6\epsilon),\qquad d_k=6\epsilon-3\cdot2^{-k}\epsilon.
\]
Then 
\[
c_k\nearrow1-6\epsilon,\qquad d_k\nearrow6\epsilon.
\]
We proceed by induction to show that, for every $\theta>0$, there exists $n_k(\theta)$ and $Y_k(\theta)$ with $\mu(Y_k(\theta))>1-\theta$ such that for all $R\geq n_k(\theta)$ and $x,y\in Y_k(\theta)$, there is $h\in \mathcal{D}_{C_R}$ satisfying
\[
\norm{h'-\id}_{\infty}<d_k,\qquad d_{C_R}(\mathcal{P}_{\ell}(x)\circ h,\mathcal{P}_{\ell}(y))\leq 1-c_k.
\]
\begin{enumerate}[label=(\roman*)]
    \item \textbf{Base Case:} For $k=0$, take $h=\id$ and $Y_0=X$.
    \item \textbf{Inductive Step:} Fix $\theta>0$ and let 
    \[
    \eta=\frac{3\epsilon}{256K^2},
    \]
    Apply the inductive hypothesis with measure loss $\eta$, using the set $Y_k(\eta)$ and threshold $n_k(\eta)$. Choose $R_i>n_k(\eta)$ in the additional assertion of Lemma~\ref{lem:nestingArgument} and let $W$ be the resulting finite set of translations. 
    
    Combining that lemma with the lattice ergodic theorem, we obtain a set $Y_{k+1}(\theta)$ with $\mu(Y_{k+1}(\theta))>1-\theta$ such that, for all sufficiently large $R$, all $x,y\in Y_{k+1}(\theta)$ and every $\w\in W$, each of the conditions
    \[
    \alpha_{\v}.y\notin Y_k(\eta),\qquad \alpha_{\v+\w}.x\notin Y_k(\eta)
    \]
    holds on less than a proportion $2\eta$ of the tile centers. Thus the cubes where at least one center point lies outside $Y_k(\eta)$ occupy less than a proportion $4\eta$ of $C_R$. Leave an outer collar region of width $\sqrt{R}$ and increase $R$ so that the untiled region's volume less than a proportion $4\eta$. Hence the total exceptional proportion is less than
    \[
    8\eta=\frac{3\epsilon}{32K^2}.
    \]

    Recall that the lemma gives one common translation vector $\w$ and good cubes occupying a positive proportion $t\geq 1/(32K^2)$. On these cubes, we can use the matches of cost less than $3\epsilon$. On the remaining cubes, we just use the inductive matches whenever both center points lie in $Y_k(\eta)$, and finally we use identity otherwise. In particular, the centered maps fix their cube boundaries and also glue to $H\in\mathcal{D}_{C_R}$ with 
    \[
    \norm{H'-\id}_{\infty}<d_k.
    \]

    Let $\chi$ be a smooth bump function which equals to one on the tiled cubes, vanishes near $\partial C_R$, and satisfying
    \[
    \norm{\nabla\chi}_{\infty}=O(R^{-1/2}).
    \]
    Recall that $\w$ belongs to a fixed finite set, the derivative error term of $g(\v)=\v+\chi(\v)\w$ goes to zero as $R\to\infty$. If $R$ is sufficiently large, we obtain that the map $h=g\circ H$ belongs to $\mathcal{D}_{C_R}$, has derivative error term less than $d_{k+1}$, and also inserts the required common translation vector on each desired cube. 

    Let $\kappa$ denote the sum of the proportions of the untiled region, and the exceptional cubes outside the good cubes, then we have 
    \[
    \kappa<3\epsilon/(32K^2).
    \]
    Recall that the good cubes contributes at most $3\epsilon t$ to the mismatch fraction. Moreover, the remaining non-exceptional cubes contribute at most $(1-t-\kappa)(1-c_k)$ proportion, and the exceptional region contributes at most $\kappa$ proportion. Hence
    \[
    \begin{aligned}
        d_{C_R}\left(\mathcal{P}_{\ell}(x)\circ h,\mathcal{P}_{\ell}(y)\right)&\leq 3\epsilon t+(1-t-\kappa)(1-c_k)+\kappa\\
        &\leq 3\epsilon t+(1-t)(1-c_k)+\frac{3\epsilon}{32K^2}\\
        &\leq \left(1-\frac{1}{32K^2}\right)(1-c_k)+\frac{3\epsilon}{16K^2}=1-c_{k+1},
    \end{aligned}
    \]
    where the last inequality follows from the facts that $t\geq 1/(32K^2)$ and $1-c_k\geq 6\epsilon$. By taking $n_{k+1}(\theta)$ sufficiently large, we complete the inductive step.
\end{enumerate}

\medskip

We choose $k_0$ so that $1-c_{k_0}<6\epsilon+\theta$. The resulting maps give $(12\epsilon+\theta,\mathcal{P}_{\ell})$-matches for all sufficiently large $R$ on a set of measure greater than $1-\theta$. The same shape comparison transfers these matches back to the original coordinates, with errors tending to zero as $\epsilon,\theta\to0$.

\medskip

Since $\{\mathcal{P}_{\ell}\}_{\ell\in\N}$ is increasing and generating, the partition approximation argument in the proof of Proposition~\ref{prop:generatingPartition} gives the same conclusion for every finite measurable partition. Hence $\alpha$ is loosely Kronecker.
\end{proof}

The following corollary of Proposition~\ref{prop:LK-Ferenczi} will be used in \S\ref{sec:horocyclicAction}:
\begin{corollary}\label{cor:LK-Ferenczi}
    Let $\alpha$ be an ergodic $\R^n$-action. Then $\alpha$ is loosely Kronecker if and only if there exists an increasing generating sequence of finite measurable partitions $\{\mathcal{P}_{\ell}\}_{\ell\in\N}$ such that for every $\ell\in\N$ and every $\epsilon > 0$, there exist a constant $c = c(\epsilon,\ell)>0$, an increasing sequence $\{R_i\}_{i\in\N}$ with $\lim_{i\to\infty}R_i=+\infty$ and sets $E_{R_i}\subset X$ with $\mu(E_{R_i})>1-\epsilon$ such that for every $x\in E_{R_i}$,
    \[ 
    \mu(B^{\alpha}_{R_i,\id}(x,\ve,\mc P_{\ell})) > c.
    \]
\end{corollary}
\begin{proof}[Proof of Corollary~\ref{cor:LK-Ferenczi}]
    Suppose that $\alpha$ is loosely Kronecker. It follows from Definition~\ref{def:LK} that for every finite measurable partition $\mathcal{P}$ and every $\epsilon>0$, there exists an $M>0$ such that for all $N>M$, there exists $E_N\subset X$ with $\mu(E_N)>1-\epsilon$ such that for every $x\in E_N$,
    \begin{equation*}
        \mu(B^{\alpha}_{N,\id}(x,\epsilon,\mathcal{P}))>1-\epsilon.
    \end{equation*}
    This finishes the proof for one direction.

    \medskip

    We now proceed to the other direction. Suppose the assumptions in the statement of the corollary hold. Fix an $\ell\in\N$ and sufficiently small $\epsilon>0$ depending only on $n$. Let $c$, $\{R_i\}_{i\in\N}$, and $\{E_{R_i}\}_{i\in\N}$ be the constants and sets provided by the hypothesis. Fix an index $i\in\N$. Let $\Lambda\subset E_{R_i}$ be a maximal collection of points such that if $x,y\in\Lambda$ and $x\neq y$, then 
    \begin{equation}\label{eq:xyfRSeparated}
        \max\left(f_{C_{R_i}}^{\alpha}(x,y,\mathcal{P}_{\ell}),f_{C_{R_i}}^{\alpha}(y,x,\mathcal{P}_{\ell})\right)\geq5\epsilon.
    \end{equation}
    We claim that for any distinct $x,y\in\Lambda$, 
    \[
    B^{\alpha}_{R_i,\id}(x,\epsilon,\mathcal{P}_{\ell})\cap B^{\alpha}_{R_i,\id}(y,\epsilon,\mathcal{P}_{\ell})=\emptyset.
    \]
    Suppose this is not the case and there exists $z\in B^{\alpha}_{R_i,\id}(x,\epsilon,\mathcal{P}_{\ell})\cap B^{\alpha}_{R_i,\id}(y,\epsilon,\mathcal{P}_{\ell})$. Then the estimates following \eqref{eq:approxTriIneq} give
    \[
    \max\left(f_{C_{R_i}}^{\alpha}(x,y,\mathcal{P}_{\ell}), f_{C_{R_i}}^{\alpha}(y,x,\mathcal{P}_{\ell})\right)\leq3\epsilon,
    \]
    which contradicts the assumption \eqref{eq:xyfRSeparated} and hence proves the claim. 

    Because $x\in E_{R_i}$ for all $x\in\Lambda$, the hypothesis guarantees that each of these disjoint $\epsilon$-$f$ balls has measure strictly greater than $c$. Since $\mu(X)=1$, it follows that $\card{\Lambda}\leq \frac{1}{c}$. Furthermore, by the maximality of $\Lambda$, the collection of balls $\{B_{R_i,\id}^{\alpha}(x,5\epsilon,\mathcal{P}_{\ell})\}_{x\in\Lambda}$ 
    must cover $E_{R_i}$. Since $\mu(E_{R_i})>1-\epsilon$, this forms a valid cover for the definition of the covering number and gives that $K_{R_i,\id}^{\alpha}(5\epsilon,\mathcal{P}_{\ell})\leq 1/c$. This, together with Proposition~\ref{prop:LK-Ferenczi}, finishes the proof of Corollary~\ref{cor:LK-Ferenczi}.
\end{proof}

\subsubsection{Proof of Theorem~\ref{thm:vanishingCriterion}}

We now return to the proof of Theorem~\ref{thm:vanishingCriterion}.

\medskip

Suppose that $\alpha$ is loosely Kronecker. By Definition \ref{def:LK}, for any $\epsilon>0$ and any finite measurable partition $\mathcal{P}$, for sufficiently large $R$, the covering number satisfies
    \[
    K_{R,\id}^{\alpha}(\epsilon,\mathcal{P})=1.
    \]
    This, together with Proposition~\ref{prop:LinearInvariance}, implies that for all $s\in\mathcal{F}$ and every $\mathcal{L}\in\GL_n(\R)$, we have
    \[
    \beta(\alpha,s,\epsilon,\mathcal{L},\mathcal{P})=0.
    \]
    By taking $\epsilon\to0$ and the supremum over $\mathcal{P}$, we obtain
    \[
    h_{\mu,s}^{\Kak}(\alpha,\mathcal{L})=0 \quad \text{for all $s\in\mathcal{F}$ and $\mathcal{L}\in\GL_n(\R)$,} 
    \]
    which concludes the proof of one direction of Theorem~\ref{thm:vanishingCriterion}.

\medskip

Conversely, suppose that there exists an $\mathcal{L}\in\operatorname{GL}(n,\mathbb{R})$ such that 
\[
h_{\mu,s}^{\Kak}(\alpha,\mathcal{L})=0 \quad \text{for all $s\in\mathcal{S}$}.
\]
By the definition of the invariant (cf.~\S\ref{sec:defKakInvBasicProp}), this implies that for all $s\in\mathcal{S}$, any $\epsilon>0$ and any finite measurable partition $\mathcal{P}$, we have 
\[
\beta(\alpha,s,\epsilon,\mathcal{L},\mathcal{P})=0.
\]
In particular, fix an arbitrary finite measurable partition $\mathcal{P}$ and $\epsilon>0$. We claim:
\begin{enumerate}[label=(\Alph*)]
    \item\label{item:boundedClaim} there exists an unbounded nondecreasing sequence $\{R_i\}_{i\in\N}$ such that $K_{R_i,\mathcal{L}}^{\alpha}(\epsilon,\mathcal{P})$ is bounded by a constant for every $i\in\N$. 
\end{enumerate}
    
Suppose this is not the case. Then
\[
\tilde{s}(R)=\inf\{\log K_{r,\mathcal{L}}^{\alpha}(\epsilon,\mathcal{P}):r\geq R\}
\]
is a nondecreasing function satisfying $\lim_{R\to+\infty}\tilde{s}(R)=+\infty$. Choose a positive, smooth, nondecreasing function $u:\R\to(0,\infty)$ tending to infinity sufficiently slowly that
\[
0\leq u'(t)\leq 1,\qquad u(t)\leq\tilde{s}(e^t),\quad \text{for all sufficiently large $t$.} 
\]
Define $s(R)=u(\log R)$. For every fixed $a>0$,
\[
\abs{\frac{s(aR)}{s(R)}-1}\leq\frac{\abs{\log a}}{s(R)}\longrightarrow0.
\]
Thus $s\in\mathcal{S}$. Moreover, for all sufficiently large $R$,
\[
s(R)\leq\tilde{s}(R)\leq\log K_{R,\mathcal{L}}^{\alpha}(\epsilon,\mathcal{P}).
\]
Consequently,
\[
\beta(\alpha,s,\epsilon,\mathcal{L},\mathcal{P})\geq\liminf_{R\to+\infty}\frac{s(R)}{s(\abs{\det\mathcal{L}}R)}=1,
\]
a contradiction. Thus, Claim~\ref{item:boundedClaim} is proved.

Recall that Proposition~\ref{prop:LinearInvariance} and Claim~\ref{item:boundedClaim} imply that there exists an unbounded nondecreasing sequence $\{\tilde{R}_i\}_{i\in\N}$ along which the covering numbers for $\id$ are bounded. Since this holds for an arbitrary partition $\mathcal{P}$, it holds in particular for every partition in a generating sequence $\{\mathcal{P}_{\ell}\}_{\ell\in\N}$. Thus, the hypothesis of Proposition~\ref{prop:LK-Ferenczi} is satisfied; hence $\alpha$ is loosely Kronecker. This completes the proof of the remaining direction of Theorem~\ref{thm:vanishingCriterion}.\qed

\subsection{Compatibility with the Kakutani invariant for flows}
Recall that in Ratner's construction of the Kakutani invariant for flows \cites{ratner1981some,ratner2017time}, the requirements imposed on the matching function are less restrictive than those in \S\ref{sec:fMetricMatching}. In this section, we demonstrate that these two definitions yield the same invariant for flows.

We begin by recalling the relevant definitions from \cites{ratner1981some,ratner2017time}:
\begin{definition}
    Let $\alpha$ be an ergodic flow acting on $(X,\mathcal{B},\mu)$ and $\mathcal{P}$ a finite measurable partition. For any $x,y\in X$ and $\epsilon, R>0$, the orbit segments $\alpha_{[0,R]}.x$ and $\alpha_{[0,R]}.y$ are said to be \textbf{one-dimensional $(\epsilon,\mathcal{P})$-matchable} if there exist a subset $A\subset[0,R]$ with $m(A)>(1-\epsilon)R$, and an increasing absolutely continuous map $h:A\to A'\subset[0,R]$ onto its image $A'$, where $m(A')\geq(1-\epsilon)R$, such that $\mathcal{P}(\alpha_t.x)=\mathcal{P}(\alpha_{h(t)}.y)$ for all $t\in A$, and the derivative satisfies
    \[
    \abs{h'(t)-1}<\epsilon\quad \text{for all $t\in A$.}
    \]
\end{definition}

\medskip

Recall that by the Ambrose-Kakutani representation theorem, the ergodic flow $\alpha$ is measurably isomorphic to a suspension flow. This implies the existence of an increasing generating sequence of partitions $\{\mathcal{P}_n\}_{n\in\N}$ such that for every $n\in\N$, there exist a sequence of positive numbers $\{\epsilon_k\}_{k\in\N}$ decreasing to zero and a sequence of measurable sets $D_k\subset X$ with $\mu(D_k)>1-10^{-k}$ satisfying the following: for every $x\in D_k$,
\begin{equation}\label{eq:goodSetPartition}
    \mathcal{P}_n(x)=\mathcal{P}_n(\alpha_t.x)\quad \text{ for all $\abs{t}\leq\epsilon_k$.}
\end{equation}
We fix $\alpha$ and the sequence $\{\mathcal{P}_n\}_{n\in\N}$ as above for the remainder of this section. The coincidence of the two invariants follows from the proposition below.
\begin{proposition}\label{prop:oneDimMatchImpliesMatch}
Fix $n\in\N$ and $\epsilon>0$. There exists a set $Z\subset X$ with $\mu(Z)>1-\epsilon$ such that for any $x,y\in Z$ and all sufficiently large $R$, the following implication holds: if the orbit segments $\alpha_{[0,R]}.x$ and $\alpha_{[0,R]}.y$ are one-dimensional $(\epsilon,\mathcal{P}_n)$-matchable, then they are $(20\sqrt{\epsilon},\mathcal{P}_n)$-matchable.
\end{proposition}

It follows from the stability estimates in  Proposition~\ref{prop:oneDimMatchImpliesMatch} and the reduction to generating partitions in  Proposition~\ref{prop:generatingPartition} that the Kakutani invariant defined in \S\ref{sec:KakutaniInvariant} is bounded above by that defined in \cites{ratner1981some,ratner2017time}. Since the reverse inequality is immediate from the definitions, we establish the coincidence of these two invariants, using the same interval and scale normalization. For $s\in\mathcal{S}$, normalization by interval length or radius gives the same value since $s(2R)/s(R)\to1$.

\medskip

We now proceed to the proof of Proposition~\ref{prop:oneDimMatchImpliesMatch}. 
The following combinatorial lemma plays a central role in the argument, allowing us to merge small error intervals into controlled safety gaps.
\begin{lemma}\label{lem:Qconstruction}
Fix $R>1$ and let $\epsilon > 0$ be sufficiently small.
Let $\{O_i\}_{i \in I}$ be a collection of disjoint open intervals of $[0,R]$, and let $\{\ell_i\}_{i\in I}$ be a family of positive numbers satisfying  $\ell_i \geq m(O_i)$ for all $i\in I$ and 
\[
\sum_{i \in I}\ell_i < 5\epsilon R.
\]
Then there exist a collection $\{Q_j\}_{j \in J}$ of disjoint open intervals of $[0,R]$ and a map $f : I \to J$ such that:
\begin{enumerate}[label=(\alph*)]
\item\label{part:covers} $O_i \subset Q_{f(i)}$ for all $i\in I$;
\item\label{part:length} $m(Q_j) = \frac{1}{\sqrt{\epsilon}}\sum_{f(i) = j} \ell_i$ for all $j\in J$.
\end{enumerate}
\end{lemma}

\begin{proof}[Proof of Lemma~\ref{lem:Qconstruction}]
Since the index set $I$ is countable, we may index the collection as $\{O_n\}_{n=1}^N$, where $N\in\N\cup\{\infty\}$. We proceed by induction to construct a sequence of families of disjoint intervals.

\begin{enumerate}[label=(\roman*)]
    \item We start with the base case $n=1$. Choose an open interval $Q_1^{(1)}$ containing $O_1$ with length $m(Q_1^{(1)})=\ell_1/\sqrt{\epsilon}$. This is possible since $\ell_1/\sqrt{\epsilon}<5\sqrt{\epsilon} R<R$.
    \item Suppose that for $n-1$, we have already constructed a collection $\mathcal{Q}_{n-1}$ of disjoint open intervals such that each $Q\in\mathcal{Q}_{n-1}$ satisfies \ref{part:length} for the subcollection of sets $\{O_k\}_{k=1}^{n-1}$ it covers. Let
    \[
    \mathcal{K}_0=\{Q\in\mathcal{Q}_{n-1}:Q\cap O_n\neq\emptyset\}
    \]
    be the collection of the interval in the current cover which intersects $O_n$.
    \begin{enumerate}
        \item If $\mathcal{K}_0=\emptyset$, we just choose a new interval $Q_{\operatorname{new}}$ that contains $O_n$ with length $\ell_n/\sqrt{\epsilon}$. If this interval is disjoint from $\mathcal{Q}_{n-1}$, we then set 
        \[
        \mathcal{Q}_{n}=\mathcal{Q}_{n-1}\cup\{Q_{\operatorname{new}}\}.
        \]
        \item If $\mathcal{K}_0\neq\emptyset$, we must merge. By defining a larger cluster $\mathcal{K}\subset\mathcal{Q}_{n-1}$ recursively to include all intervals related to $O_n$. However, recall that $\mathcal{Q}_{n-1}$ is a finite collection of disjoint intervals, we can define $W$ as the connected component of the set $O_n\cup\bigcup_{Q\in\mathcal{Q}_{n-1}}Q$ containing $O_n$. Let 
        \[
        \mathcal{K}=\{Q\in\mathcal{Q}_{n-1}:Q\subset W\}.
        \]
        Construct a single new interval $Q_{\operatorname{new}}$ containing $W$ such that 
        \[
        m(Q_{\operatorname{new}})=\frac{\ell_n}{\sqrt{\epsilon}}+\sum_{Q\in\mathcal{K}}m(Q).
        \]
        Notice that the sets in $\mathcal{K}$ are disjoint, we obtain that the measure of their union is the sum of their measures. Moreover, it follows from $O_n$ connects with them that we have strict sub-additivity for the geometric union
        \[
        m(W)<m(\bigcup\mathcal{K})+m(O_n)\leq m(Q_{\operatorname{new}}).
        \]
        This in particular guarantees the existence of $Q_{\operatorname{new}}$. We thus define the new family as 
        \[
        \mathcal{Q}_n=(\mathcal{Q}_{n-1}\setminus\mathcal{K})\cup\{Q_{\operatorname{new}}\}.
        \]
    \end{enumerate}

    In either case, if the enlarged interval meets another old interval, we absorb it, add its length to the prescribed length, and enlarge again. Notice that there are finitely many old intervals and the total prescribed length is bounded by $R$, this process will terminate with a disjoint family in finite steps. This completes the inductive step. The total length of the system is bounded by 
    \[
    \frac{1}{\sqrt{\epsilon}}\sum_{i\in I}\ell_i<5\sqrt{\epsilon} R<R,
    \]
    ensuring we stay within $[0,R]$. This completes the proof for the case where $I$ is a finite set.

    \item We now address the situation $I=\N$. Observe that the construction is monotonic: for every $n$, each interval $Q\in\mathcal{Q}_n$ is contained in a unique interval $Q'\in\mathcal{Q}_{n+1}$. For each $O_k$, let $Q^{(n)}(O_k)$ denote the unique interval in $\mathcal{Q}_n$ containing $O_k$ (defined for $n\geq k$). This sequence is increasing, i.e., $Q^{(k)}(O_k)\subset Q^{(k+1)}(O_k)\subset\ldots$.

    We define the limit collection $\{Q_j\}_{j\in J}$ inductively. Let
    \[
    Q_1=\bigcup_{n=1}^{\infty}Q^{(n)}(O_1)\text{ and }
    A_1 = \set{ i \in \N : O_i \subset Q_1}.
    \]
    Since $Q^{(n)}(O_1)$ are open intervals, their increasing union $Q_1$ is an open interval, whose measure satisfies the required sum condition by convergence.

    Next let $i_2$ be the smallest element of $\N \setminus A_1$. Define 
    \[
    Q_2=\bigcup_{n=i_2}^{\infty}Q^{(n)}(O_{i_2})\text{ and }
    A_2 = \set{ i \in \N : O_i \subset Q_2}.
    \]
    We repeat this process inductively to define the collection $\{Q_j\}_{j\in J}$ and the partition sets $\{A_j\}_{j\in J}$. 
    
    Suppose for contradiction that $Q_j\cap Q_{j'}\neq\emptyset$ for $j\neq j'$. Since $Q_j$ and $Q_{j'}$ are unions of open sets, there exists a point $x$ and an integer $N$ such that $x\in Q^{(N)}(O_{i_j})$ and $x\in Q^{(N)}(O_{i_{j'}})$. But by construction, intersecting intervals in $\mathcal{Q}_N$ are merged into a single interval and thus $Q_{j}=Q_{j'}$, a contradiction. 

    Thus the limiting intervals $\{Q_j\}_{j\in J}$ are disjoint and the map $f:I\to J$ defined by $f(i) = j$ if $i \in A_j$ satisfies the required condition (cf.~\ref{part:length}).
\end{enumerate}
\end{proof}

\medskip

Equipped with Lemma~\ref{lem:Qconstruction}, we are now ready to prove Proposition~\ref{prop:oneDimMatchImpliesMatch}.
\begin{proof}[Proof of Proposition \ref{prop:oneDimMatchImpliesMatch}]
Fix $n\in\N$. We may assume $20\sqrt{\epsilon}\leq1$, since otherwise the identity gives the conclusion. Choose $k$ sufficiently large such that 
\[
\max(10^{-k+10},\epsilon_k)<\epsilon
\]
where $\epsilon_k$ is the sequence associated with $\mathcal{P}_n$ as in \eqref{eq:goodSetPartition}. Let $D_k$ be the corresponding good set with $\mu(D_k)>1-10^{-k}$ as in \eqref{eq:goodSetPartition}. By the Pointwise Ergodic Theorem and Egorov's theorem, there exists a set $Z$ with $\mu(Z)>1-\epsilon/10$ and a positive integer $N_1$ such that for every $R\geq N_1$ and every $x\in Z$, 
\begin{equation}\label{eq:ergodicErrorEstimates}
m(\{t\in[0,R]:\alpha_t.x\in D_k\})\geq (1-\epsilon/10) R.
\end{equation}

Suppose that $x,y\in Z$ and $R\geq N_1$. Assume that the orbit segments $\alpha_{[0,R]}.x$ and $\alpha_{[0,R]}.y$ are one-dimensional $(\epsilon,\mathcal{P}_n)$-matchable. This implies there exists an increasing absolutely continuous map 
\[
h:A\to h(A)\subset[0,R]
\]
defined on a domain $A\subset[0,R]$ with 
\[
m(A)>(1-\epsilon)R\text{ and }m(h(A))\geq(1-\epsilon)R,
\]
such that for all $t\in A$
\begin{equation*}
    \begin{gathered}
        \mathcal{P}_n(\alpha_t.x)=\mathcal{P}_n(\alpha_{h(t)}.y),\qquad 
    \abs{h'(t)-1}<\epsilon.
    \end{gathered}
\end{equation*}

Choose a compact subset $\tilde{A}\subset A\cap(0,R)$ with $h(\tilde{A})\subset(0,R)$ such that
\[
m(\tilde{A})\geq (1-2\epsilon)R,\qquad m(h(\tilde{A}))>(1-2\epsilon)R.
\]
Adjoin the values $h(0)=0$, $h(R)=R$, and extend $h$ linearly across the complementary intervals of $\tilde{A}\cup\{0,R\}$. The set
\[
E=[0,R]\setminus (\tilde{A}\cup\{0,R\})
\]
is open and satisfies $m(E)\leq2\epsilon R$. We can write $E$ as a countable union of disjoint open intervals,  $E=\bigcup_{i\in I}O_i$. For each $i\in I$, define 
\[
\ell_i=\max\{m(O_i),m(h(O_i))\},
\]
where $h$ now denotes this extension. Since the complementary image intervals have total length $R-m(h(\tilde{A}))$, we have
\begin{equation}\label{eq:badImageMeasure}
\sum_{i\in I}\ell_i\leq m(E)+R-m(h(\tilde{A}))\leq5\epsilon R.
\end{equation}

We now apply Lemma~\ref{lem:Qconstruction} to the collection $\{O_i\}_{i\in I}$ and lengths $\{\ell_i\}_{i\in I}$. This yields a collection of disjoint open intervals $\{Q_j\}_{j\in J}$ covering $E$ such that 
\begin{equation}\label{eq:badSetMeasure}
m\left(\bigcup Q_j\right)=\frac{1}{\sqrt{\epsilon}}\sum_{i\in I}\ell_i\leq \frac{5\epsilon R}{\sqrt{\epsilon}}=5\sqrt{\epsilon}R.
\end{equation}

Define a new map $\bar{h}:[0,R]\to[0,R]$ by:
\begin{equation}\label{eq:h1definition}
\bar{h}(t)=\left\{
      \begin{array}{ll}
        h(t), & \hbox{if $t\notin\bigcup_jQ_j$;} \\
        \text{linear function mapping $Q_j$ onto} & \hbox{if $t\in\bigcup_jQ_j$;}\\
        \text{$[\inf_{s\in Q_j}h(s),\sup_{s\in Q_j}h(s)]$.}
      \end{array}
    \right.
\end{equation}
The map $\bar{h}$ is continuous and strictly increasing. We now estimate its derivative on the intervals $Q_j$. For a fixed $j$, let $K_j=\{i\in I:O_i\subset Q_j\}$. The slope of $\bar{h}$ on $Q_j$ is given by
\[
\bar{h}'(t)=\frac{m(h(Q_j))}{m(Q_j)}=\frac{m(h(Q_j\setminus E))+\sum_{i\in K_j}m(h(O_i))}{m(Q_j\setminus E)+\sum_{i\in K_j}m(O_i)}.
\]
Using the bounds $\abs{h'(t)-1}<\epsilon$ on the good set $Q_j\setminus E$ and $\ell_i\geq \max(m(O_i),m(h(O_i)))$, one can show that 
\begin{equation*}\label{eq:h1derivative}
\begin{gathered}
\abs{\bar{h}'(t)-1}<4\sqrt{\epsilon},\\
\end{gathered}
\end{equation*}
Consequently, we have a reasonable control for the derivative on the set $\cup Q_j$.

Notice that we also have
\[
\mathcal{P}_n(\alpha_{h(t)}.y)=\mathcal{P}_n(\alpha_{\bar{h}(t)}.y)\text{ for }t\in[0,R]\setminus(\bigcup Q_j).
\]
This together with \eqref{eq:ergodicErrorEstimates} gives that
\begin{equation*}\label{eq:badHValue}
m(A')\geq (1-8\sqrt{\epsilon})R,
\end{equation*}
where $A'$ is defined as
\[
A'=\{t\in (0,R)\setminus(\bigcup Q_j):\alpha_{\bar{h}(t)}.y\in D_k\}.
\]

Moreover, notice that $\bar{h}$ is an absolutely continuous increasing Lipschitz function defined from $[0,R]$ to itself. Modifying $\bar{h}-\id$ in arbitrarily small endpoints neighborhoods and then mollifying, we obtain $\tilde{h}\in \mathcal{D}_{[0,R]}$ such that for every $t\in[0,R]$
\begin{equation}\label{eq:denseCInfinity}
\abs{\tilde{h}(t)-\bar{h}(t)}\leq\epsilon^2_k,\quad \abs{\tilde{h}'(t)-1}<8\sqrt{\epsilon},\quad \tilde{h}(0)=0,\quad \tilde{h}(R)=R.
\end{equation}

We now verify the matching condition. For every $t\in A'$, we have $\alpha_{\bar{h}(t)}.y\in D_k$. By the definition of  $\epsilon_k$ and $D_k$ (cf.~\eqref{eq:goodSetPartition}), the symbolic name is constant for time shifts up to $\epsilon_k$. 
Hence the estimates \eqref{eq:denseCInfinity} imply that for every $t\in A'$
\[
\mathcal{P}_n(\alpha_{\bar{h}(t)}.y)=\mathcal{P}_n(\alpha_{\tilde{h}(t)}.y).
\]
Thus the symbols match on the set $A'$. Choosing the smoothing sufficiently small and using $\tilde{h}^{-1}$ in the orientation of Definition~\ref{def:epsilonPmatchable}, the derivative and mismatch estimates give total error less than $20\sqrt{\epsilon}$. We conclude that $\alpha_{[0,R]}.x$ and $\alpha_{[0,R]}.y$ are $(20\sqrt{\epsilon},\mathcal{P}_n)$-matchable. This finishes the proof of Proposition~\ref{prop:oneDimMatchImpliesMatch}.
\end{proof}

\section{Topological Kakutani balls and Horocyclic \texorpdfstring{$\R^n$}{Rn}-actions}\label{sec:horocyclicAction}

In this section, we outline an approach to show that an action is loosely Kronecker by controlling the decay rates (or non-decay rates) of certain topologically defined objects. That is, we establish a way to show the assumptions of Proposition \ref{prop:LK-Ferenczi} assuming good properties of a measure and topological properties.

\begin{definition}
    Let $\alpha : \R^n \curvearrowright (X,d)$ be a continuous action on a complete metric space $(X,d)$ by uniformly continuous maps. We say that there is an {\it $(\ve,R)$-topological matching} between $x,y \in X$ if there exists a $C^1$-function $\tau_{x,y} : [-R,R]^n \to [-(1+\ve)R,(1+\ve)R]^n$ such that
    \begin{itemize}
        \item $d(\tau_{x,y}(a)\cdot y,a\cdot x) < \ve$ and
        \item $\norm{d\tau_{x,y}(a) - \id} < \ve$
    \end{itemize}
    for all $a \in [-R,R]^n$.  
    Let $\TKak(x,\ve,R)$ denote the set of points $y \in X$ for which an $(\ve,R)$-topological matching exists.
\end{definition}

Topological matching is often easier to detect and compute than the Kakutani matching which appeared in the previous section. We establish a first lemma about rescalable matchings.

\begin{definition}
    We call an $(\ve,R)$-matching $\tau_{x,y}$ {\bf rescalable} if $\tau_{x,y}|_{[-r,r]^n}$ is an $(\ve,r)$-matching for all $r \in [R/2,R]$.
\end{definition}

\begin{lemma}
\label{lem:orbit-open}
    If $\ve$ is sufficiently small and there is a rescalable $(\ve,R)$-topological matching between $x,y \in X$, then there is a $(2\ve,(1-2\ve)R)$-topological matching between $x$ and $\alpha_a(y)$ for every $a \in [-\ve R/2,\ve R/2]^n$.
\end{lemma}

\begin{proof}
    Let $\tau_{x,y}$ denote the rescalable matching between $x$ and $y$. Then define $\tau_{x,\alpha_a(y)}(b) = \tau_{x,y}(b)-a$. Then $\tau_{x,\alpha_a(y)}$ has the same derivative as $\tau_{x,y}$, and the distance estimates follow by construction. Furthermore, if $b \in [-(1-2\ve)R,(1-2\ve)R]^n$, 
    \begin{eqnarray*}
        \tau_{x,\alpha_a(y)}(b) & \in & \tau_{x,y}([-(1-2\ve)R,(1-2\ve)R]^n) -a  \\
        & \subset & [-(1-2\ve)(1+\ve)R,(1-2\ve)(1+\ve)R]^n -a  \\
        & \subset & [-(1-\ve/2-2\ve^2)R,(1-\ve/2-2\ve^2)R]^n  \\
        & \subset & [(1-2\ve)(1+2\ve)R,(1-2\ve)(1+2\ve)R]^n.
    \end{eqnarray*}
    It follows that there is a $(2\ve,(1-2\ve)R)$-matching.
\end{proof}

For this reason, our goal is to relate the size $\TKak(x,\ve,R)$ with the matching balls $B_{R,\id}(x,\ve,\mc P)$ when the partition $\mc P$ has nice properties.

\begin{definition}
    Let $\R^n \curvearrowright X$ be a continuous action of $\R^n$ by uniformly continuous maps on a locally compact complete metric space $X$ preserving a probability measure $\mu$. An increasing sequence of finite Borel partitions $\mc P_k$ is said to be {\it neatly generating} if
    \begin{enumerate}
    \item there exists $C_k \in \mc P_k\cup\{\emptyset\}$ (which we call the {\it cusp}) such that $\mu(C_k) \to 0$ as $k \to \infty$,
    \item $\mu(\partial \mc P_k) = 0$, where $\partial \mc P_k = \bigcup_{P \in \mc P_k} \partial P$ and $\partial P$ denotes the closure of $P$ minus its interior,
    \item for every sufficiently small $\delta>0$, if $d(x,y)<\delta$ and $x,y\notin\partial_{\delta}\mathcal{P}_k\coloneq\{z\in X:d(z,\partial\mathcal{P}_k)<\delta\}$, then $x$ and $y$ belong to the same atom of $\mathcal{P}_k$, and
    \item $\mc P_k$ forms a sequence of generating partitions for $(X,\mu)$ considered as a measure space.
\end{enumerate}
\end{definition}

Any complete Riemannian manifold of finite volume admits a neatly generating sequence of partitions by adding a triangulation of $\overline{B(x_0,k)} \setminus B(x_0,k-1)$, and refining previous triangulations with mesh tending to zero as $k$ increases. 
\begin{lemma}
\label{lem:top-to-measurable}
    Fix a continuous action $\alpha : \R^n \curvearrowright (X,\mu)$ by uniformly continuous maps preserving an ergodic probability measure $\mu$ and a sequence of neatly generating partitions $\mc P_k$. Let $\delta_i,\ve_i>0$ be sequences of real numbers such that $\delta_i,\ve_i \to 0$. Assume there exist sets $A_{i,R} \subset X$ of measure at least $1-\ve_i$ such that $\mu(\TKak(x,\delta,R)) > c_i(\delta)>0$ for every $x \in A_{i,R}$, $0<\delta < \delta_i$ and sufficiently large $R > 0$, where $c_i(\delta)$ is independent of $R$.
    
    Then for sufficiently small $\ve > 0$ and sufficiently large $i,k \in \N$, there exists $c' = c'(i,k,\ve)>0$ such that for sufficiently large $R$,  $\mu(B_{R,\id}(x,\ve,\mc P_k)) > c'(i,k,\ve)$ for a set of measure at least $1- 2\ve_i$.
\end{lemma}

%By Lemma \ref{lem:BowKak-nondecay}, the assumptions of Lemma \ref{lem:top-to-measurable} are satisfied for horocylic actions, so Lemma \ref{lem:top-to-meas-horo} follows from Lemma \ref{lem:top-to-measurable}.

\begin{proof}
    For small $\ve >0$ and large $i,k \in \N$, we will show that for $x$ in a set of measure $1-\ve_i$, $\TKak(x,\delta,R) \subset B_{R,\id}(x,\ve,\mc P_k)$ for some fixed $\delta$ and sufficiently large $R$. Since we have a lower bound on $\mu(\TKak(x,\delta,R))$ for $x$ in a set of measure $1 - \ve_i$, the result follows.

    Fix a sufficiently small $\ve > 0$, an $i$ with $\epsilon_i<\epsilon/10$ and $k \in \N$. Since $\partial \mc P_k$ has measure 0, the measures of $\partial_\eta\mc P_k$ tend to 0 as $\eta$ goes to 0. Pick a sufficiently small $\eta$ such that $\mu(\partial_\eta \mc P_k) < \ve_i/2$.
    
    Choose $\delta < \min\set{\eta/4,\delta_i}$ sufficiently small. Then if $y \in \TKak(x,\delta,R)$, there exists a matching function $\tau_y : [-R,R]^n \to \R^n$ such that $d(\tau_y(a)\cdot y,a\cdot x) < \eta/4$ and $\norm{d\tau_y(a) - \id} < \delta$ for all $a \in [-R,R]^n$. Let $S_x \subset [-R,R]^n$ denote the set of $a$ for which $a\cdot x \in \partial_\eta\mc P_k$. By ergodicity, we obtain a set $A'_i \subset X$ with $\mu(A'_i) > 1-\ve_i$ such that for sufficiently large $R$ and for all $x \in A_i'$,
    \[
    m(S_x)/m([-R,R]^n) < \ve_i.
    \]
    If $a \not\in S_x$, then it follows from our construction that $\tau_y(a)\cdot y$ and $a\cdot x$ belong to the same atom of $\mc P_k$. Since 
    \[
    \norm{\tau_y-\id}_{\infty}=O_n(\delta R),
    \] 
    we can restrict to a slightly smaller cube whose image belongs to $[-R,R]^n$, and approximate the restricted map in $C^1$-topology by a smooth map by preserving the derivative bound and keeping the orbit distance bound below $\eta/2$. 
    
    It follows from Lemma \ref{lem:matchingExtension} that we have a $h\in\mathcal{D}_{[-R,R]^n}$ with $\norm{dh-\id}<\epsilon/10$, which agrees with the smoothed map outside a neighborhood of the boundary of relative volume less than $\epsilon/10$. Notice that the coding names agree there outside $S_x$. We then use $h^{-1}$ and its Jacobian bound to obtain the matching convention of Definition~\ref{def:epsilonPmatchable}, which in particular with total error less than $\epsilon$. Thus we have $y \in B_{R,\id}(x,\ve,\mc P_k)$ whenever $x \in A_{i,R} \cap A'_i$, and the conclusion follows with $c'=c_i(\delta)$.
\end{proof}

Combining Lemma \ref{lem:top-to-measurable} with Corollary \ref{cor:LK-Ferenczi}, we immediately obtain

\begin{corollary}
\label{cor:top-LK}
     If $\alpha$ satisfies the assumptions of Lemma \ref{lem:top-to-measurable}, $\alpha$ is loosely Kronecker.
\end{corollary}

\subsection{Horocyclic actions and homogeneous Kakutani balls}

We now build an even more well-behaved object, the {\it homogeneous Kakutani balls}. These will be explicitly defined through the geometry of homogeneous spaces and have precise estimates on their measure. We will then relate them to the topological Kakutani balls considered above.

Throughout this section, we will assume that Lie groups are equipped with right-invariant Riemannian metrics. Recall $U \in \Lie(G)$ is called {\it $\ad$-unipotent} if all eigenvalues of $\ad(U)$ are 0. Call a subalgebra $\mf u \subset \Lie(G)$ $\ad$-unipotent if every $U \in \mf u$ is $\ad$-unipotent.

\begin{definition}\label{def:horocyclicAction}
    Let $G$ be a Lie group and $X \in \Lie(G)$ be such that $\ad_X$ is $\R$-semisimple. The {\it horocyclic subalgebra of $X$} is the sum of the positive eigenspaces of $\ad_X$. A subalgebra $\mf u \subset \Lie(G)$ is called {\it horocyclic} if it is the horocyclic subalgebra for some $\R$-semisimple $X \in \Lie(G)$. The subgroup $\exp(\mf u) \subset G$ is called the {\it horocyclic subgroup}, and the left action of $\exp(\mf u)$ on a finite volume discrete quotient $G /\Gamma$ is called a {\it horocyclic action}.
\end{definition}

Note that all horocyclic subalgebras are unipotent, but not all unipotent subalgebras are horocyclic.
Indeed, the following lemma is a first step which was established in \cite{KanigowskiKundeVinhageWeiSlowHigher22}*{Lemma~7.2}:

\begin{lemma}
\label{lem:conformal}
    If $G$ is simple and $\mf u \subset \Lie(G)$ is a nonzero  abelian horocyclic subalgebra, then the $\R$-semisimple renormalizing element $X$ can be chosen such that:

    \begin{enumerate}
        \item $X$ belongs to some $\R$-split Cartan subalgebra $\mf a \subset \mf g$.
        \item $\mf u$ is normalized by $\mf a$.
        \item There is a set of simple roots $\Delta_s$ and some $\chi_0 \in \Delta_s$ such that if $\chi\in \Delta_s \setminus{\{\chi_0\}}$, then $\chi(X) = 0$.
        \item The stable and unstable subalgebras of $X$ have the same dimensions.
        \item ${\ad_X}|_{\mf u} = \id$.
        \item ${\ad_X}|_{\mf u^-} = -\id$, where $\mf u^-$ is the sum of the negative eigenspaces of $\ad_X$.
    \end{enumerate}
\end{lemma}

For semisimple $G$, sum these elements over the simple factors on which $\mathfrak{u}$ is nonzero, as in \cite{KanigowskiKundeVinhageWeiSlowHigher22}*{Corollary~7.4}. Given $X$ as above, we denote the positive eigenspace by $\mf u^+ = \mf u$, the negative eigenspace by $\mf u^-$, and the zero eigenspace by $\mf u^0$. Note that by Lemma~\ref{lem:conformal}, $\ad_X$ is conformal on each of its eigenspaces, and $\dim(\mf u^+) = \dim(\mf u^-)$.

\begin{introthm:horocyclic-Lk}
If $G$ is a connected semisimple Lie group and $U \subset G$ is an abelian horocyclic subgroup, and $\Gamma <G$ is a  lattice, then the homogeneous action $\alpha : U \curvearrowright G / \Gamma $ defined by $\alpha(u)\cdot g\Gamma = (ug)\Gamma$ is loosely Kronecker whenever it is ergodic for Haar probability measure.
\end{introthm:horocyclic-Lk}

We prove Theorem~\ref{introthm:horocyclic-Lk} using the criterion established in Corollary \ref{cor:top-LK}. We first establish canonical matching functions. Let $\mathcal{I}_\Gamma(\ve)$ be the set of points in $G /\Gamma$ with injectivity radius greater than $\ve$. Note that 
\[
G/\Gamma = \bigcup_{\ve >0} \mc I_\Gamma(\ve)
\]
is an increasing union of open subsets such that $\vol(\mc I_\Gamma(\ve)) \to \vol(G/\Gamma)$. In particular, by partitioning $\mc I_\Gamma(\ve)$ and adding $G/\Gamma \setminus \mathcal{I}_\Gamma(\ve)$ as the cusp, it follows that $G/\Gamma$ has a neatly generating sequence of partitions.

\medskip

We do not provide explicit matching functions, but show how to extend a  canonical matching on any compact cell to a larger cell. We first recall some features of horocyclic actions coming from Lemma~\ref{lem:conformal}. Indeed, the renormalizing element $X$ is semisimple and $\ad_X$ has three eigenspaces $\mf u^+ = \mf u$ (with eigenvalue 1), $\mf u^-$ (with eigenvalue $-1$) and $\mf u^0$ (with eigenvalue 0). It follows that
\[ [\mf u^+,\mf u^0] = \mf u^+ \qquad [\mf u^+,\mf u^-] \subset \mf u^0\]

We accordingly let $U^{0/-}$ be the connected subgroup with Lie algebra $\mathfrak{u}^-\oplus\mathfrak{u}^0$, and $U^+ = \exp(\mf u^+)$ denote the unstable subgroup, whose translations form the action we aim to study. If $g \in U^{0/-}$ is sufficiently close to $e$, there exists a {\it center-stable holonomy} $\bar{\tau}_g:U^+\to U^+$ defined near the compact cell $U_1^+$ below,  uniquely defined by
\[ 
\bar{\tau}_g(u)g = g'u \mbox{ for some }g' \in U^{0/-}.
\]

Let $U^+_\delta$ denote a box $\exp([-\delta,\delta]^n) \subset U^+$ in some basis of $\mf u^+$, using the sup norm in these coordinates. Define $\HKak^{0,-}(1,\delta)$ (denoting {\it homogeneous Kakutani}) to be the set of all $g$ in a sufficiently small neighborhood of $e$ in $U^{0/-}$ such that $\bar{\tau}_g$ is defined on $U^+_1 \subset U^+$ and both of the following hold: 
\begin{equation}
\label{eq:unipotent-katkutni-close} 
d(\bar{\tau}_g(u)g,u) <\delta \mbox{ and }\norm{d\bar{\tau}_g(u)-\id} < \delta \mbox{ for all }u \in U^+_1.
\end{equation}
Then $\HKak^{0,-}(1,\delta)$ is an open neighborhood of $e$ in $U^{0/-}$ since $\bar{\tau}_g|_{U^+_1} \to \id$ in $C^1$ as $g \to e$. Then $\bar{\tau}_g$ is a rescalable $(\delta,1)$-topological matching since $\bar{\tau}_g(e) = e$, so using the estimate on the norm implies 
\[
\bar{\tau}_g([-r,r]^n) \subset [-(1+\delta)r,(1+\delta)r]^n,\qquad r\in[1/2,1].
\]
By the translation in Lemma \ref{lem:orbit-open} and the fact that $U^+$ is transverse to $U^{0/-}$, for every $\ve >0$, choose $\delta=\delta(\ve)>0$ sufficiently small, with $2\delta<\ve$, so that
\[
\HKak(1,\ve) = U^+_\delta \HKak^{0,-}(1,\delta)\subset B_{G}(e,\epsilon/2).
\]
This is a neighborhood of $e$ with positive Haar measure $m_\ve$. For every $wg \in \HKak(1,\ve)$, the map
\[
\tau_{wg}(v) = \bar{\tau}_g(v) - w.
\]
is an $(\ve,1)$-matching between $e$ and $wg$.

For simplicity of notation, we let $\Conj_g : G \to G$ denote the conjugation map on $G$, 
\[
\Conj_g(h) = ghg^{-1}.
\]
Furthermore, denote $X^t = \exp(tX)$ and $r = \log R$ for $R\geq 1$. Then define 
\[
\HKak(R,\ve) = \Conj_{X^r}\HKak(1,\ve),
\]
and if $wg = \Conj_{X^r}(w'g')$ for $w'g' \in \HKak(1,\ve)$, let 
\[
\tau_{wg}(v) = \Conj_{X^r} \of \tau_{w'g'} \of \Conj_{X^{-r}}(v) = X^r\tau_{w'g'}(X^{-r}vX^r)X^{-r}.
\]

\begin{lemma}
\label{lem:BowKak-nondecay}
    For all $h\Gamma\in G/\Gamma$ and $R\geq 1$, $\HKak(R,\ve)h\Gamma\subset \TKak(h\Gamma,\ve,R)$.
\end{lemma}
\begin{proof}
We claim that the matching function $\tau_{wg}$ defined above is a topological matching. Indeed, $\Conj_{X^r}$ is a homothety of $\mf u$, and hence commutes with $d\tau_{w'g'}$ at every point. 

In coordinates, $\tau_{wg}(v)=R\tau_{w'g'}(v/R)$, so the derivative and range estimates are unchanged. It similarly satisfies the closeness condition. On a fixed small neighborhood of $e$ in $U^{0/-}$, conjugation by $X^r$ has a Lipschitz constant $C$ independent of $r\geq 0$, since its differential fixes $\mathfrak{u}^0$ and contracts $\mathfrak{u}^-$. Choose $\delta$ so that $C\delta<\epsilon$. Writing $g''=\bar{\tau}_{g'}(X^{-r}uX^r)g'(X^{-r}uX^r)^{-1}\in U^{0/-}$, we obtain
\begin{eqnarray*}
    d(\tau_{wg}(u)wgh,uh)  &=  &d(\tau_{wg}(u)wg,u) \\
     & = & d(\tau_g(u)g,u)  \\
    & = & d(X^r \tau_{g'}(X^{-r}uX^r)g'X^{-r},u) \\
    & = & d(X^r(g''X^{-r}uX^r)X^{-r},u) \\
    & = & d(X^rg''X^{-r},e) \leq Cd(g'',e)<C\delta<\epsilon.
\end{eqnarray*} 
\end{proof}

When $G / \Gamma$ is compact, Lemma \ref{lem:BowKak-nondecay} implies Theorem~\ref{introthm:horocyclic-Lk} by Corollary \ref{cor:top-LK} since 
\[
m_G(\HKak(R,\ve)) = m_G(\HKak(1,\epsilon))=m_\ve
\]
is independent of $R$ by unimodularity, where $m_G$ is normalized so that $G/\Gamma$ has volume one. Indeed, if $\ve$ is less than the injectivity radius of the compact quotient, then $\HKak(1,\epsilon)$ is injective at $X^{-r}h\Gamma$. Since
\[
\HKak(R,\epsilon)h\Gamma=X^{r}(\HKak(1,\epsilon)X^{-r}h\Gamma),
\]
we have $\mu(\HKak(R,\ve)h\Gamma) = m_{\epsilon}$ for all $h \in G$, and one may take $A_{i,R} = G / \Gamma$ in the assumptions of Lemma \ref{lem:top-to-measurable}.

\medskip

To finish the proof of Theorem~\ref{introthm:horocyclic-Lk} when $G/\Gamma$ is not compact, we need to build the sets $A_{i,R}$ more carefully. Let $C_i$ be the set of points with injectivity radius at most $1/i$ and let
\[
A_{i,R}=X^{\log R}(G/\Gamma\setminus C_i),\qquad \delta_i=1/i,\qquad \epsilon_i=\mu(C_i)+1/i.
\]
If $x\in A_{i,R}$ and $0<\delta<\delta_i$, then $\HKak(1,\delta)$ is injective at $X^{-\log R}x$. The same identity gives
\[
\mu(\TKak(x,\delta,R))\geq m_{\delta}>c_i(\delta)>0,\qquad c_i(\delta)=m_{\delta}/2,
\]
independently of $R$. Since $\mu(A_{i,R})>1-\epsilon_i$ and $\epsilon_i\to0$, the assumptions of Lemma~\ref{lem:top-to-measurable} are satisfied, and we conclude Theorem~\ref{introthm:horocyclic-Lk}.\qed

\appendix

\section{Nadler's Abramov formula for \texorpdfstring{$\R^n$}{Rn}-actions}\label{sec:AppendixAbramovFormula}

In this appendix, we provide a proof of the Abramov formula for $\R^n$-actions, following the approach of Nadler~\cite{NadlerAbramovNflow78}. We start with some preliminaries (cf.~\S\ref{sec:AppEntropyTilingsPartitions}) and auxiliary lemmas (cf.~\S\ref{sec:AppReductionAuxiliaryLemmas} and \S\ref{sec:AppMainProp}), then proceed to the proof of Theorem~\ref{thm:AbramovMainNadler} in \S\ref{sec:AppAbramovMainNadler}.

\subsection{Entropy for \texorpdfstring{$\R^n$}{Rn}-actions, Reduction, and Stable Partitions}\label{sec:AppEntropyTilingsPartitions}
\subsubsection{Entropy for \texorpdfstring{$\R^n$}{Rn}-actions}
In this subsection, we recall some basic definitions and properties of entropy for $\R^n$-actions. For more details and background, we refer to \cites{FeldmanREntropy80Israel,OrnsteinWeissEntropyIsomorphism}. 

\medskip

Throughout this appendix, let $(X,\mathcal{B},\mu)$ be a standard non-atomic probability space, and $\alpha$ a measurable, measure-preserving, ergodic, and free $\R^n$-action on $(X,\mathcal{B},\mu)$.

\medskip

For a finite measurable partition $\mathcal{P}$, $\v\in\R^n$, and $x\in X$, we define the coding map
\begin{equation}\label{eq:codingMap}
\alpha_x^{\mathcal{P}}(\v)=\mathcal{P}(\alpha_{\v}.x), 
\end{equation}
where $\mathcal{P}(\alpha_{\v}.x)$ denotes the atom of $\mathcal{P}$ that contains $\alpha_{\v}.x$.

Let $S\subset\R^n$ be a bounded measurable set with positive measure. We define the average Hamming metric over $S$ as
\begin{equation}\label{eq:HammingMetric}
    d_{S}^{\alpha}(\alpha_x^{\mathcal{P}},\alpha_y^{\mathcal{P}})=\frac{1}{m(S)}\int_{S}d\left(\alpha_x^{\mathcal{P}}(\v),\alpha_y^{\mathcal{P}}(\v)\right)dm(\v).
\end{equation}
Coding maps in an expression with another action in the superscript are formed using that action. Here $d$ is the discrete metric
        \[
        d(w_1,w_2)=\left\{
      \begin{array}{ll}
        1, & \hbox{if $w_1\neq w_2$;} \\
        0, & \hbox{if $w_1=w_2$.}
      \end{array}
    \right.
        \]

Let $\{F_k\}_{k\in\N}$ be a F\o lner sequence in $\R^n$. For any given $\epsilon>0$, a measurable subset $W\subset X$ is called a \textbf{$\mathcal{P}$-$\epsilon$-$F_k$-ball for $\alpha$} if for all $x,y\in W$
\[
d_{F_k}^{\alpha}(\alpha_x^{\mathcal{P}},\alpha_y^{\mathcal{P}})\leq\epsilon.
\]

For fixed $\epsilon$, we consider the set of possible growth rates:
\begin{equation} \label{eq:EntropyDef}
    \begin{aligned}
        H_{\epsilon}(\alpha,\mathcal{P}) \coloneqq \{h:&\text{ for all $\delta>0$, if $k$ is large enough one} \\
        &\text{ can cover a subset of $X$ of measure $1-\delta$ by } \\
        &\text{ fewer than $e^{hm(F_k)}$ $\mathcal{P}$-$\epsilon$-$F_k$-balls for $\alpha$}\}.
    \end{aligned}
\end{equation}
It was shown in \cite{OrnsteinWeissEntropyIsomorphism} that $H_{\epsilon}(\alpha,\mathcal{P})$ is non-empty. We define 
\[
h_{\mu}(\alpha,\mathcal{P},\epsilon)=\inf H_{\epsilon}(\alpha,\mathcal{P}).
\]
Then the entropy of the action $\alpha$ is defined as
\[
h_{\mu}(\alpha)=\sup_{\mathcal{P}}\lim_{\epsilon\to 0}h_{\mu}(\alpha,\mathcal{P},\epsilon).
\]

\begin{remark}\label{rem:remGeneralEntropy}
    \begin{enumerate}
        \item The general definition of entropy for amenable group actions is often referred to as \textbf{spatial cluster entropy} (cf.~\cite{OrnsteinWeissEntropyIsomorphism}*{II.\S.6}). However, when the acting group is abelian, the spatial cluster entropy coincides with spatial entropy; hence, we can use the simple definition given above.
        \item For $\Z^n$-actions, entropy can also be introduced using the classical language of information theory. We refer the reader to \cite{KatznelsonWeissCommuting72} for details.
    \end{enumerate}
\end{remark}

The following is a generalization of the Shannon-McMillan theorem:
\begin{theorem}[\cite{OrnsteinWeissEntropyIsomorphism}*{\S II.4~Theorem~5}]\label{thm:GeneralShannonMcMillan}
    Given $\epsilon,\delta>0$, if $F_k$ is sufficiently large, there exists a subset $E\subset X$ with $\mu(E)>1-\delta$ such that:
    \begin{enumerate}[label=(\roman*)]
        \item no $\mathcal{P}$-$\epsilon$-$F_k$-ball for $\alpha$ in $E$ has measure greater than $e^{-(h-\delta)m(F_k)}$;
        \item $E$ can be covered by $\leq e^{(h+\delta)m(F_k)}$ $\mathcal{P}$-$\epsilon$-$F_k$-balls for $\alpha$; 
    \end{enumerate}
    where $h=h(\alpha,\mathcal{P},\epsilon)$.
\end{theorem}

Theorem~\ref{thm:GeneralShannonMcMillan} implies the following useful observation regarding the measure threshold in the definition of entropy.
\begin{remark}\label{rem:modifiedEntropyDef}
In the definition of $H_{\epsilon}(\alpha,\mathcal{P})$, the requirement to cover a $1-\delta$ portion of the space $X$ for every $\delta>0$ can be relaxed. Given any fixed constant $c\in(0,1)$, if we define
\begin{equation*}
    \begin{aligned}
        \tilde{H}_{\epsilon}(\alpha,\mathcal{P})=\{h:&\text{ if $k$ is sufficiently large, one can cover}\\ 
        &\text{ a subset of $X$ of measure $c$ by fewer} \\
        &\text{ than $e^{hm(F_k)}$ many $\mathcal{P}$-$\epsilon$-$F_k$-balls for $\alpha$}\},
    \end{aligned}
\end{equation*}
then 
\[
\inf \tilde{H}_{\epsilon}(\alpha,\mathcal{P})=\inf H_{\epsilon}(\alpha,\mathcal{P}).
\]
\end{remark}

The following definitions will be central to our proof:
\begin{definition}\label{def:HammingBall}
    For $x\in X$, $r>0$, a subset $S\subset \R^n$, and a finite measurable partition $\mathcal{P}$, we define the \textbf{Hamming ball} as
    \[
    B^{H,\alpha}_r(x,S,\mathcal{P})=\{y\in X:d_S^{\alpha}(\alpha^{\mathcal{P}}_x,\alpha^{\mathcal{P}}_y)<r\}.
    \]
\end{definition}
\begin{definition}
    Let $\mathcal{P}$ be a finite measurable partition of $X$, $S\subset\R^n$ and $r,c>0$. A family of subsets of $X$ is said to be a \textbf{$(\mathcal{P},S,r)$-$c$ family for $\alpha$} if each set in the family is contained in $B^{H,\alpha}_r(x,S,\mathcal{P})$ for some $x\in X$ and the total measure of the union of these sets is greater than $c$.
\end{definition}

\subsubsection{Reduction to even time changes}
Let $\rho$ be a reparametrization cocycle for $\alpha$, and $\tau$ its inverse cocycle. Since Proposition~\ref{prop:jRhoProp} implies that $J(\rho):\R^n\to\R^n$ is a linear map, we may restate Theorem~\ref{thm:AbramovMainNadler} in the following equivalent form:
\begin{theorem}
    Let $\alpha$ be a measurable, measure-preserving, ergodic, and free $\R^n$-action on a standard non-atomic probability space $(X,\mathcal{B},\mu)$, and $\alpha^{\tau}$ a tempered time change of $\alpha$. Then 
    \[
    h_{\mu}(\alpha)= h_{\mu^{\tau}}(\alpha^{\tau})\cdot\abs{\det J(\rho)},
    \]
    where $\mu^{\tau}$ is the invariant measure of $\alpha^{\tau}$ induced by the time change defined in equation~\eqref{eq:NewMeasureDensityEquation}.
\end{theorem}

Our goal is to reduce our proof to the case of \textbf{even time changes} (cf.~Definition~\ref{def:kakutani}), i.e.,  $J(\rho)=\id$. This reduction relies on the following proposition describing the entropy scaling under linear coordinate changes. Note that in this case, $\mu^{\tau}=\mu$.

\begin{proposition}\label{prop:linearTimeChange}
    Let $\alpha$ be a measurable, measure-preserving, ergodic, and free $\R^n$-action on a standard non-atomic probability space $(X,\mathcal{B},\mu)$. Suppose $\mathcal{L}:\R^n\to\R^n$ is a nondegenerate linear transformation and $\tilde{\alpha}$ satisfies $\tilde{\alpha}_{\v}.x=\alpha_{\mathcal{L}(\v)}.x \text{ for $\mu$-a.e. $x\in X$ and all $\v\in\R^n$.}$
    Then
    \[
    h_{\mu}(\alpha)= h_{\mu}(\tilde{\alpha})\cdot\abs{\det \mathcal{L}^{-1}}.
    \]
\end{proposition}
\begin{proof}[Proof of Proposition~\ref{prop:linearTimeChange}]
    Let $\{\mathcal{P}_k\}_{k\in\N}$ be a generating increasing sequence of finite measurable partitions. Recall that $C_N=[-N,N]^n$. 
    
    Fix a $k\in\N$ and $\epsilon>0$. By the change of variables formula for integration, for any $x,y\in X$, the Hamming distances satisfy:
    \[
    d_{C_N}^{\alpha}(\alpha_x^{\mathcal{P}_k},\alpha_y^{\mathcal{P}_k})=d_{\mathcal{L}^{-1}(C_N)}^{\tilde{\alpha}}(\alpha_x^{\mathcal{P}_k},\alpha_y^{\mathcal{P}_k}),
    \]
    where the distance on the left is computed with respect to the action $\alpha$ over the domain $C_N$, and the right with respect to $\tilde{\alpha}$ over $\mathcal{L}^{-1}(C_N)$.
    
    Consequently, a collection $\mathcal{U}$ is a $(\mathcal{P}_k,C_N,\epsilon)$-$0.9$ family for $\alpha$ if and only if $\mathcal{U}$ is also a $(\mathcal{P}_k,\mathcal{L}^{-1}(C_N),\epsilon)$-$0.9$ family for $\tilde{\alpha}$. This observation, combined with  Remark~\ref{rem:modifiedEntropyDef}, implies that 
    \[
    h_{\mu}(\alpha,\mathcal{P}_k,\epsilon)=h_{\mu}(\tilde{\alpha},\mathcal{P}_k,\epsilon)\cdot\abs{\det\mathcal{L}^{-1}}.
    \]
    Letting $\epsilon\to0$ and taking the supremum over partitions yields
    \[
    h_{\mu}(\alpha)=h_{\mu}(\tilde{\alpha})\cdot\abs{\det\mathcal{L}^{-1}},
    \]
    which concludes the proof of Proposition~\ref{prop:linearTimeChange}.
\end{proof}

\begin{remark}\label{rem:identityJrho}
    By applying Proposition~\ref{prop:linearTimeChange} with $\mathcal{L}=J(\rho)^{-1}$, we can reduce the general case to a system with even time changes. Thus we will assume from now on that our reparametrization map $\rho$ satisfies $J(\rho)=\id$.
\end{remark}

\subsubsection{Stable Partitions}
In this subsection, we recall results from \cite{RudolphRectangularTilings86} regarding rectangular tilings of $\R^n$, which are essential for constructing the partitions used in our later proofs.

\begin{definition}[\cite{RudolphRectangularTilings86}] \label{def:Tiles}
\begin{enumerate}
    \item Fix irrational numbers $\alpha_1,\ldots,\alpha_n$. Let $\mathscr{B}=\{B_1,\ldots,B_{2^n}\}$ denote the collection of $n$-dimensional rectangles of the form
    \[
    \prod_{i=1}^{n}[0,d_i),
    \]
    where each side length $d_i$ is either $\alpha_i$ or $1$. The elements of $\mathscr{B}$ are referred to as \textbf{basic tiles}.
    \item A \textbf{tiling} of $\R^n$ is a partition of $\R^n$ into sets $\{T_i\}_{i\in\N}$, where each tile satisfies
    \[
    T_i=\v_i+B_{j(i)},
    \]
    for some vector $\v_i\in\R^n$ and some basic tile $B_{j(i)}\in\mathscr{B}$. The vector $\v_i$ is called the \textbf{base point} of the tile $T_i$, and we denote the tiling by $\mathscr{T}$.
\end{enumerate}
\end{definition}

\begin{definition}[\cite{RudolphRectangularTilings86}]
A \textbf{tiling} for the action $\alpha$ is a representation of $X$ as
    \[
    X_0=\bigcup_{i=1}^{2^n} (Y_i\times B_i),
    \]
    where $X_0\subset X$ is a $\alpha$-invariant subset of full measure, $\{Y_i\}_{i=1}^{2^n}$ are disjoint Lebesgue spaces, and each $B_i$ is a basic tile equipped with the standard Lebesgue structure of $\R^n$.
\end{definition}

\begin{remark}\label{rem:tilingMove}
    The action on the tiling is defined naturally: if $(y,\v_0)\in Y_i\times B_i$ and $\v\in\R^n$ are such that $\v_0+\v\in B_i$, then
    \[
    \alpha_{\v}.(y,\v_0)=(y,\v_0+\v).
    \]
\end{remark}

By \cite{RudolphRectangularTilings86}*{Lemma~1.8, Corollary~1.11, and Theorem~3.9}, we obtain that any such action admits a tiling. 
\begin{theorem}[\cite{RudolphRectangularTilings86}]\label{thm:RudolphTiling2n}
Any measurable, measure-preserving and free $\R^n$-action $\alpha$ can be measurably tiled by rectangles of the $2^n$ shapes described in Definition~\ref{def:Tiles}. 
\end{theorem}

\medskip

We now have the necessary ingredients to construct a sequence of generating partitions $\{\mathcal{P}_k\}_{k\in\N}$ with specific stability properties. To formalize this, we introduce the concept of stable partition.
\begin{definition}\label{def:stablePartition}
    Let $\eta>0$ and $\mathcal{P}$ a finite measurable partition of $(X,\mathcal{B},\mu)$. We say that $\mathcal{P}$ is an \textbf{$\eta$-stable partition} for the $\R^n$-action $\alpha$ with parameters $(A,\delta,N)$ if there exist a set $A\subset X$ with $\mu(A)>1-\eta$, and constants $\delta,N>0$ satisfying the following conditions:
\begin{enumerate}[label=(\roman*)]
\item\label{item:stableProp1} For almost every $x\in X$, 
\[
\frac{m(\{\v\in C_N:\alpha_{\v}.x\in A\})}{m(C_N)}\geq 1-5\eta,
\]
    \item\label{item:stableProp2} For every $x\in A$, the $\delta$-neighborhood of the orbit of $x$ remains within the same partition atom:
    \begin{equation*}
    \{\alpha_{\v}.x:\v\in C_{\delta}\}\subset\mathcal{P}(x), 
    \end{equation*}
where $\mathcal{P}(x)$ denotes the atom of $\mathcal{P}$ containing $x$.
\end{enumerate}
\end{definition}

The next proposition guarantees the existence of a generating sequence of $\eta$-stable partitions for $\R^n$-actions.

\begin{proposition}\label{prop:goodPartition}
    Let $\alpha$ be a free and ergodic $\R^n$-action on $(X,\mathcal{B},\mu)$. There exists a generating sequence of finite measurable partitions $\{\mathcal{P}_k\}_{k\in\N}$ such that for every fixed $k\in\N$ and $\eta\in(0,1)$, there exist constants $\delta\in(0,1)$ and $T\in\N$, and a measurable set $A\subset X$ with $\mu(A)\geq 1-\eta$ such that for every $N\geq T$, $\mathcal{P}_k$ is an $\eta$-stable partition for $\alpha$ with parameters $(A,\delta,N)$. Here $T$ is independent of $k$.
\end{proposition}
\begin{proof}
    Recall that Theorem~\ref{thm:RudolphTiling2n} provides a tiling of a full measure invariant set $X_0\subset X$ by basic tiles $\{B_i\}_{i=1}^{2^n}$ and disjoint Lebesgue spaces $\{Y_i\}_{i=1}^{2^n}$. Denote the base space $Y=\bigcup Y_i$ and let $\mathcal{Q}$ be a measurable partition of $Y$ defined by
    \[
    \mathcal{Q}=\{Y_1,\ldots,Y_{2^n}\}.
    \]

    Let $\mathcal{P}'_k$ be an increasing generating sequence of finite measurable partitions of $Y$ refining $\mathcal{Q}$. We subdivide each basic tile $B_i$ into $(100k)^n$ congruent sub-rectangles $\{B_{i,j}\}_j$ by dividing each side into $100k$ equal segments. We define the partition $\mathcal{P}_k$ on $X$ by:
    \[
    \mathcal{P}_k=\{P'\times B_{i,j}:i=1,\ldots,2^n,~j=1,\ldots,(100k)^n,~P'\in\mathcal{P}'_k,~P'\subset Y_i\}.
    \]
    The sequence $\{\mathcal{P}_k\}_{k\in\N}$ is clearly generating.

    \medskip

    Fix $\eta\in (0,1)$ and $k\in\N$. We construct $A$ by shrinking the tiles. Let $\v_{i,j}$ be the center of $B_{i,j}$ and define the contraction $C_{i,j}^{\lambda}(\v)$ by
    \[
    C_{i,j}^{\lambda}(\v)=\lambda(\v-\v_{i,j})+\v_{i,j},\quad \forall \lambda\in\R^+,\v\in\R^n.
    \]
    Choose $\eta'>0$ small enough such that 
    \[
    (1-\eta')^n>1-\eta\text{ and }(1+\eta'/100)^n<2.
    \]    
    Define the ``good core'' of each tile as 
    \[
    B'_{i,j}=C_{i,j}^{1-\eta'}(B_{i,j}),
    \]
    and let
    \begin{equation}\label{eq:AketaDef}
    A=\bigcup_{i=1}^{2^n}\bigcup_{j=1}^{(100k)^n}\left(Y_i\times B'_{i,j}\right).
    \end{equation}
    Set 
    \[
    \delta=\min_{i,j}\left(\eta'l_{\min}(B_{i,j})/2\right),
    \]
    where $l_{\min}(B_{i,j})$ denotes the minimal side length of the rectangle $B_{i,j}$. By construction, for any $x\in A$, the $\delta$-ball of the orbit remains inside a single rectangle $B_{i,j}$ and fixed base $P'$, and hence satisfies
    \[
    \mathcal{P}_k(x)=\mathcal{P}_k(\alpha_{\v}.x),\qquad \forall\v\in C_{\delta}.
    \]
    Moreover, $\mu(A)=(1-\eta')^n\geq1-\eta$. Consequently, we verify the Condition~\ref{item:stableProp2} of the stable partition (cf.~Definition~\ref{def:stablePartition}).
    
    \medskip

    To verify Condition~\ref{item:stableProp1}, let $L_{\max}$ be the maximal side length of the basic tiles. Set
    \[
    T=\frac{100L_{\max}}{\eta'}.
    \]
    
    Let $X_0$ be the full measure invariant set from Theorem~\ref{thm:RudolphTiling2n}. Fix $x\in X_0$ and $N\geq T$. It follows from Remark~\ref{rem:tilingMove} that the orbit set 
    \[
    O_N(x)=\{\alpha_{\v}.x:\v\in C_N\}
    \]
    is a union of full or partial tiles. Define
    \[
    D_N(x)=\{\v\in C_N:\alpha_{\v}.x\in A\}.
    \]
    Notice that the set $C_N\setminus D_N(x)$ consists of two parts: the regions $B_{i,j}\setminus B'_{i,j}$, and the boundary effects where $C_N$ cuts a tile.

    Notice that the density of the good set $B'_{i,j}$ belong to every sub-rectangle $B_{i,j}$ is exactly $(1-\eta')^n$, the density of $D_N(x)$ inside any \textbf{fully covered} tile satisfies 
    \[
    (1-\eta')^n>1-\eta. 
    \]
    As for the boundary effect, its contribution is controlled by the ratio $L_{\max}/N$, which is in particular less than $\eta$ by our choice of $N$. Altogether, we have:
    \[
    \frac{m(D_N(x))}{m(C_N)}\geq (1-\eta')^n\cdot\frac{(N-L_{\max})^n}{N^n}\geq 1-2\eta.
    \]
    Hence, we finish the estimates and they satisfy the required bound in Condition~\ref{item:stableProp1} of Definition~\ref{def:stablePartition} for all $x\in X_0$, which finish the proof.
\end{proof}

\medskip

The following lemma addresses the stability of the stable partitions under small perturbations of the time change.

\begin{lemma}\label{lem:stablePertQkphi}
Let $r,\delta>0$. Suppose $\alpha^{\tau}$ is a $10^{-2}$-tempered time change of $\alpha$, and $\mathcal{P}$ is an $r$-stable partition \textbf{for $\alpha^{\tau}$} with parameters $(A,\delta,N)$. Choose $\eta\in(0,n^{-1}\min(r,\delta)/200)$ sufficiently small such that
\[
1-(1-\eta)^n<r/200.
\]
If $\alpha^{\tilde{\tau}}$ is an $\eta/N$-tempered time change cohomologous to $\alpha^{\tau}$ through $\psi$, that is  $\alpha^{\tilde{\tau}}_{\v}\circ\psi=\psi\circ\alpha^{\tau}_{\v}$, then the partition $\psi(\mathcal{P})$ is a $2r$-stable partition \textbf{for $\alpha$} with parameters $(\psi(A),\delta/2,N)$.
\end{lemma}
\begin{proof}[Proof of Lemma~\ref{lem:stablePertQkphi}]
It follows from the cohomology relation that the system $(\alpha^{\tilde{\tau}},\psi(\mathcal{P}))$ has the stability properties of $(\alpha^{\tau},\mathcal{P})$. More precisely, since $\psi_*\mu^{\tau}=\mu^{\tilde{\tau}}$, the partition $\psi(\mathcal{P})$ is an $r$-stable partition for $\alpha^{\tilde{\tau}}$ with parameters $(\psi(A),\delta,N)$.

We transfer these properties to $\alpha$ using the ``closeness'' of the time changes. Since $\alpha^{\tilde{\tau}}$ is $\eta/N$-tempered, for $\mu$-a.e. $x\in X$ and all $\v\in C_N$:
    \begin{equation*}
    \norm{\tilde{\tau}(x,\v)-\v}\leq n\eta.
    \end{equation*}

\textbf{1. Density.} Since $d\mu^{\bar{\tau}}/d\mu=1+O_n(\eta/N)$, we have for sufficiently small $\eta$,
\[
\mu(\psi(A))>1-2r.
\]
Notice that the sets $C_N$ and the reparametrization image 
\[
\{\tilde{\tau}(x,\v):\v\in C_N\}
\]
differ in symmetric difference by a boundary term with size $O(\eta)$. Thus, the large density of visiting $\psi(A)$ under $\alpha^{\tilde{\tau}}$ implies large density under $\alpha$,  with an additional error bounded by $r/200$. Here the first action is parametrized by time $\tilde{\tau}(x,\v)$ and the second action is parametrized by $\v$. This satisfies Condition~\ref{item:stableProp1} of Definition~\ref{def:stablePartition} for $\alpha$ with parameter  $2r$, $\psi(A)$ and $N$.

\medskip

\textbf{2. Stability. } Let $x\in \psi(A)$. Stability for $\alpha^{\tilde{\tau}}$ implies 
\begin{equation*}
    \{\alpha^{\tilde{\tau}}_{\v}.x:\v\in C_{\delta}\}\subset\psi(\mathcal{P})(x).
    \end{equation*}
    For any vector $\w\in C_{\delta/2}$, let $\s$ be the vector such that $\tilde{\tau}(x,\s)=\w$. The derivative bound gives
    \[
    \norm{\s}_{\infty}\leq(1+\sqrt{n}\eta/N)\delta/2<\delta.
    \]
    Therefore,
    \[
    \alpha_{\w}.x=\alpha^{\tilde{\tau}}_{\s}.x\in \{\alpha^{\tilde{\tau}}_{\v}.x:\v\in C_{\delta}\}\subset\psi(\mathcal{P})(x).
    \]
    This verifies Condition~\ref{item:stableProp2} of Definition~\ref{def:stablePartition} with parameters $\delta/2$, $\psi(A)$ and $N$.
\end{proof}

\subsection{Auxiliary Lemmas}\label{sec:AppReductionAuxiliaryLemmas}

In this section, we establish the necessary preparations for the proof of Theorem~\ref{thm:AbramovMainNadler}. We work on a common invariant conull set $X_0$ on which the cocycle identities and estimates below hold.

\subsubsection{A local name comparison lemma}
Recall that Lemma~\ref{lemma:goodBasis} guarantees the existence of a basis $\mathcal{E}=\{\e_1,\ldots,\e_n\}$ of $\R^n$ such that the transformation $\alpha_{\e_i}$ is totally ergodic for every $i=1,\ldots,n$. Combining this with Proposition~\ref{prop:linearTimeChange}, we may assume without loss of generality that $\mathcal{E}$ is the standard basis of $\R^n$. 

Let $\mathcal{P}$ be a finite measurable partition of $X$. Recall the notation $\alpha_x^{\mathcal{P}}$ from \eqref{eq:codingMap}. Fix an arbitrary atom $P_0\in\mathcal{P}$ to serve as a default value. We define the \textbf{reparametrized name} $\alpha^{\mathcal{P}}_{x,\rho,N}:\R^n\to\mathcal{P}$ by
\begin{equation}\label{eq:alphaRhoName}
    \alpha^{\mathcal{P}}_{x,\rho,N}(\v)=\left\{
      \begin{array}{ll}
        \alpha^{\mathcal{P}}_x(\tau(x,\v)), & \hbox{if $\v\in \rho(x,C_N)\cap C_N$;} \\
        \text{$P_0$}, & \hbox{otherwise.}
      \end{array}
    \right.
\end{equation}
For $x,y\in X$, the Hamming distance between their reparametrized names over $C_N$ is defined by
\begin{equation*}
    \begin{gathered}
        d_{C_N}^{\alpha}(\alpha^{\mathcal{P}}_{x,\rho,N},\alpha^{\mathcal{P}}_{y,\rho,N})=\frac{1}{m(C_N)}\int_{C_N}d\left(\alpha^{\mathcal{P}}_{x,\rho,N}(\v),\alpha^{\mathcal{P}}_{y,\rho,N}(\v)\right)dm(\v),
    \end{gathered}
\end{equation*}
where $d(\cdot,\cdot)$ is the discrete metric on $\mathcal{P}$.

\medskip

The following lemma is an analogue of \cite{NadlerAbramovNflow78}*{Lemma~$3$}. Recall that an $R$-tempered time change $\alpha^{\tau}$ satisfies the bound (cf.~Definition~\ref{def:temperedTimeChange}):
\[
    \norm{\partial_{\v} \tau(x,\0)-\int_X\tau(y,\cdot)d\mu^{\tau}(y)}<R,\quad\norm{\partial_{\v} \rho(x,\0)-\int_X\rho(y,\cdot)d\mu(y)}<R.
    \]
\begin{lemma}\label{lem:nameComparison}
    Let $\alpha$ be an ergodic and free $\R^n$-action on $(X,\mathcal{B},\mu)$, $\epsilon\in(0,1)$ and $N\geq 1$. Assume that $\mathcal{P}$ is an $\epsilon/200$-stable partition for $\alpha$ with parameters $(A,\delta,N)$. Set $\eta=10^{-100}n^{-1}\min\left(\delta^{10},\epsilon^{10}\right)$. If $\alpha^{\tau}$ is an $\eta/N$-tempered time change, then for all $x,y\in X_0$ satisfying
    \[
    d_{C_N}^{\alpha}(\alpha_x^{\mathcal{P}},\alpha^{\mathcal{P}}_y)<\eta,
    \]
    we have
    \[
    d_{C_N}^{\alpha}(\alpha^{\mathcal{P}}_{x,\rho,N},\alpha^{\mathcal{P}}_{y,\rho,N})<\epsilon.
    \]
\end{lemma}
\begin{proof}
    Fix $\epsilon$, $\mathcal{P}$, $\eta$, $N$ and $\alpha^{\tau}$ as in the assumption. 
    
    \medskip

Let $X_0$ be the common invariant conull set specified above. Let $x,y\in X_0$ satisfy $d_{C_N}^{\alpha}(\alpha^{\mathcal{P}}_x,\alpha^{\mathcal{P}}_y)<\eta$. We define the following sets:
\begin{enumerate}[label=(\alph*)]
    \item Let $W_1=\{\v\in C_N:\alpha^{\mathcal{P}}_x(\v)=\alpha^{\mathcal{P}}_y(\v)\}$. The condition on the Hamming distance implies
    \[
    m(W_1)\geq(1-\eta)m(C_N).
    \]

    \item Let $W_2=\{\v\in C_N:\alpha_{\v}.x,\alpha_{\v}.y\in A\}$. By Definition~\ref{def:stablePartition} (specifically Condition~\ref{item:stableProp1} for the parameter $N$), we have
    \[
    m(W_2)\geq (1-\epsilon/20)m(C_N).
    \]

    \item \label{item:c} Let 
    \[
    W_3=\rho(x,W_1\cap W_2)\cap\rho(y,W_1\cap W_2)\cap C_N.
    \]
    From the definition of $W_1$ and $W_2$, we have 
    \[
    m(W_1\cap W_2)\geq (1-\epsilon/10)m(C_N).
    \]
    Since $\alpha^{\tau}$ is $\eta/N$-tempered, we have $\norm{\rho(z,\v)-\v}<n\eta$ for all $\v\in C_N$ and $z\in\{x,y\}$. This, together with the assumption that $J(\rho)=\id$ (cf.~Proposition~\ref{prop:jRhoProp}~\ref{item:JrhoErgodicSum}), implies that the mass lost at the boundary of $C_N$ is negligible. Specifically, our choice of $\eta$ guarantees
    \begin{equation*}
        m(W_3)\geq(1-\epsilon)m(C_N).
    \end{equation*} 
\end{enumerate}

    \medskip
    
    We now proceed to the proof. Suppose $\v\in W_3$. By definition, 
    \[
    \tau(x,\v),\tau(y,\v)\in W_1\cap W_2.
    \]
    Since these points are in $W_1$, they have the same name: 
    \begin{equation}\label{eq:xyNHBDname1}
         \mathcal{P}(\alpha_{\tau(x,\v)}.x) =\mathcal{P}(\alpha_{\tau(x,\v)}.y).
    \end{equation}
    Furthermore, since they are in $W_2$, their orbits are stable. Specifically, for any $\w\in C_{\delta}$
\begin{equation}\label{eq:xyNHBDname2}
        \begin{aligned}
            \mathcal{P}(\alpha_{\tau(x,\v)}.x)& =\mathcal{P}(\alpha_{\w+\tau(x,\v)}.x),\\
            \mathcal{P}(\alpha_{\tau(y,\v)}.y)& =\mathcal{P}(\alpha_{\w+\tau(y,\v)}.y).
        \end{aligned}
    \end{equation}
    Notice that for $\v\in C_N$, the closeness of the time change implies
    \[
    \norm{\tau(x,\v)-\tau(y,\v)}\leq2\eta\leq \delta. 
    \]
    Combining this proximity with the stability condition \eqref{eq:xyNHBDname2}, we deduce that the point $\alpha_{\tau(y,\v)}.y$ lies in the same partition atom as $\alpha_{\tau(x,\v)}.y$, which matches $\alpha_{\tau(x,\v)}.x$ (cf.~\eqref{eq:xyNHBDname1}). Thus:
    \[
    \alpha^{\mathcal{P}}_{x,\rho,N}(\v)=\alpha^{\mathcal{P}}_{y,\rho,N}(\v) \text{ for all } \v \in W_3.
    \]
    Consequently, 
    \begin{equation*}
        \begin{aligned}
            d_{C_N}^{\alpha}(\alpha^{\mathcal{P}}_{x,\rho,N},\alpha^{\mathcal{P}}_{y,\rho,N})\leq\frac{m(C_N\setminus W_3)}{m(C_N)}<\epsilon.
        \end{aligned}
    \end{equation*}   
    This concludes the proof of Lemma~\ref{lem:nameComparison}.
\end{proof}

\subsubsection{A global name comparison lemma}
Lemma~\ref{lem:nameComparison} provides name comparison at a fixed local scale $N$. To obtain a global comparison on a much larger scale $M$, we must ``paste'' these local estimates together. The following lemma achieves this by assuming additional closeness of the reparameterization maps on a discrete lattice of $\R^n$. 

\medskip

We first fix some notations.
\begin{enumerate}[label=(\roman*)]
    \item\label{item:GammaN} Given $N\in\N$, let $\Gamma_N$ be the lattice subgroup consisting of the centers of a tiling of $\R^n$ by copies of the cube $C_N$. 
    \item\label{item:WSetDef} Let $\mathcal{P}$ be a finite measurable partition of $X$. For $y_1,y_2\in X$, $s\in(0,1)$ and $M,N\in\N$, we define the set of ``good local matching times'' as:
\begin{equation*}%\label{eq:W1Def}
    W_{N,M}^{\mathcal{P}}(y_1,y_2,s)=\left\{\v\in C_{M}\cap\Gamma_N:d_{C_{N}}^{\alpha}(\alpha_{\alpha_{\v}.y_1}^{\mathcal{P}},\alpha_{\alpha_{\v}.y_2}^{\mathcal{P}})<\sqrt{s}\right\}.
\end{equation*}
\end{enumerate}

\begin{lemma}\label{lemma:discreteCounting}
    Fix $r\in(0,1)$. Let $\mathcal{P}$ be an $r/2000$-stable partition for $\alpha$ with parameters $(A,\delta,N)$. Assume that $\eta\in(0,10^{-1000}n^{-1}\min(\delta^{10},r^{100}))$ is sufficiently small such that
    \[
    (1-2\eta)^n\geq 1-r/200.
    \]
    Then there exists $M'\in\N$ such that for every $M\geq M'$ the following holds: if $\alpha^{\tau}$ is an $\eta/N$-tempered time change, then for all $y_1,y_2\in X_0$ satisfying 
    \[
    d_{C_M}^{\alpha}(\alpha_{y_1}^{\mathcal{P}},\alpha_{y_2}^{\mathcal{P}})<\eta^2/100
    \]
    and
    \[
    \norm{\rho(y_1,\v)-\rho(y_2,\v)}< \delta/8\quad\text{for all } \v\in W, 
    \]
    where $W\subset W_{N,M}^{\mathcal{P}}(y_1,y_2,\eta^2/100)$ and 
    \[
    \card{C_M\cap\Gamma_N}-\card{W}\leq(\eta/5)\card{C_M\cap\Gamma_N},
    \]
    we have
    \[
    d_{C_M}^{\alpha}(\alpha_{y_1,\rho,M}^{\mathcal{P}},\alpha_{y_2,\rho,M}^{\mathcal{P}})<r/2.
    \]
\end{lemma}
\begin{proof}[Proof of Lemma~\ref{lemma:discreteCounting}]
We first specify the constants. Choose $M'$ sufficiently large such that for all $M\geq M'$
    \[
    m(C_{M-N})/m(C_M)\geq 1-\eta^2/100.
    \]
    Fix an $M\geq M'$, and let $W$ be as in the hypothesis.
    
\medskip

    By the definition of $W$, for every $\v\in W$, we have
    \begin{equation*}
        d_{C_{N}}^{\alpha}(\alpha_{\alpha_{\v}.y_1}^{\mathcal{P}},\alpha_{\alpha_{\v}.y_2}^{\mathcal{P}})<\eta/10.
    \end{equation*}
    Applying Lemma~\ref{lem:nameComparison} with $\epsilon=r/10$ to the points $\alpha_{\v}.y_1$ and $\alpha_{\v}.y_2$, and using the choice of $\eta$, we obtain the following local reparametrization comparison for $\v\in W$
    \begin{equation}\label{eq:y1y2CloseDifferentTime}
    d_{C_{N}}^{\alpha}(\alpha^{\mathcal{P}}_{\alpha_{\v}.y_1,\rho,N},\alpha^{\mathcal{P}}_{\alpha_{\v}.y_2,\rho,N})<r/10.
\end{equation}

Now we adjust the time parameter to align the orbits. For a fixed $\v\in W$, define 
\begin{equation}\label{eq:DefV'}
    \v'=\tau(y_2,\rho(y_1,\v)).
\end{equation}
Since $\alpha^{\tau}$ is an $\eta/N$-tempered time change for $\alpha$, the derivative bound makes $\t\mapsto\tau(y,\t)$ globally $2$-Lipschitz. Thus:
\begin{equation}\label{eq:vv'Difference}
    \begin{aligned}
        \norm{\v'-\v}=&\norm{\tau(y_2,\rho(y_1,\v))-\tau(y_2,\rho(y_2,\v))}\\
        \leq &2\norm{\rho(y_1,\v)-\rho(y_2,\v)}< \delta/4.
    \end{aligned}
\end{equation}

\medskip

We now compare the names at the shifted time moments. Let 
\[
\w\in \rho(\alpha_{\v}.y_2,C_{N})\cap\rho(\alpha_{\v'}.y_2,C_{N})\cap C_{N}.
\]
Then the values of the reparametrized names satisfy:
\begin{equation}\label{eq:v'NameCompare}
    \alpha_{\alpha_{\v'}.y_2}^{\mathcal{P}}(\tau(\alpha_{\v'}.y_2,\w))=\alpha_{\alpha_{\v}.y_2}^{\mathcal{P}}(\tau(\alpha_{\v'}.y_2,\w)+\v'-\v).
\end{equation}
By using the cocycle identity and \eqref{eq:DefV'}, we now can estimate the time drift inside the argument:
\begin{equation}\label{eq:taukDiffPt}
    \begin{aligned}
        &\tau(\alpha_{\v'}.y_2,\w)-\tau(\alpha_{\v}.y_2,\w)\\
        =&\tau(\alpha^{\tau}_{\rho(y_2,\v')}.y_2,\w)-\tau(\alpha^{\tau}_{\rho(y_2,\v)}.y_2,\w)\\
        =&\tau(\alpha^{\tau}_{\rho(y_1,\v)}.y_2,\w)-\tau(\alpha^{\tau}_{\rho(y_2,\v)}.y_2,\w)\\
        =&\tau(y_2,\w+\rho(y_1,\v))-\tau(y_2,\rho(y_1,\v))-\tau(y_2,\w+\rho(y_2,\v))\\&+\tau(y_2,\rho(y_2,\v)).
    \end{aligned}
\end{equation}
It follows from the $2$-Lipschitz property of $\tau(y_2,\cdot)$ and the assumption on $\rho$-closeness that
\begin{equation*}
\begin{aligned}
     \norm{\tau(y_2,\w+\rho(y_1,\v))-\tau(y_2,\w+\rho(y_2,\v))}\leq &2\norm{\rho(y_1,\v)-\rho(y_2,\v)}\\
     < &\delta/4.
\end{aligned}
\end{equation*}
Similarly, 
\begin{equation*}
    \begin{aligned}
        &\norm{\tau(y_2,\rho(y_1,\v))-\tau(y_2,\rho(y_2,\v))}\leq 2\norm{\rho(y_1,\v)-\rho(y_2,\v)}< \delta/4.
    \end{aligned}
\end{equation*}
These inequalities together with \eqref{eq:taukDiffPt} imply that
\begin{equation*}
    \begin{aligned}
     \norm{\tau(\alpha_{\v'}.y_2,\w)-\tau(\alpha_{\v}.y_2,\w)}< \delta/2.
    \end{aligned}
\end{equation*}
Combining this with \eqref{eq:vv'Difference} and \eqref{eq:v'NameCompare}, the total time shift is small:
\begin{equation*}
    \norm{\tau(\alpha_{\v'}.y_2,\w)+\v'-\v-\tau(\alpha_{\v}.y_2,\w)}<\delta,
\end{equation*}
Since $\mathcal{P}$ is an $r/2000$-stable partition with parameters $(A,\delta,N)$, this small shift implies that whenever $\alpha_{\tau(\alpha_{\v}.y_2,\w)}.\alpha_{\v}.y_2\in A$,
\begin{equation}\label{eq:mostWareSame}
    \alpha_{\alpha_{\v}.y_2}^{\mathcal{P}}(\tau(\alpha_{\v}.y_2,\w))=\alpha_{\alpha_{\v'}.y_2}^{\mathcal{P}}(\tau(\alpha_{\v'}.y_2,\w)).
\end{equation}
The set of such valid $\w$ is large. Since $\alpha^{\tau}$ is an $\eta/N$-tempered time change for $\alpha$, 
\begin{equation*}
    \begin{aligned}
        m(\rho(\alpha_{\v}.y_2,C_{N})&\cap\rho(\alpha_{\v'}.y_2,C_{N})\cap C_{N})\\ \geq& (1-2\eta)^nm(C_{N})\geq(1-r/100)m(C_{N}).
    \end{aligned}
\end{equation*}
The additional parameters excluded by the condition of belonging to $A$ have relative measure at most $(1+O_n(\eta/N))r/400$, by stability and the Jacobian bound. Therefore \eqref{eq:mostWareSame} implies:
\[
d_{C_{N}}^{\alpha}(\alpha^{\mathcal{P}}_{\alpha_{\v}.y_2,\rho,N},\alpha^{\mathcal{P}}_{\alpha_{\v'}.y_2,\rho,N})<r/10.
\]
Combining this with \eqref{eq:y1y2CloseDifferentTime} and observing that $\alpha_{\v'}.y_2=\alpha^{\tau}_{\rho(y_1,\v)}.y_2$, we obtain for all $\v\in W$:
\begin{equation}\label{eq:CloseTilkdePhi}
    d_{C_{N}}^{\alpha}\left(\alpha^{\mathcal{P}}_{\alpha^{\tau}_{\rho(y_1,\v)}.y_1,\rho,N},\alpha^{\mathcal{P}}_{\alpha^{\tau}_{\rho(y_1,\v)}.y_2,\rho,N}\right)<r/5.
\end{equation}
This estimate shows that the reparametrized names match well on the block $C_N$ anchored at time $\rho(y_1,\v)$ in the \textbf{reparametrized $\R^n$-action $\alpha^{\tau}$}.

Finally, we aggregate these local estimates to the scale $M$. By the hypothesis on $W$, the number of ``bad'' lattice points is small:
\[
\card{C_M\cap \Gamma_N}-\card{W}\leq (\eta/5)\card{C_M\cap\Gamma_N}.
\]
We decompose the global distance on $C_M$ into the sum over ``good'' lattice points ($\v\in W$) where the error is $<r/5$, and ``bad'' lattice points where the error is at most $1$. Use the disjoint sets 
\[
E_{\v}=\rho(y_1,\v+C_N)=\rho(y_1,\v)+\rho(\alpha_{\v}.y_1,C_N).
\]
Each differs from $\rho(y_1,\v)+C_N$ by relative measure $O(\eta/N)$. Including the boundary and truncation errors, for sufficiently large $M$ we obtain
\[
d_{C_M}^{\alpha}(\alpha_{y_1,\rho,M}^{\mathcal{P}},\alpha_{y_2,\rho,M}^{\mathcal{P}})<\eta/5+r/5+O_{n}(\eta/N+N/M)<r/2.
\]
This concludes the proof of Lemma~\ref{lemma:discreteCounting}.
\end{proof}

\subsubsection{A combinatorial counting lemma}\label{sec:combinatorialLemma}
By assuming additional closeness on a lattice of $\R^n$ for the  reparametrization cocycles, Lemma~\ref{lemma:discreteCounting} establishes that pairs Hamming-close under the action $\alpha$ are also Hamming-close under $\alpha^{\tau}$. We now investigate the additional error terms in the counting of Hamming balls for $\alpha^{\tau}$ arising from these closeness assumptions. This, together with Lemma~\ref{lemma:discreteCounting}, provides a comprehensive preparation for the proof of Theorem~\ref{thm:AbramovMainNadler}. 

This section corresponds to \cite{NadlerAbramovNflow78}*{Lemma~5} and the proof is a modification of the original argument. 

\medskip

We first introduce some necessary definitions. In the following items, $N,M\in\N$, and recall that $\Gamma_N=2N\Z^n$.
\begin{enumerate}[label=(\roman*)]
    \item Given $\delta>0$, we define the fine lattice:
\begin{equation*}
     \Gamma(\delta)=2\delta\Z^n.
\end{equation*}
The \textbf{$\epsilon$-neighborhood} of a subset $Y\subset\Gamma(\delta)$ is defined as
\[
C_{\epsilon}(Y)=\{\v\in\R^n:\exists\gamma\in Y\text{ s.t. }\v\in C_{\epsilon}(\gamma)\}.
\]
    \item  For $\eta\in(0,1)$, a subset $W\subset C_M\cap\Gamma_N$ is \textbf{$\eta$-good} if 
    \[
    \card{W}\geq (1-\eta)\card{C_M\cap \Gamma_N}.
    \]
    \item Let $\alpha^{\tau}$ be a tempered time change, $\delta>0$, and $W\subset C_M\cap \Gamma_N$. Given $x\in X_0$, a map $Y:W\to \Gamma(\delta)$ has \textbf{configuration $W$} with respect to $x$, denoted $\Conf(Y)=(W,\rho,x)$, if 
    \[
    \rho(x,\v)\in C_{\delta}(Y(\v)),\qquad \v\in W.
    \]
    \item Let $\mathcal{Y}_{\delta}(W,\rho,x)$ be defined as
    \[
    \mathcal{Y}_{\delta}(W,\rho,x)=\{Y:W\to\Gamma(\delta):~\Conf(Y)=(W,\rho,x)\},
    \]
    and let $ \mathcal{Y}_{\delta}(W,\rho)=\bigcup_{x\in X_0} \mathcal{Y}_{\delta}(W,\rho,x)$.
\end{enumerate}

\begin{lemma}\label{lemma:latticeConfiguration}
Let $\epsilon,\delta\in(0,1)$. There exists constants $N'\in\N$ and $s'\in(0,1)$ such that for any integer $N\geq N'$ and any $s\in(0,s')$, there exists a constant $M'\in\N$ such that for all integer $M\geq M'$, the following holds: if $\alpha^{\tau}$ is a $2$-tempered time change and $W\subset C_M\cap\Gamma_N$ is an $s$-good subset, then
\[
\card{\mathcal{Y}_{\delta}(W,\rho)}\leq \exp(\epsilon\cdot m(C_M)/2).
\]
\end{lemma}
\begin{proof}[Proof of Lemma~\ref{lemma:latticeConfiguration}]
Fix $\epsilon,\delta\in(0,1)$. We choose $N'$ large enough such that for every $N\geq N'$, $N$ dominates terms of order $\log(\delta^{-1})$ and $\log N$, specifically:
\[
N^n>10^3n\epsilon^{-1}\log(\delta^{-1}N).
\]
We choose $s'$ small enough such that for every $s\in(0,s')$,
\[
10^4\delta^{-n}s<\epsilon/2.
\]
Fix an integer $N\geq N'$ and $s\in(0,s')$, we choose $M$ sufficiently large to mitigate boundary effects:
     \[
      m(C_{M-N})/m(C_{M})\geq 1-s.       
     \]

Fix $M$ large. We view $V=\Gamma_N\cap C_M$ as a graph with edges connecting adjacent points along the coordinate axes. Since $W$ is $s$-good, we can select a root $\w_1\in W\cap C_{M\sqrt[n]{2s}}$. Traverse a spanning tree of $V$ from $\w_1$, using every edge twice, and connect successive first visits to vertices of $W$. This gives a tree $\mathcal{T}$ on $W$ with total edge length at most $4N\card{V}$.

We estimate the number of valid configurations by determining the position of the image and then iteratively determining the images of descendants.

\textbf{1. Root Choice.} The root $\w_1$ is close to the origin. Since $J(\rho)=\id$ and $\alpha^{\tau}$ is $2$-tempered, $\rho(x,\cdot)$ is $3$-Lipschitz and fixes $\0$. The number of possible lattice point $\gamma\in\Gamma(\delta)$ satisfies
\[
\rho(x,\w_1)\in C_{\delta}(\gamma)
\]
is uniformly bounded in $x$ by 
\begin{equation*}%\label{eq:rho_kw_1}
   b_{\operatorname{root}} \leq C_n\delta^{-n}s\cdot m(C_M),
\end{equation*}
where $C_n$ depends only on $n$.

\textbf{2. Inductive Steps.} Let $\w,\w'\in W$ be the adjacent vertices in the tree $\mathcal{T}$ (i.e. $\w'$ is a child of $\w$). It follows that there is a path of at most $n$ coordinate parallel segments connecting $\w$ to $\w'$, with vertices in $\Gamma_N$. Define this path as $\{q_0,\ldots,q_{n}\}$ with $q_0=\w$ and $q_n=\w'$. We obtain from the Cauchy-Schwarz inequality that,
\[
\sum_{i=0}^{n-1}\norm{q_{i+1}-q_i}\leq\sqrt{n}\sqrt{\sum_{i=0}^{n-1}\norm{q_{i+1}-q_i}^2}=\sqrt{n}\norm{\w-\w'}.
\]

Notice that the map $\rho(x,\cdot)$ is $3$-Lipschitz, we have
\[
\norm{\rho(x,\w')-\rho(x,\w)}\leq 3\sqrt{n}\norm{\w-\w'}.
\]
Given the recorded value of $Y(\w)$, then the number of possible elements $\gamma\in\Gamma(\delta)$ such that 
\[
\rho(x,\w')\in C_{\delta}(\gamma)
\]
is bounded by the number of $\delta$-boxes that covers a ball of this radius:
\begin{equation*}%\label{eq:jStepBound}
   b(\w,\w')\leq C_n\norm{\w-\w'}^n \delta^{-n}.
\end{equation*}

\textbf{3. Total Count.} Notice that the total number of configurations is bounded by the product over all edges in the spanning tree is the set of edges of $\mathcal{T}$:
\begin{equation*}
    \begin{aligned}
        \card{\mathcal{Y}_{\delta}(W,\rho)}\leq b_{\operatorname{root}}\cdot\prod_{(\w,\w')\in E(\mathcal{T})}b(\w,\w').
    \end{aligned}
\end{equation*}
By taking the logarithm and using $\card{W}\leq\card{V}\leq 2m(C_M)/m(C_N)$ for large $M$:
\begin{equation}\label{eq:ArhoWest1}
    \begin{aligned}
        &\log\left(\card{\mathcal{Y}_{\delta}(W,\rho)}\right)\\
        \leq&\log b_{\operatorname{root}}+\sum_{(\w,\w')\in E(\mathcal{T})}\log b(\w,\w')\\
        \leq& \log(C_n\delta^{-n}s\cdot m(C_M))+\sum_{(\w,\w')\in E(\mathcal{T})}\left(\log(C_n\delta^{-n})+n\log\norm{\w-\w'}\right).
    \end{aligned}
\end{equation}
Combining the total edge length bound with Jensen's inequality, we obtain that for sufficiently large $M$,
\begin{equation*}
    \begin{aligned}
       \sum_{(\w,\w')\in E(\mathcal{T})}\log \norm{\w-\w'}\leq&(\card{W}-1)\log\left(\frac{4N\card{V}}{\card{W}-1}\right)\\
       \leq &\card{V}\log(16N)
    \end{aligned}
\end{equation*}
where the last step follows from $\card{W}\geq \card{V}/2$. Substituting this estimate back to \eqref{eq:ArhoWest1}, we obtain:
\begin{equation*}
    \begin{aligned}
        &\log\left(\card{\mathcal{Y}_{\delta}(W,\rho)}\right)\\
        \leq &\log(C_n\delta^{-n}s\cdot m(C_M))+\frac{2m(C_M)}{m(C_{N})}\left(\log(C_n\delta^{-n})+n\log(16N)\right).
    \end{aligned}
\end{equation*}
Notice that for $M$ sufficiently large, the first term divided by $m(C_M)$ goes to zero. The remaining normalized terms are $O_n(N^{-n}\log(2\delta^{-1}N))$ and are less than $\epsilon/4$ after enlarging $N'$ if necessary. Altogether, the total log-cardinality is bounded by $\epsilon\cdot m(C_M)/2$ and hence we finish the proof of Lemma~\ref{lemma:latticeConfiguration}.
\end{proof}

\subsection{The main proposition}\label{sec:AppMainProp}
The following proposition is the most important step for proving Abramov's entropy formula for time changes of $\R^n$-action. It establishes a comparison between the covering numbers for the original action $\alpha$ and the tempered time change $\alpha^{\tau}$. This result follows from the combination of the stability of partitions (Lemma~\ref{lem:stablePertQkphi}), name comparisons (Lemma~\ref{lemma:discreteCounting}), and the lattice counting bounds (Lemma~\ref{lemma:latticeConfiguration}). 

This proposition corresponds to \cite{NadlerAbramovNflow78}*{Proposition~4}. 

\medskip

We adopt the following setup:
\begin{itemize}
    \item $\alpha$ is a measurable, measure-preserving, ergodic, and free $\R^n$-action on the non-atomic standard probability space $(X,\mathcal{B},\mu)$.
    \item $\alpha^{\tau}$ is a $10^{-3}$-tempered time change of $\alpha$.
    \item $\{\mathcal{P}_k\}_{k\in\N}$ is a generating sequence of finite measurable partitions for $\alpha^{\tau}$, obtained via Proposition~\ref{prop:goodPartition}.
\end{itemize}

\begin{proposition}\label{prop:NadlerMainProp}
    Given sequences $\{\epsilon_k\}_{k\in\N}\subset(0,1)$ and $\{r_k\}_{k\in\N}\subset(0,1)$, 
    there exist a sequence $\{s_k\}_{k\in\N}$ with $s_k\in(0,r_k)$ and a sequence of finite measurable partitions $\{\mathcal{Q}_k\}_{k\in\N}$ such that for every $k\in\N$, the following holds for all sufficiently large $M$: if $\mathcal{C}_k$ is a $(\mathcal{Q}_k,C_M,s_k)$-$(1-\epsilon_k)$ covering family for $\alpha$, then there exists a $(\mathcal{P}_k,C_M,r_k)$-$(1-6\epsilon_k)$ covering family $\mathcal{D}_k$ for $\alpha^{\tau}$ such that 
    \begin{equation*}
        \begin{aligned}
        \operatorname{Card}(\mathcal{D}_k)\leq e^{\epsilon_k m(C_M)}\operatorname{Card}(\mathcal{C}_k).
        \end{aligned}
    \end{equation*}
\end{proposition}
\subsubsection{Specification of Constants and Sets}
    Fix $k\in\N$. We specify the constant and set for the proof in the following logical order:
    \begin{enumerate}[label=(\roman*)]
        \item\label{item:choicerk} \textbf{Stability Precision ($r'_k$):} Let $r'_k\in(0,10^{-9}r_k)$ be small enough so that
        \begin{equation*}
            \begin{gathered}
            1-(1-r'_k)^n<10^{-9}r_k,\\
    -6r'_k\log(2r'_k)-(1-2r'_k)\log(1-4r'_k)<\epsilon_k/2.
            \end{gathered}
        \end{equation*}
        
        \item\label{item:mostPointsAreStable} \textbf{Partition Stability ($T_k,\delta_k,A_k$):} Apply  Proposition~\ref{prop:goodPartition} to $\alpha^{\tau}$ with $\eta=r'_k$. Let $T_k$, $\delta_k$ and $A_k$ be the resulting parameters. Thus for every $N\geq T_k$, $\mathcal{P}_k$ is an $r'_k$-stable partition for $\alpha^{\tau}$ with parameters $(A_k,\delta_k,N)$.
    \[
    \eta_k=10^{-3000}n^{-1}\big(\min(\delta_k,r'_k)\big)^{100}.
    \]
    \item \label{item:Nk}\textbf{Lattice Constants ($N_k,s_k,M'_k$):} Apply Lemma~\ref{lemma:latticeConfiguration} with $\epsilon=\epsilon_k$ and $\delta=\delta_k/(64\sqrt{n})$. Let $N'_k$ and $s'_k$ be the constants provided by the lemma. We set our working parameters as:
    \[
    N_k=2\max(N'_k,T_k),\qquad s_k=\min(\eta_k^2/400,(s'_k)^2/16).
    \]
    Now let $M'_k$ be the threshold constant provided by Lemma~\ref{lemma:latticeConfiguration} corresponding to  $N_k$ and $2\sqrt{s_k}$.
    
    \item\label{eq:JrhokStable}\textbf{Perturbed Action ($\tau_k,\psi_k,\mathcal{Q}_k$):} By the equivalence theorem (cf.~Theorem~\ref{thm:kakutani-equivalences}), there exists an $\eta_k/N_k$-tempered time change $\alpha^{\tau_k}$ of $\alpha$ which is cohomologous to $\alpha^{\tau}$ via an isomorphism $\psi_k$. Specifically:
\begin{equation}\label{eq:rhokPerturb}
    \begin{gathered}
        \norm{\partial_{\v}\rho_{k}(x,\0)-\id}<\eta_k/N_k \text{ for $\mu$-a.e. $x\in X$,}\\
        \alpha^{\tau_k}_{\v}.\psi_{k}(x)=\psi_{k}(\alpha^{\tau}_{\v}.x)\text{ for all $\v\in\R^n$ and $\mu$-a.e. $x\in X$.}
    \end{gathered}
\end{equation}
Define the image partition $\mathcal{Q}_k=\psi_k(\mathcal{P}_k)$.

    \item\label{item:GoodAprSetE}\textbf{Good Set $G(R)$:} Recall that $\rho$ and $\rho_k$ are the reparametrization cocycles for $\alpha^{\tau}$ and $\alpha^{\tau_k}$ respectively. Define the set of points with good ergodic approximation: 
    \begin{equation*}
        \begin{gathered}
            G(R)=\{x\in X:\norm{\rho(x,\v)-\v}+\norm{\rho_k(x,\v)-\v}<r'_k\norm{\v}\text{ if $\norm{\v}\geq R$}\}.
        \end{gathered}
    \end{equation*}
    By Proposition~\ref{prop:jRhoProp}~\ref{item:JrhoErgodicSum}, there exists $M''_k\in\N$ such that for all $M\geq M''_k$
    \[
    \mu(G(M))>1-\epsilon_k.
    \]
    \item\label{item:choiceN}  We consider $M\in\N$ sufficiently large:
    \[
    M\geq \max(M'_k,M''_k).
    \]
    \end{enumerate}

\subsubsection{Covering Comparison Lemma}

The core of the proposition is the following lemma, which bounds the covering number of $\alpha$ by $\alpha^{\tau_k}$.
\begin{lemma}\label{lem:HammingTimeChangeCounting}
    For every $x\in X$, there exist an integer $p\leq\exp(\epsilon_km(C_M))$ and points $y_1,\ldots,y_p\in X_0$ such that
    \[
    B^{H,\alpha}_{s_k}(x,C_{M},\mathcal{Q}_k)\cap X_0\subset \bigcup_{j=1}^pB^{H,\alpha^{\tau_k}}_{r_k/2}(y_j,C_{M},\mathcal{Q}_k).
    \]
\end{lemma}
\begin{proof}[Proof of Lemma~\ref{lem:HammingTimeChangeCounting}]
Note that $\mathcal{P}_k$ is an $r'_k$-stable partition for $\alpha^{\tau}$ with parameters $(A_k,\delta_k,N_k)$, and $\alpha^{\tau_k}$ is an $\eta_k/N_k$-tempered time change of $\alpha$ that is cohomologous to $\alpha^{\tau}$ through $\psi_k$. These, together with the choice of $\eta_k$ (cf.~\ref{item:mostPointsAreStable} on \pageref{item:mostPointsAreStable}) and Lemma~\ref{lem:stablePertQkphi}, give that $\mathcal{Q}_k$, i.e. $\psi_k(\mathcal{P}_k)$, is a $2r'_k$-stable partition for $\alpha$ with parameters $(\psi_k(A_k),\delta_k/2,N_k)$.

Fix $x\in X$. Suppose $y\in B^{H,\alpha}_{s_k}(x,C_{M},\mathcal{Q}_k)\cap X_0$. This implies their $\alpha$-names are $s_k$-close:
\begin{equation*}%\label{eq:xyHammingClose}
    d^{\alpha}_{C_{M}}(\alpha_{x}^{\mathcal{Q}_k},\alpha_{y}^{\mathcal{Q}_k})<s_k.
\end{equation*}
Consider the coarse lattice $\Gamma_{N_k}=2N_k\Z^n$. We define the set of grid points where the local names match:
\begin{equation*}%\label{eq:W1Def}
 W_{N_k,M}^{\mathcal{Q}_k}(x,y,s_k)=\left\{\v\in C_{M}\cap\Gamma_{N_k}:d_{C_{N_k}}^{\alpha}(\alpha_{\alpha_{\v}.x}^{\mathcal{Q}_k},\alpha_{\alpha_{\v}.y}^{\mathcal{Q}_k})<\sqrt{s_k}\right\}.
 \end{equation*}
By the Markov inequality, for sufficiently large $M$ and every such $y$, the set of matching grid points is very dense:
\begin{equation*}%\label{eq:CardWBound}
  \card{W_{N_k,M}^{\mathcal{Q}_k}(x,y,s_k)}\geq (1-2\sqrt{s_k})\card{C_M\cap\Gamma_{N_k}}.  
\end{equation*}
Let $\mathcal{W}$ be the collection of all such subsets $W$ arising from $y\in B^{H,\alpha}_{s_k}(x,C_{M},\mathcal{Q}_k)\cap X_0$:
\[
\mathcal{W}=\{W_{N_k,M}^{\mathcal{Q}_k}(x,y,s_k):y\in B^{H,\alpha}_{s_k}(x,C_{M},\mathcal{Q}_k)\cap X_0\}.
\]
The cardinality of $\mathcal{W}$ is bounded by the number of ways to choose the ``bad'' grid points. Using Stirling's approximation and Condition~\ref{item:choicerk} on $r'_k$ (and thus $s_k$):
\begin{equation}\label{eq:HammingCloseComBound}
\begin{aligned}
     &\card{\mathcal{W}}\leq \sum_{j=0}^{\lfloor2\sqrt{s_k}\card{C_M\cap\Gamma_{N_k}}\rfloor}\left(\begin{array}{cc}
         \card{C_{M}\cap\Gamma_{N_k}}  \\
         j
    \end{array}\right) \\ \leq &\exp\Big(m(C_M)\big(-6\sqrt{s_k}\log(2\sqrt{s_k})-(1-2\sqrt{s_k})\log(1-4\sqrt{s_k})\big)\Big)\\
     \leq &\exp\Big(\epsilon_km(C_M)/2\Big).
\end{aligned}
\end{equation}

\medskip     

Now, fix a specific profile $W\in\mathcal{W}$. We count how many distinct reparametrizations are compatible with this profile. Assign to each $y$ its configuration on $W$ by a fixed choice of nearest point in $\Gamma(\delta_k/(64\sqrt{n}))$. If $y_1$, $y_2$ have the same profile and configuration, then
\[
\begin{gathered}
    d_{C_M}^{\alpha}(\alpha_{y_1}^{\mathcal{Q}_k},\alpha_{y_2}^{\mathcal{Q}_k})<2s_k<\eta_k^2/100,\\
    d_{C_{N_k}}^{\alpha}(\alpha_{\alpha_{\v}.y_1}^{\mathcal{Q}_k},\alpha_{\alpha_{\v}.y_2}^{\mathcal{Q}_k})<2\sqrt{s_k}<\eta_k/10,\qquad \v\in W,\\
    \norm{\rho_k(y_1,\v)-\rho_k(y_2,\v)}\leq \delta_k/16,\qquad \v\in W.
\end{gathered}
\]
Since $\card{(C_M\cap\Gamma_{N_k})\setminus W}<(\eta_k/5)\card{C_M\cap\Gamma_{N_k}}$, Lemma~\ref{lemma:discreteCounting} with $r=r_k/2$ and $\delta=\delta_k/2$, gives
\[
 d_{C_M}^{\alpha}(\alpha_{y_1,\rho_k,M}^{\mathcal{Q}_k},\alpha_{y_2,\rho_k,M}^{\mathcal{Q}_k})<r_k/4.
\]
The two truncation regions have total relative measure $O_n(\eta_k/N_k)<r_k/4$, so the actual $\alpha^{\tau_k}$-Hamming distance is less than $r_k/2$. Recall from \S\ref{sec:combinatorialLemma} that $\mathcal{Y}_{\delta_k/(64\sqrt{n})}(W,\rho_k)$ records the configurations as $y$ varies. By Lemma~\ref{lemma:latticeConfiguration}, and the choices of $N_k$, $s_k$ and $M$,
\[
\card{\mathcal{Y}_{\delta_k/(64\sqrt{n})}(W,\rho_k)}\leq \exp(\epsilon_k\cdot m(C_M)/2).
\]
Choosing one point in each nonempty class gives the required balls. Combining the number of profiles (cf.~\eqref{eq:HammingCloseComBound}) and the number of configurations (cf.~above estimates), the total number of $\alpha^{\tau_k}$-balls required to cover the $\alpha$-ball is bounded by the product:
\[
p\leq \card{\mathcal{W}}\cdot \max_{W\in\mathcal{W}}\card{\mathcal{Y}_{\delta_k/(64\sqrt{n})}(W,\rho_k)}\leq \exp(\epsilon_k\cdot m(C_M)).
\]
This concludes the proof of Lemma~\ref{lem:HammingTimeChangeCounting}.
\end{proof}

\subsubsection{Proof of Proposition~\ref{prop:NadlerMainProp}}
Let $\mathcal{C}_k$ be a $(\mathcal{Q}_k,C_M,s_k)$-$(1-\epsilon_k)$ covering family for $\alpha$. By the construction in Lemma~\ref{lem:HammingTimeChangeCounting}, for each $C\in\mathcal{C}_k$, we can cover $C\cap G(M)\cap X_0$ with at most $\exp(\epsilon_k\cdot m(C_M))$ balls of the form $B^{H,\alpha^{\tau_k}}_{r_k}(y,C_{M},\mathcal{Q}_k)$, choosing their centers in $C\cap G(M)\cap X_0$. Note that we use radius $r_k$ (instead of $r_k/2$) here to absorb the boundary errors arising from reparametrization.

Indeed, for $x,y\in G(M)$ (the good set from~\ref{item:GoodAprSetE} on p.~\pageref{item:GoodAprSetE}), the choice of $M$ (cf.~\ref{item:GoodAprSetE} and \ref{item:choiceN}) implies that the regions where their $\rho_k$ and $\rho$-reparametrization map into $C_M$ overlap significantly with $C_M$. Specifically, the choice of $M$ ensures:
\[
m\big(\rho_k(x,C_M)\cap\rho(x,C_M)\cap\rho_k(y,C_M)\cap\rho(y,C_M)\cap C_M\big)\geq (1-r_k/2)m(C_M).
\]
Consequently, if their reparametrized names are close:
\[
d_{C_M}^{\alpha}(\alpha^{\mathcal{Q}_k}_{x,\rho_k,M},\alpha^{\mathcal{Q}_k}_{y,\rho_k,M})<r_k/2,
\]
then the actual $\alpha^{\tau_k}$-Hamming distance on $C_M$ satisfies:
\[
    d_{C_{M}}^{\alpha^{\tau_k}}\big(\alpha^{\mathcal{Q}_k}_{x},\alpha^{\mathcal{Q}_k}_{y}\big)<r_k.
\]

Collecting these balls for all $C\in\mathcal{C}_k$, we obtain a family $\mathcal{D}'_k$ for $\alpha^{\tau_k}$ with cardinality:
\[
\card{\mathcal{D}'_k}\leq e^{\epsilon_km(C_M)} \card{\mathcal{C}_k}.
\]
We now transfer this cover to $\alpha^{\tau}$ using the isomorphism $\psi_k$. Since $\alpha^{\tau_k}_{\v}.\psi_{k}(x)=\psi_{k}(\alpha^{\tau}_{\v}.x)$ for all $\v\in\R^n$ and $\mu$-a.e. $x\in X$, the map $\psi_k^{-1}$ sends the orbits of $\alpha^{\tau_k}$ to orbits of $\alpha^{\tau}$. Define the family $\mathcal{D}_k$ as
\[
\mathcal{D}_k=\{\psi_k^{-1}(B):B\in\mathcal{D}'_k\}.
\]
Note that since $\mathcal{Q}_k=\psi_k(\mathcal{P}_k)$, we have
\[
\psi_k^{-1}\big(B^{H,\alpha^{\tau_k}}_{r_k}(y,C_{M},\mathcal{Q}_k)\big)=B^{H,\alpha^{\tau}}_{r_k}(\psi_k^{-1}(y),C_{M},\mathcal{P}_k).
\]

Finally, we verify the measure. The union of the original balls covers a set of $\mu$-measure at least $1-\epsilon_k$. Restricting to $G(M)\cap X_0$ loses at most another $\epsilon_k$. Since $J(\rho_k)=\id$, we have
\[
\frac{d\mu^{\tau_k}}{d\mu}=\abs{\det D_{\v}\rho_k(x,\0)}\leq (1+\eta_k/N_k)^n<2.
\]
Using $(\psi_k)_*\mu^{\tau}=\mu^{\tau_k}$, we obtain
\[
\mu^{\tau}\big(X\setminus\cup\mathcal{D}_k\big)=\mu^{\tau_k}\big(X\setminus\cup\mathcal{D}'_k\big)\leq 2\mu\big(X\setminus(\cup\mathcal{C}_k\cap G(M)\cap X_0)\big)\leq 6\epsilon_k.
\]
Thus, $\mathcal{D}_k$ is the required cover for $\alpha^{\tau}$ and we complete the proof of Proposition~\ref{prop:NadlerMainProp}\hfill $\Box$

\subsection{Proof of Theorem~\ref{thm:AbramovMainNadler}}\label{sec:AppAbramovMainNadler}
 By Proposition~\ref{prop:barelyLinTempControl} and invariance of entropy under measurable isomorphism, we may assume that $\alpha^{\tau}$ is $10^{-3}$-tempered. Let $\{\mathcal{P}_k\}_{k\in\N}$ be an increasing generating sequence of finite measurable partitions for $\alpha^{\tau}$, obtained from Proposition~\ref{prop:goodPartition} by passing to the subsequence indexed by $k!$ and relabeling.

 The inequality $h_{\mu^{\tau}}(\alpha^{\tau})\leq h_{\mu}(\alpha)$ if $h_{\mu}(\alpha)=\infty$. Otherwise, fix $\delta>0$ and a finite number $H<h_{\mu^{\tau}}(\alpha^{\tau})$. Choose $\{r_k\}_{k\in\N}\subset(0,1)$ such that
 \begin{equation}\label{eq:hPDrchoice}
      r_k\to 0,\qquad h_{\mu^{\tau}}(\phi^{\tau},\mathcal{P}_k,2r_k)>H.
 \end{equation}
for all sufficiently large $k$. Let $\{\epsilon_k\}_{k\in\N}\subset(0,\min(\delta,10^{-2}))$ tend to zero. 

 We invoke Proposition~\ref{prop:NadlerMainProp}. Fixing $\{r_k\}_{k\in\N}$, $\{\epsilon_k\}_{k\in\N}$, and $\{\mathcal{P}_k\}_{k\in\N}$, the proposition provides a sequence $\{s_k\}_{k\in\N}\subset(0,1)$ and a sequence of finite measurable partitions $\{\mathcal{Q}_k\}_{k\in\N}$ for $\alpha$.

For sufficiently large $M$, choose  a $(\mathcal{Q}_k,C_M,s_k)$-$(1-\epsilon_k)$ covering family $\mathcal{C}_k$ consisting of $\mathcal{Q}_k$-$s_k/2$-$C_M$-balls and satisfying
\[
\log\card{\mathcal{C}_k}\leq m(C_M)(h_{\mu}(\alpha,\mathcal{Q}_k,s_k/2)+\delta).
\]
Such a family exists by Theorem~\ref{thm:GeneralShannonMcMillan}. Proposition~\ref{prop:NadlerMainProp} implies the existence of a $(\mathcal{P}_k,C_M,r_k)$-$(1-6\epsilon_k)$ covering family $\mathcal{D}_k$ for $\alpha^{\tau}$ such that 
    \begin{equation}\label{eq:inequalityCandD}
        \begin{aligned}
        \operatorname{Card}(\mathcal{D}_k)\leq e^{\epsilon_k m(C_M)}\operatorname{Card}(\mathcal{C}_k).
        \end{aligned}
    \end{equation}

\medskip

Every member of $\mathcal{D}_k$ has Hamming diameter less than $2r_k$. Since their union has measure greater than $1-6\epsilon_k$, Theorem~\ref{thm:GeneralShannonMcMillan} gives, for sufficiently large $M$,
\[
h_{\mu^{\tau}}(\phi^{\tau},\mathcal{P}_k,2r_k)-2\delta\leq\frac{\log\card{\mathcal{D}_k}}{m(C_M)}.
\]
Combining this with \eqref{eq:hPDrchoice} and \eqref{eq:inequalityCandD}, we obtain:
\begin{equation}\label{eq:hCNinequality}
    m(C_M)\cdot (H-2\delta)\leq \epsilon_km(C_M)+\log\card{\mathcal{C}_k}.
\end{equation}
By our choice of $\mathcal{C}_k$ and the definition of entropy, 
\[
\log\card{\mathcal{C}_k}\leq m(C_M)(h_{\mu}(\alpha,\mathcal{Q}_k,s_k/2)+\delta)\leq m(C_M)(h_{\mu}(\alpha)+2\delta).
\]
Substituting this back into \eqref{eq:hCNinequality} and dividing by $m(C_M)$:
\begin{equation*}
    H-2\delta\leq h_{\mu}(\alpha)+2\delta+\epsilon_k.
\end{equation*}
Letting $k\to\infty$, then $\delta\to0$ and finally  
$H\nearrow h_{\mu^{\tau}}(\alpha^{\tau})$, we obtain the upper bound:
\[
h_{\mu^{\tau}}(\alpha^{\tau})\leq h_{\mu}(\alpha).
\]

To prove the reverse inequality, we observe that the argument is symmetric. The action $\alpha$ can be viewed as a time change of $\alpha^{\tau}$ via the inverse cocycle. Since when $J(\rho)=\id$, the same argument, again using a cohomologous $10^{-3}$-tempered representative, gives $h_{\mu^{\tau}}(\alpha^{\tau})\geq h_{\mu}(\alpha)$.

Thus, we conclude
\[
h_{\mu}(\alpha)=h_{\mu^{\tau}}(\alpha^{\tau}).
\]
This finishes the proof of Theorem~\ref{thm:AbramovMainNadler} in the case of $J(\rho)=\id$. Combining this with the linear time change reduction (cf.~Proposition~\ref{prop:linearTimeChange}) completes the proof of Theorem~\ref{thm:AbramovMainNadler} in full generality.

\bibliographystyle{plain} 
\bibliography{refs}

\end{document}